\documentclass{amsart}

\usepackage{amsmath,amssymb,amsthm,amsbsy}
\usepackage{epsfig}

\newtheorem{thm}{Theorem}[section]
\newtheorem{lem}[thm]{Lemma}

\theoremstyle{definition}

\theoremstyle{remark}

\numberwithin{equation}{section}

\newcommand{\C}{\mathbb C}
\newcommand{\R}{\mathbb R}
\newcommand{\E}{\mathbb E}

\newcommand{\wt}{\widetilde}

\newcommand{\m}[1]{m^{(#1)}(z)}
\newcommand{\mpl}[1]{m_+^{(#1)}(z)}
\newcommand{\mni}[1]{m_-^{(#1)}(z)}

\begin{document}

\title[Longest Increasing Subsequence]
{On the Distribution of the Length of the Longest Increasing Subsequence 
of Random Permutations} 

\author{Jinho Baik}
\address{Courant Institute of Mathematical Sciences, New York}

\curraddr{}

\email{baik@cims.nyu.edu}
\thanks{The authors would like to acknowledge
many extremely useful and enlightening conversations with Persi Diaconis and
Andrew Odlyzko.  Special thanks are due 
to Andrew Odlyzko and Eric Rains for providing us with
the results of their Monte Carlo simulations. }


\author{Percy Deift}
\address{Courant Institute of Mathematical Sciences, New York}
\email{deift@cims.nyu.edu}
\thanks{The work of the second author was supported in part by NSF
grant \#DMS-9500867.}

\author{Kurt Johansson}
\address{Royal Institute of Technology, Sweden}
\email{kurtj@math.kth.se}
\thanks{The work of the third
author was supported in part by the Swedish Natural Research Council
(NFR)}

%    General info

\date{July 2, 1998 ; March 25, 1999 (revised)}

\dedicatory{}


\begin{abstract}
  The authors consider the length, $l_N$, of the length of the longest 
increasing subsequence of a random permutation of $N$ numbers. 
The main result in this paper is a proof that the distribution function 
for $l_N$, suitably centered and scaled, converges to the 
Tracy-Widom distribution [TW1] of the largest eigenvalue of a random 
GUE matrix. 
The authors also prove convergence of moments.
The proof is based on the steepest decent method for Riemann-Hilbert 
problems, introduced by Deift and Zhou in 1993 [DZ1] in the context 
of integrable systems.
The applicability of the Riemann-Hilbert technique depends, in turn, 
on the determinantal formula of Gessel [Ge] for the Poissonization 
of the distribution function of $l_N$.
\end{abstract}

\maketitle

%%  =============   my file below   =====

%%\tableofcontents
%%\newpage

\section{{\bf Introduction}}
\label{s-intro}

Let $S_N$ be the group of permutations of $1,2,\dots,N$.
If $\pi\in S_N$, we say that $\pi(i_1),\cdots,\pi(i_k)$
is an increasing subsequence in $\pi$ if $i_1<i_2< \cdots <i_k$ and
$\pi(i_1) < \pi(i_2) < \cdots < \pi(i_k)$.
Let $l_N(\pi)$ be the length of the longest increasing subsequence.
For example, if $N=5$ and $\pi$ is the permutation 
$5\ \ 1\ \ 3\ \  2\ \ 4$ (in one-line notation 
: thus $\pi(1)=5$, $\pi(2)=1$, $\dots$), then the longest
increasing subsequences are $1\ \ 2\ \ 4$ and $1\ \ 3\ \ 4$, and $l_N(\pi)=3$.
Equip $S_N$ with uniform distribution, 
\begin{equation*}
   q_{n,N}=Prob(l_N\le n)=\frac{f_{N,n}}{N!}, 
\end{equation*}
where $f_{N,n}=\# \text{(permutations $\pi$ in $S_N$ with $l_N\le n$)}$.
The goal of this paper is to determine the asymptotics of 
$q_{n,N}$ as $N\to\infty$.
This problem was raised by Ulam in the early 60's \cite{Ul}, 
and on the basis of Monte Carlo simulations, he 
conjectured that the limit
\begin{equation}
   c\equiv \lim_{N\to\infty} \frac1{\sqrt{N}} \E_N(l_N)
\end{equation}
exists.   
(Here $\E_N(\cdot)$ denotes the expectation value with respect
to the distribution function $q_{n,N}$.)
The problem of proving the existence of this limit and the computation 
of $c$ has became known as ``Ulam's problem''.
An argument of Erd\"os and Szekeres \cite{ES} shows that 
$\E_N(l_N)\geq \frac12\sqrt{N-1}$, so that if the limit 
exists, then $c\geq\frac12$.
Subsequent numerical work by Baer and Brock \cite{BB} 
in late 60's suggested that value of $c$ is 2.
The existence of the limit was rigorously 
established by Hammersley \cite{Ha} in 1972.
In \cite{LS}, Logan and Shepp proved that $c\geq 2$ and simultaneously 
Vershik and Kerov \cite{VK} (see also \cite{VK2}) showed that $c=2$, 
thus settling Ulam's problem.
Alternative proofs of Ulam's problem are due to Aldous and Diaconis 
\cite{AD}, Sepp\"al\"ainen \cite{Se} and Johansson \cite{Johan}.
Over the years, various conjectures have been made concerning 
the variance $Var(l_N)$ of $l_N$, and Monte Carlo simulations of 
Odlyzko and Rains beginning in 1993, indicated that 
\begin{equation}\label{e-c0}
   \lim_{N\to\infty} \frac1{N^{1/3}} Var(l_N) = c_0
\end{equation}
for some numerical constant $c_0 \sim 0.819$. 
Also Odlyzko and Rains computed $\E(l_N)$ to higher order 
and found 
\begin{equation}\label{e-c1}
   \lim_{N\to\infty} \frac{\E(l_N)-2\sqrt{N}}{N^{1/6}} =c_1
\end{equation}
where $c_1 \sim -1.758$.
Further historical information on Ulam's problem, together with some 
discussions of the methods used by various authors, can be found 
in \cite{AD} and \cite{OR}

\bigskip
Before stating our results, we need to define the Tracy-Widom 
distribution \cite{TW} (see below).
Let $u(x)$ be the solution of the Painlev\'e II (PII) equation,
\begin{equation}\label{e-int1}
  u_{xx}=2u^3+xu, \quad\text{and}\quad 
u\sim -Ai(x)\quad\text{as}\quad x\to\infty,
\end{equation}
where $Ai$ is the Airy function.
The (global) existence and uniqueness of this solution was first 
established in \cite{HM} : 
the asymptotics as $x\to\pm\infty$ are,  
\begin{equation}\label{e-asyPII}
  \begin{split}
   &u(x)=-Ai(x)+O\biggl(\frac{e^{-(4/3)x^{3/2}}}{x^{1/4}}\biggr) 
\quad \text{as}\quad x\to\infty,\\
   &u(x)= -\sqrt{\frac{-x}2} \biggl(1+O\bigl(\frac1{x^2}\bigr)\biggr) 
\qquad \text{as}
\quad x\to -\infty,
  \end{split}
\end{equation}
(see, for example, \cite{HM}, \cite{IN}, \cite{DZ}).
Recall \cite{AS} that $Ai(x)\sim 
\frac{e^{-(2/3)x^{3/2}}}{2\sqrt{\pi}x^{1/4}}$ as $x\to\infty$.
Define the Tracy-Widom distribution 
\begin{equation}\label{e-F}
  F(t) = \exp \biggl(-\int_t^\infty (x-t)u^2(x) dx\biggr).
\end{equation}
From~\eqref{e-asyPII} and~\eqref{e-F}, $F'(t)>0$, 
$F(t)\to 1$ as $t\to +\infty$ and $F(t)\to 0$ 
as $t\to -\infty$, so that $F$ is indeed a distribution function.
Our first result concerns the convergence of $l_N$ in distribution 
after appropriate centering and scaling. 
\begin{thm}\label{thm1}
  Let $S_N$ be the group of all permutations of $N$ numbers
with uniform distribution and let $l_N(\pi)$ be the length of
the longest increasing subsequence of $\pi\in S_N$.
   Let $\chi$ be a random variable whose distribution function is $F$.
   Then, as $N\to\infty$, 
\begin{equation*}
   \chi_N\equiv\frac{l_N-2\sqrt{N}}{N^{1/6}}\to\chi 
\qquad\text{in distribution,}
\end{equation*}
i.e.
\begin{equation*}
   \lim_{N\to\infty} Prob\biggl(\chi_N\equiv
\frac{l_N-2\sqrt{N}}{N^{1/6}}\le t\biggr)
=F(t) \quad\text{for all $t\in\R$.}
\end{equation*}
\end{thm}
In order to show that the moments of $\chi_N$ converge to the
corresponding moments of $\chi$ we need estimates for the distribution
function $F_N(t)$ of $\chi_N$ for large $|t|$. From the large
deviation formulas for $l_N$ (see below), we expect that $F_N(t)$
(resp., $1-F_N(t)$) 
should go to zero rapidly as $t\to-\infty$ (resp., $t\to+\infty$). 
In fact, we will prove that, for
$M>0$ sufficiently large, there are positive constants $c$ and $C(M)$
such that
\begin{equation}\label{e-J1.6a'}
F_N(t)\le C(M)e^{ct^3}     %%(1.6a')
\end{equation}
if $-2N^{1/3}\le t\le -M$, and
\begin{equation}\label{e-J1.6b'}
1-F_N(t)\le C(M)e^{-ct^{3/5}}   %% (1.6b')
\end{equation}
if $M\le t\le N^{5/6}-2N^{1/3}$.
Together with Theorem~\ref{thm1}  these estimates yield 
\begin{thm}\label{thm2}
   For any $m=1,2,3,\cdots$, we have
\begin{equation*}
   \lim_{N\to\infty} \E_N(\chi_N^m)= \E(\chi^m),
\end{equation*}
where $\E(\cdot)$ denotes expectation with respect to the distribution
function $F$. In particular,
 \begin{equation}\label{e-Jcor1}
      \lim_{N\to\infty} \frac{Var(l_N)}{N^{1/3}} =
\int_{-\infty}^{\infty} t^2dF(t) -
\biggl( \int_{-\infty}^{\infty} tdF(t) \biggr)^2.
  \end{equation}
and
\begin{equation}\label{e-Jcor2}
  \lim_{N\to\infty} \frac{\E_N(l_N)-2\sqrt{N}}{N^{1/6}} =
\int_{-\infty}^{\infty} tdF(t).
\end{equation}
\end{thm}
If one solves the Painlev\'e II equation~\eqref{e-int1} numerically 
(see, \cite{TW}), and then 
computes the integrals on the RHS of the formulae of ~\eqref{e-Jcor1} 
and~\eqref{e-Jcor2}, one obtains the values 0.8132 and -1.7711 
which agree with 
$c_0$ and $c_1$ in~\eqref{e-c0} and~\eqref{e-c1} respectively, 
up to two decimal places.

\bigskip
The distribution function $F(t)$ in Theorems \ref{thm1} and \ref{thm2} 
first arose in the work of Tracy and Widom on 
the Gaussian Unitary Ensemble (GUE) of random matrix theory.
In this theory (see, e.g., \cite{Me}), one considers $N\times N$ 
hermitian matrix $M=(M_{ij})$ with probability density
\begin{equation*}
   Z_N^{-1}e^{-tr(M^2)}dM
=Z_N^{-1}e^{-tr(M^2)} \biggl(\prod_{i=1}^{N} dM_{ii}\biggr) 
\prod_{i=1}^{N} d(ReM_{ij}) d(ImM_{ij}),
\end{equation*}
where $Z_N$ is the normalization constant.
In \cite{TW}, Tracy and Widom showed that 
as the size of the hermitian matrices increases, 
the distribution of the (properly centered and scaled) largest 
eigenvalue of a random GUE matrix 
converges precisely to $F(t)$ !
In other words, properly centered and scaled, the length of the longest 
increasing subsequence for a permutation $\pi\in S_N$, behaves statistically 
for large $N$ like the largest eigenvalue of a random GUE matrix 
(see the Appendix for an intuitive argument).
In \cite{TW}, the authors also computed the distribution functions of 
the second, third, $\cdots$ largest eigenvalues of such random matrices, 
and the question arises whether such distribution functions describe 
the statistics of quantities identifiable in the random permutation 
context.

Recall the Robinson-Schensted correspondence 
(see, e.g., \cite{Sa}, and also Section 5.1.4 in \cite{Kn}) 
which establishes a bijection $\pi\mapsto (P(\pi),Q(\pi))$ from $S_N$ 
to pairs of Young tableaux with shape$(P(\pi))$  $=$shape$(Q(\pi))$. 
Under this correspondence, the number of boxes in the first row 
of $P(\pi)$ (equivalently $Q(\pi)$) is precisely $l_N(\pi)$, 
(see, \cite{Sa}, \cite{Kn}).
In other words, the results on $l_N$ can be rephrased as results on  
the statistics of the number of boxes in the first row of Young tableaux.
Monte Carlo simulations of Odlyzko and Rains \cite{OR} indicate that $\wt{l}_N$, 
the number of boxes in the second row of $P(\pi)$ 
(equivalently $Q(\pi)$) behaves 
statistically for large $N$, like the second largest eigenvalue of a 
random GUE matrix.
More precisely, their simulations indicate that 
\begin{equation*}
   \lim_{N\to\infty} \frac{\E_N(\wt{l}_N)-2\sqrt{N}}{N^{1/6}}
=-3.618, 
\end{equation*}
and 
\begin{equation*}
   \lim_{N\to\infty} \frac{Var(\wt{l}_N)}{N^{1/3}} = 0.545.
\end{equation*}
These values agree, once again, to two decimal places with 
the mean and variance of the suitably centered 
and scaled second largest eigenvalue of a GUE matrix, 
as computed in \cite{TW}.
Presumably, the number of boxes in the third row of $P(\pi)$ should 
behave statistically like the third largest eigenvalue of a GUE matrix 
as $N\to\infty$, etc. 
In recent work \cite{BDJ}, the authors have shown that this 
conjecture is indeed true for the second row.
Also, beautiful results of Okounkov \cite{Ok}, 
using arguments from combinatorial topology, have now provided an elegant 
basis for understanding the relationship between the statistics 
of Young tableaux and the eigenvalues of random matrices.
Over the last year, many other intriguing results have been obtained 
on a variety of problems arising in mathematics and mathematical 
physics, which are closely related to, or motivated by, 
the longest increasing subsequence problem. 
We refer the reader to 
\cite{TW2}, \cite{Bo}, \cite{Jo2, Jo3} and \cite{BR}.

\bigskip
As in \cite{Johan}, we consider the Poissonization $\phi_n(\lambda)$ 
of $q_{n,N}$, 
\begin{equation}\label{e-int4}
  \phi_n(\lambda) \equiv \sum_{N=0}^{\infty} 
\frac{e^{-\lambda}\lambda^N}{N!} q_{n,N}.
\end{equation}
The function $\phi_n(\lambda)$ is a distribution function (in $n$) of
a random variable $L(\lambda)$ coming from a superadditive process
introduced by Hammersley in \cite{Ha}, and used by him to show that the
limit (1.1) exists. The random variable $L(\lambda)$ is defined as
follows.
Consider a homogeneous rate one Poisson process in
the plane and let $L(\lambda)$ denote the maximum number of points in
an up-right (increasing) path through the points starting at $(0,0)$
and ending at $(\sqrt{\lambda},\sqrt{\lambda})$. For more details see
\cite{AD} and \cite{Se2}, and for a generalization to the non-homogeneous case
see \cite{DeZe}. Theorem~\ref{thm1}  and~\ref{thm2} 
hold for the random variable
$L(\lambda)$ as $\lambda\to\infty$.
Referring to the ``de-Poissonization'' 
Lemmas~\ref{lem-dep2} and ~\ref{lem-dep3} 
below, we see that it is easy to recover the asymptotics of $q_{n,N}$ 
as $N\to\infty$ from the knowledge of $\phi_n(\lambda)$ for 
$\lambda\sim N$. 
In other words, in order to compute the asymptotics of $l_N$, 
we must investigate 
the double scaling limit of $\phi_n(\lambda)$ when 
$\lambda\to\infty$ and $1\le n\le N\sim\lambda$, and this is 
the technical thrust of the paper.

To this end we use the following representation for $\phi_n(\lambda)$, 
\begin{equation}\label{e-phi}
    \phi_n(\lambda) 
    = e^{-\lambda} D_{n-1}(\exp(2\sqrt{\lambda}\cos{\theta})), 
\end{equation}
where $D_{n-1}$ denotes the $n\times n$ Toeplitz determinant
with weight function $f(e^{i\theta})=\exp(2\sqrt{\lambda}\cos{\theta})$
on the unit circle, (see, e.g.,\cite{Sz}).
The above formula follows from work of Gessel in \cite{Ge} using well known 
results about Toeplitz determinants.
As noted in \cite{Johan}, the formula can also be proved using 
the following representation for $q_{n,N}$, $1\le n\le N$, 
discovered by \cite{OPWW}, 
\begin{equation}\label{e-rains}
   q_{n,N}=\frac{2^{2N}N!}{(2N)!}
\int_{[-\pi,\pi]^n} \bigl( \sum_{j=1}^{n} \cos\theta_j \bigr)^{2N}
\prod_{1\le j<k\le n} |e^{i\theta_j}-e^{i\theta_k}|^2
\frac{d^n\theta}{(2\pi)^nn!}.
\end{equation}
In addition, an earlier result of Diaconis and Shahshahani (\cite{DS}) shows 
that the above formula~\eqref{e-rains} is true also in the case 
$n>N$ when $q_{n,N}\equiv 1$. 
Inserting~\eqref{e-rains} into ~\eqref{e-int4}, we obtain 
\begin{equation}\label{e-DiS}
    \phi_n(\lambda) = e^{-\lambda} \frac{1}{(2\pi)^nn!} 
\int_{[-\pi,\pi]^n}\exp(2\sqrt{\lambda}\sum_{j=1}^{n} \cos{\theta_j}) 
\prod_{1\leq j<k \leq n} |e^{i\theta_j}-e^{i\theta_k}|^2 d^n\!\theta, 
\end{equation}
which is precisely~\eqref{e-phi} by standard methods in the theory of 
Toeplitz determinants (see, \cite{Sz}).
An additional proof of~\eqref{e-phi} can be found in \cite{GWW}, and 
also an alternative derivation of formula~\eqref{e-rains} is given 
in \cite{Rains}.
For the convenience of the reader we provide (yet another) 
proof of~\eqref{e-phi} 
in the Appendix to this paper.

Using the integral representation~\eqref{e-phi}, 
Johansson (\cite{Johan}) proved 
the following bound for $\phi(\lambda)$ :     
for any given $\epsilon >0$, there exist $C$ and $\delta>0$ such that 
\begin{equation}\label{e-J1.11}
  \begin{split}
    &0\le\phi_n(\lambda)\le Ce^{-\delta\lambda}
\quad\text{if}\quad (1+\epsilon)n< 2\sqrt{\lambda},\\
    &0\le 1- \phi_n(\lambda) \le \frac{C}{n}
\quad\text{if}\quad (1-\epsilon)n> 2\sqrt{\lambda}.
  \end{split}
\end{equation}
This information and the de-Poissonization Lemma~\ref{lem-dep2}  
are enough to give a new proof (\cite{Johan}) that
\begin{equation}\label{e-J1.11'}
  \lim_{N\to\infty} L_N/2\sqrt{N} =1.
\end{equation}
The first estimate in~\eqref{e-J1.11} is a consequence of the following lower
tail large deviation formula for $\phi_n(\lambda)$,
\begin{equation}\label{e-J1.12'}
\lim_{\lambda\to\infty}\frac
1{\lambda}\phi_{[x\sqrt{\lambda}]}(\lambda)=              %%  (1.12')
-1+2x-\frac 34x^2-\frac{x^2}2\log\frac 2x\equiv -U(x),
\end{equation}
if $x<2$.
For the upper tail Sepp\"al\"ainen in \cite{Se2} used
the interacting particle system implicitly
introduced by Hammersley in \cite{Ha} to show that
\begin{equation}\label{e-J1.13}
  \lim_{\lambda\to\infty} \frac1{\sqrt{\lambda}}
\log\bigl(1-\phi_{[x\sqrt{\lambda}]}(\lambda)\bigr) = -2x\cosh^{-1}(x/2)
+2\sqrt{x^2-4}\equiv -I(x)
\end{equation}
if $x>2$. We note that Hammersley's interacting particle system was also used
earlier by Aldous and Diaconis in \cite{AD}.
The super-additivity of the process described above implies that we
actually have, see \cite{Se2} and also \cite{Ki},
\begin{equation}\label{e-J1.13'}
1-\phi_{[xM]}(M^2)\le e^{-MI(x)},    %%(1.13')
\end{equation}
if $M$ is a positive integer and $x\ge 2$. This estimate can be used
to show~\eqref{e-J1.6b'}, but in this paper we will give an independent proof
of~\eqref{e-J1.6b'}.
The large deviation formula~\eqref{e-J1.13} implies, via a de-Poissonization
argument, that for $x>2$,
\begin{equation}\label{e-J1.14}
  \lim_{N\to\infty} \frac1{\sqrt{N}} \log Prob\bigl( l_N>x\sqrt{N}\bigr)
=-I(x).
\end{equation}
For the lower tail the large deviation formula for $l_N$ is not the
same as for $L(\lambda)$, the Poissonized case. Deuschel and Zeitouni
in \cite{DeZe2} use combinatorial and variational
ideas from Logan and Shepp \cite{LS} to prove that
\begin{equation}\label{e-z50}
  \lim_{N\to\infty} \frac1{N} \log Prob\bigl( l_N<x\sqrt{N}\bigr)
=-H(x),
\end{equation}
if $0<x<2$, where
\begin{equation}\label{e-J15}
  H(x)=-\frac12 +\frac{x^2}8+\log\frac{x}2
-\bigl(1+\frac{x^2}4\bigr)\log\bigl(\frac{2x^2}{4+x^2}\bigr).
\end{equation}
For the lower tail we have no analogue of~\eqref{e-J1.13'}. The rate functions
$U$ and $H$ are related via a Legendre transform, see \cite{Se2}.
The above results show clearly that the distribution function for
$l_N$ is sharply concentrated in the region
$\{ (2-\epsilon)\sqrt{N} < l_N < (2+\epsilon)\sqrt{N} \}$ for
any $\epsilon >0$, and they can be used to see heuristically that the
variance for $l_N$ should be of order $N^{1/3}$, see \cite{Ki}.

\bigskip
As is well known (see, \cite{Sz}) the Toeplitz determinant $D_{n-1}$ 
in~\eqref{e-phi} 
is intimately connected with the polynomials  
$p_n(z;\lambda)=\kappa_n(\lambda)z^n+\cdots$, 
which are orthonormal with respect to 
the weight $f(e^{i\theta})\frac{d\theta}{2\pi}=
\exp(\sqrt{\lambda}(z+z^{-1})) \frac{dz}{2\pi iz}$
on the unit circle, 
\begin{equation}
   \int_{-\pi}^{\pi} p_n(e^{i\theta})\overline{p_m(e^{i\theta})} 
f(e^{i\theta})\frac{d\theta}{2\pi}
=\delta_{n,m} \quad\text{for $n,m\geq 0$}.
\end{equation}
The leading coefficient $\kappa^2_n(\lambda)$ can be expressed
in terms of Toeplitz determinants,
\begin{equation}
   \kappa_n^2(\lambda)=\frac{D_{n-1}(\lambda)}{D_n(\lambda)}
\end{equation}
where $D_n(\lambda)=D_n(\exp(2\sqrt{\lambda}\cos\theta))$.
But by Szeg\"o's strong limit theorem (\cite{Sz2})
for Toeplitz determinants, 
$\lim_{n\to\infty}D_{n}(\lambda)=e^{\lambda}$, and hence 
\begin{equation}
  \log\phi_n(\lambda) =  \sum_{k=n}^{\infty} \log{\kappa^2_{k}(\lambda)}.
\end{equation}
Therefore, if one can control the large $k,\lambda$
behavior of $\kappa^2_{k}(\lambda)$
for all $k\geq n$, one will control the large 
$n,\lambda$ behavior of $\phi_n(\lambda)$.

The key point in our analysis is that
$\kappa^2_{k}(\lambda)$ can be expressed in terms of 
the following Riemann-Hilbert Problem (RHP) :
Let $\Sigma$ be the unit circle oriented counterclockwise.
Let $Y(z;k+1,\lambda)$ be the $2\times 2$ matrix-valued function satisfying
\begin{equation}\label{e-Y}
  \begin{cases}
     Y(z;k+1,\lambda) \quad\text{is analytic in}\quad \C-\Sigma,\\
     Y_+(z;k+1,\lambda)=Y_-(z;k+1,\lambda) \begin{pmatrix}
1&\frac{1}{z^{k+1}}e^{\sqrt{\lambda}(z+z^{-1})}\\0&1 \end{pmatrix}
\quad\text{on}\quad \Sigma,\\
     Y(z;k+1,\lambda) z^{-(k+1)\sigma_3}
=I+O(\frac1{z}) \quad\text{as}\quad z\to\infty,
  \end{cases}
\end{equation}
where $Y_+$ and $Y_-$ denote the limit from inside and outside of
the circle respectively, and $\sigma_3=
\bigl( \begin{smallmatrix} 1&0\\0&-1 \end{smallmatrix} \bigr)$, 
so that $z^{-(k+1)\sigma_3}=
\bigl( \begin{smallmatrix} z^{-(k+1)}&0\\
0&z^{k+1} \end{smallmatrix} \bigr)$.
Here $I$ is the $2\times 2$ identity matrix.
This RHP has a unique solution  
(see~\eqref{e-beg1} below), and the fact of the matter is that
\begin{equation}\label{e-Y21}
   \kappa^2_{k}(\lambda)=-Y_{21}(0;k+1,\lambda)
\end{equation}
where $Y_{21}(0;k,\lambda)$ is the $(21)$-entry of the solution
$Y$ at $z=0$.
In \cite{DZ1} and \cite{DZ}, Deift and Zhou introduced a steepest descent 
type method to compute the asymptotic behavior of RHP's containing 
large oscillatory and/or exponentially growing/decaying factors as in~\eqref{e-Y}.
This method was further extended in \cite{DVZ2} and eventually placed 
in a very general form by Deift, Zhou and Venakides in \cite{DVZ}, 
making possible the analysis of the limiting behavior 
of a large variety of asymptotic problems in pure and applied mathematics 
(see, e.g., \cite{DIZ}).
As we will see, the application of this method to~\eqref{e-Y} makes it possible 
to control the large $k,\lambda$ behavior of $\kappa_k^2(\lambda)$. 
The calculation in this paper have many similarities to the computations 
in [DKMVZ], where the authors use the steepest descent method to obtain 
Plancherel-Rotach type asymptotics for polynomials orthogonal with respect of 
varying weights, $e^{-NV(x)}dx$ on the real line, and hence to prove  
universality for a class of random matrix models.
The Riemann-Hilbert formulation of the theory of orthogonal 
polynomials on the line is due to Fokas, Its and Kitaev (\cite{FIK}) : 
the RHP~\eqref{e-Y} is an adaptation of the construction in \cite{FIK} 
to the case of orthogonal polynomial with respect to a weight 
on the unit circle.

\bigskip
This paper is arranged as follows. 
In Section~\ref{s-PII}, we discuss of some of the basic theory of RHP's 
and also provide some information on the RHP associated with the PII equation.
This information will be used in the construction of an approximate solution,  
i.e. a \emph{parametrix},  
for the RHP~\eqref{e-Y} in  
subsequent sections. The appearance of the 
PII equation in the limiting distribution 
$F(t)$ for $\chi_N$  
originates in this construction of the parametrix.
A connection of $\phi_n(\lambda)$ to Toda lattice and the Painlev\'e III equation 
is presented in Section~\ref{s-TL}.
Section~\ref{s- Begin} is the starting point for the 
analysis of the RHP~\eqref{e-Y}.
In this section, ~\eqref{e-Y} is transformed into an equivalent 
RHP via  
a so-called $g$-function.
The role of $g$-function, first introduced in \cite{DZ}, and then 
analyzed in full generality in \cite{DVZ}, 
is to replace exponentially growing terms in a RHP by oscillatory or 
exponentially decreasing terms.
It turns out that in the case of~\eqref{e-Y}, as in \cite{DKMVZ}, the 
$g$-function can be constructed in terms of an associated equilibrium 
measure $d\mu(s)$ as follows, 
\begin{equation}
   g(z)\equiv\int_{\Sigma} \log(z-s) d\mu(s).
\end{equation}
The measure $d\mu$ is the unique minimizer of the following variational problem : 
\begin{equation}\label{e-variational}
   E^{V}=\inf \{I^{V}(\wt{\mu}) : \text{$\wt{\mu}$ is a probability measure on the 
unit circle $\Sigma$}\}
\end{equation}
where 
\begin{equation}
   I^{V}(\wt{\mu})= \iint_{\Sigma\times\Sigma} \log|s-w|^{-1} 
d\wt{\mu}(s) d\wt{\mu}(w) +\int_{\Sigma} V(s) d\wt{\mu}(s)
\end{equation}
and $V(s)=-\sqrt{\lambda}(s+s^{-1})$. 
The variational problem~\eqref{e-variational} describes the equilibrium configuration 
of electrons, say, confined to the unit circle with Coulomb interactions, 
and acted on by an external field $V$.
It turns out that the support of the equilibrium measure depends  
critically on the quantity 
\begin{equation}
   \gamma=\frac{2\sqrt{\lambda}}{k+1}.
\end{equation}
We need to distinguish these two cases, $\gamma\le 1$ and $\gamma>1$.
As noted by Gross and Witten (\cite{GW}), and also by Johansson 
(\cite{Johan}), the point $\frac{2\sqrt{\lambda}}{k+1}=1$ corresponds 
to a (third order) phase transition for a statistical system 
with partition function~\eqref{e-DiS}.
The first case, when $\gamma=\frac{2\sqrt{\lambda}}{k+1}\le 1$, 
is discussed in 
Section~\ref{s-K1}, and the second case,  
when $\gamma=\frac{2\sqrt{\lambda}}{k+1}>1$, is discussed 
in Section~\ref{s-K2}.
The principal results of the above two sections are summarized in 
Lemmas~\ref{lem-K1} and ~\ref{lem-K2}. 
We obtain full asymptotics of $\kappa^2_{n}(\lambda)$ 
for $n, \lambda >0$ when $n, \lambda \to\infty$.
In Section~\ref{s-phi}, by summing up $\kappa^2_{k}(\lambda)$
for all $k\geq n$, we obtain the asymptotics of $\phi_n(\lambda)$ 
in Lemma~\ref{lem-phi}.
The relation between $\phi_n(N)$ and $q_{n,N}$ 
(de-Poissonization) is discussed in 
Section~\ref{s- dep}.
Finally, the proofs of Theorem~\ref{thm1} and~\ref{thm2} are given in 
Section~\ref{s-pf}.

\bigskip
\textbf{Notational remarks} : The primary variables in this paper are $n,N$, 
and $\lambda$. 
The letters $C$, $c$ denote general positive constants.
Rather than introducing many such constants $C_1,C_2,\cdots,c_1,c_2,\cdots,$ 
we always interpret $C,c$ in a general way. 
For example, we write $|f(x)|\le 2C|g(x)|+e^c|h(x)|\le C(|g(x)|+|h(x)|)$, etc. 
We will also use certain auxiliary positive parameters 
$M,M_1 ,M_2 ,\cdots,M_7$. 
If a constant depends on some of these parameters, we indicate 
this explicitly, for example, $C(M_2 ,M_4)$. 
In addition to the standard big $O$ notation, we also use a notation 
$O_M$. 
Thus $f=O(\frac1{n^{1/3}})$ means $|f|\le \frac{C}{n^{1/3}}$, where 
$C$ is independent of $M,M_1,\cdots$.
On the other hand, $f=O_M(\frac1{n^{1/3}})$ means 
$|f|\le \frac{C(M,M_1,\cdots)}{n^{1/3}}$, 
where $C(M,M_1,\cdots)$ depends on at least one of the parameters 
$M,M_1,\cdots$.

In the estimates that follow we will often claim that an inequality 
is true ``as $n\to\infty$''. 
For example, in~\eqref{e-z35} below, we say that 
\begin{equation*}
   |\log{\phi_n(\lambda)}| \le C\exp\biggl(-c(n+1)
\bigl(1-\frac{2\sqrt{\lambda}}{n+1}\bigr)^{3/2}\biggr),
\end{equation*}
as $n\to\infty$.
This mean that there exists a number $n_0$, say, which may 
depend on \emph{all} the other relevant constants in the problem, 
such that the inequality is true for $n\geq n_0$, etc.
(For this particular inequality the only other parameter is 
$M_5$, but it turns out that the constants $C$, $c$ can be 
chosen independent of $M_5$ (see below).)

%%%=====================secPII

\section{{\bf Riemann-Hilbert Theory}}
\label{s-PII}

In this section, we first summarize some basic facts about RHP's in general,   
and then discuss the RHP for the PII equation in some detail.
Basic references for RHP's are \cite{CG}, \cite{GK}, and 
the material on PII is taken from \cite{DZ}.

\begin{figure}[ht]
 \centerline{\epsfig{file=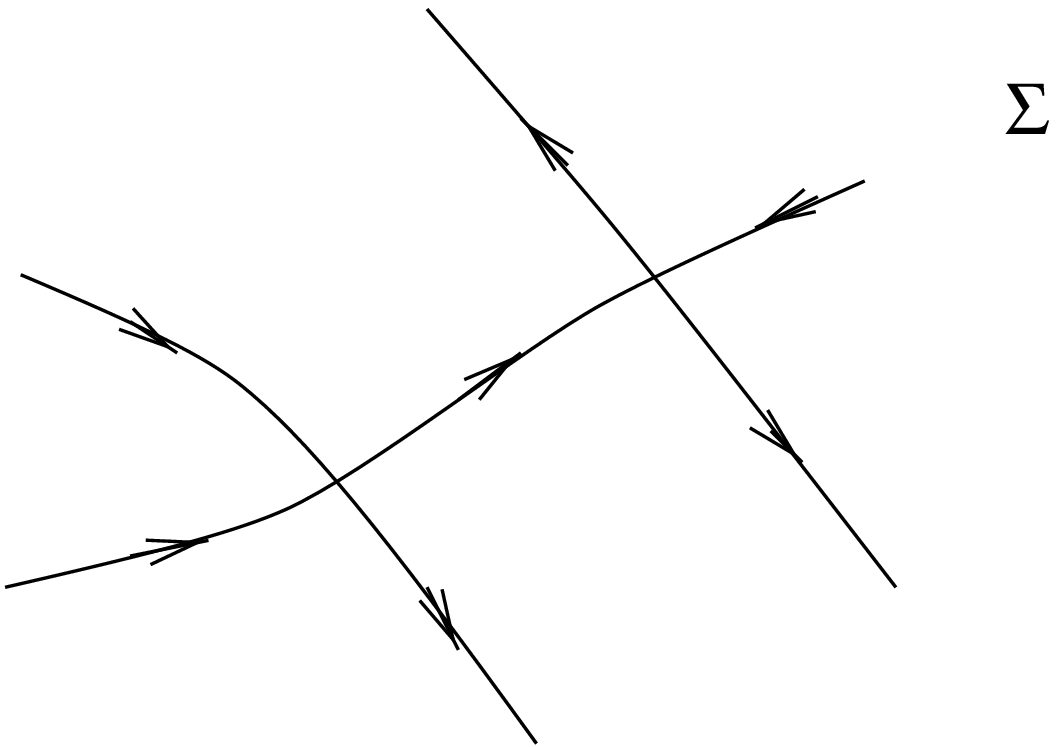, width=6cm}}
 \caption{}\label{fig-rh1}
\end{figure}

\begin{figure}[ht]
 \centerline{\epsfig{file=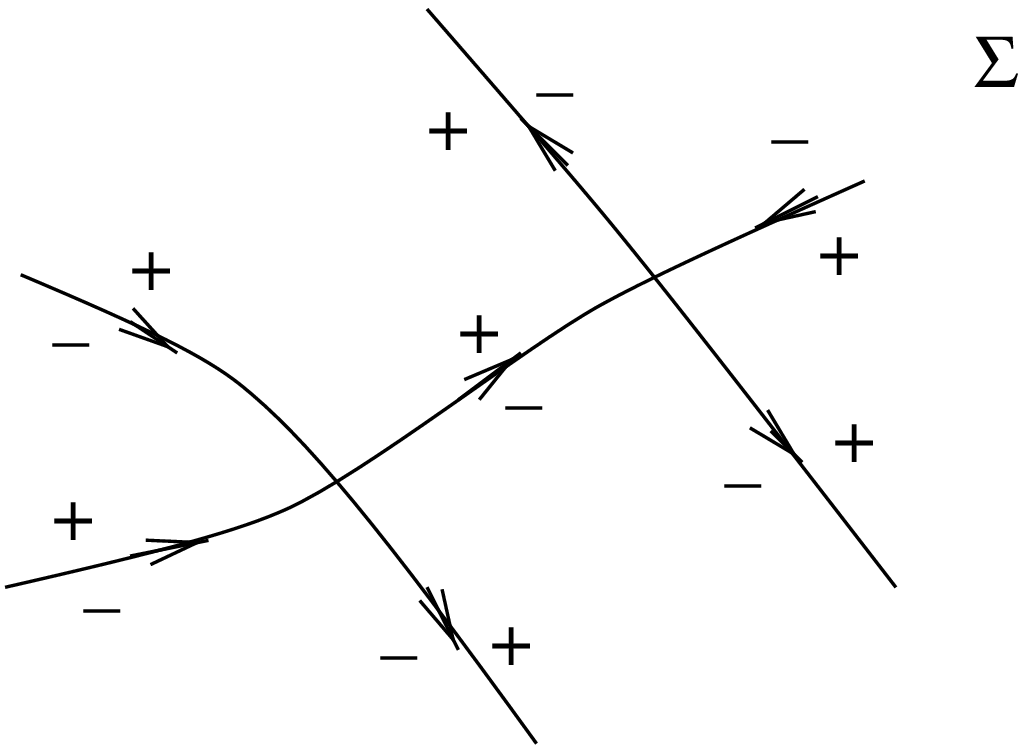, width=6cm}}
 \caption{}\label{fig-rh2}
\end{figure}

\medskip
Let $\Sigma$ be an oriented curve in the plane (see, 
for example, Figure~\ref{fig-rh1}).
By convention, the $(+)$-side (resp., $(-)$-side) of an arc in $\Sigma$ 
lies to the left (resp., right) as one traverses the arc in the direction 
of the orientation.
Thus, corresponding to Figure~\ref{fig-rh1}, we have Figure~\ref{fig-rh2}.
Let $\Sigma_0=\Sigma-\{\text{points of self-intersection}\}$ and $v$ 
be a smooth map from $\Sigma_0\to Gl(n,\C)$, for some $n$.
If $\Sigma$ is unbounded, we require that $v(z)\to I$ as $z\to\infty$ 
along $\Sigma$.
The RHP $(\Sigma,v)$ consists of the following (see, e.g., \cite{CG}) : 
establish the existence and uniqueness of an $n\times n$ matrix valued 
function $Y(z)$ (the \emph{solution} of the RHP $(\Sigma,v)$) 
such that 
\begin{equation}\label{e-tec1}
  \begin{cases}
     Y(z) \qquad\text{is analytic $\C-\Sigma$,}\\
     Y_+(z)=Y_-(z)v(z) \quad z\in\Sigma_0,\\
     Y(z)\to I \qquad\text{as $z\to\infty$.}
  \end{cases}
\end{equation}
Here $Y_{\pm}(z)=\lim_{z'\to z}Y(z')$ where 
$z'\in (\pm)-\text{side of $\Sigma$}$.
The precise sense in which these boundary values are attained, and also 
the precise sense in which $Y(z)\to I$ as $z\to\infty$, are 
technical matters 
that should be specified for any given RHP $(\Sigma,v)$. 
In \emph{this} paper, by a solution $Y$ of a RHP $(\Sigma,v)$, we always mean 
that 
\begin{equation}\label{e-tec2}
  \begin{split}
   &\text{$Y(z)$ is analytic in $\C-\Sigma$ and continuous up to the 
boundary}\\
&\text{(including the points in $\Sigma-\Sigma_0$) in each component.}\\ 
&\text{The jump relation $Y_+(z)=Y_-(z)v(z)$ is taken in the sense of}\\
&\text{continuous boundary values, and $Y(z)\to I$ as $z\to\infty$ means}\\
&\text{$Y(z)=I+O\bigl( \frac1{|z|} \bigr)$ uniformly as 
$z\to\infty$ in $\C-\Sigma$.}
  \end{split}
\end{equation}
Given $(\Sigma,v)$, the existence of $Y$ under appropriate 
technical assumptions on $\Sigma$ and $v$, is in general 
a subtle and difficult question.
However, for the RHP~\eqref{e-Y}, and hence for all RHP's obtained by 
deforming~\eqref{e-Y} (see, e.g.,~\eqref{e-m1orig}), we will prove the 
existence of $Y$ directly by construction (see, Lemma~\ref{lem-exi}) : 
uniqueness, as we will see, is a simple matter.

The solution of a RHP $(\Sigma, v)$ can be expressed in terms of 
the solution of an associated singular integral equation on $\Sigma$ 
(see,~\eqref{e-sam4},~\eqref{e-tec4} below) as follows.
Let $C_\pm$ be the Cauchy operators
\begin{equation}\label{e-sam3}
   (C_\pm f)(z)=\lim_{z'\to z_\pm} \int_{\Sigma}
\frac{f(s)}{s-z'} \frac{ds}{2\pi i}, \quad z\in\Sigma,
\end{equation}
where $z'\to z_\pm$ denotes the non-tangential limit from the $\pm$ side
of $\Sigma$ respectively.
A useful reference for Cauchy operators on curves which may have points of 
self-intersection is \cite{GK}.
Under mild assumptions on $\Sigma$, which will always be satisfied 
for the curves that arise in this paper, the non-tangential limits 
in~\eqref{e-sam3} will exist pointwise a.e.~on $\Sigma$.
Furthermore, if $f\in L^p(\Sigma,|dz|)$, $1<p<\infty$, then 
the boundary values (appropriately interpreted at the points 
$\Sigma-\Sigma_0$ of self-intersection) 
of $\int_{\Sigma} \frac{f(s)}{s-z} \frac{ds}{2\pi i}$ 
are also taken in the sense of $L^p$ and 
$\|C_\pm f\|_{L^p(\Sigma,|dz|)} \le c_p\|f\|_{L^p(\Sigma,|dz|)}$.
A simple calculation shows that  
\begin{equation}\label{e-tec5}
  C_+-C_-=1.
\end{equation}

Let 
\begin{equation}\label{e-sam2}
    v=b^{-1}_-b_+ \equiv (I-w_-)^{-1}(I+w_+)
\end{equation}
be any factorization of $v$.
We assume $b_\pm$, and hence $w_\pm$, are smooth on $\Sigma_0$, and 
if $\Sigma$ is unbounded, we assume $b_\pm(z)\to I$ as $z\to\infty$ 
along $\Sigma$.
Define the operator
\begin{equation}\label{e-tec3}
   C_w(f)\equiv C_+(fw_-)+C_-(fw_+).  
\end{equation}
By the above discussion, if $w_\pm\in L^{\infty}(\Sigma,|dz|)$, then 
$C_w$ is bounded from $L^2(\Sigma,|dz|)\to L^2(\Sigma,|dz|)$.
Suppose that the equation 
\begin{equation}\label{e-sam4}
   (1-C_w)\mu=I \quad\text{on $\Sigma$}
\end{equation}
has a solution $\mu\in I+L^2(\Sigma)$,
Or more precisely, suppose $\mu-I\in L^2(\Sigma)$ solves 
\begin{equation}\label{e-tec4}
   (1-C_w)(\mu-I)=C_wI=C_+(w_-)+C_-(w_+),
\end{equation}
which is a well-defined 
equation in $L^2(\Sigma)$ provided that 
$w_\pm\in L^\infty\cap L^2(\Sigma,|dz|)$.
Then the solution of the RHP~\eqref{e-tec1} is given by 
(see, \cite{CG},\cite{BC}) 
\begin{equation}\label{e-sam5}
   Y(z)=I+\int_{\Sigma} \frac{\mu(s)(w_+(s)+w_-(s))}{s-z}
\frac{ds}{2\pi i}, \quad z\notin\Sigma.
\end{equation}
Indeed for a.e. $z\in\Sigma$, from~\eqref{e-sam4} and~\eqref{e-tec5}, 
\begin{equation*}
  \begin{split}
     Y_+(z)&=I+C_+(\mu(s)(w_+(s)+w_-(s)))\\
     &=I+C_w(\mu) + (C_+-C_-)(\mu w_+)\\
     &=\mu+\mu w_+\\
     &=\mu(z)b_+(z),
  \end{split}
\end{equation*}
and similarly $Y_-(z)=\mu(z)b_-(z)$, so that $Y_+(z)=Y_-(z)b_-^{-1}(z)b_+(z)
=Y_-(z)v(z)$ a.e. on $\Sigma$.
Under the appropriate regularity assumptions on $\Sigma$ and $v$, 
one then shows that $Y(z)$ solves the RHP $(\Sigma, v)$ 
in the sense of~\eqref{e-tec2}.

As indicated, the above approach to the RHP goes through for any 
factorization $v=(I-w_-)^{-1}(I+w_+)$. 
Different factorization may be used at different points in the 
analysis of any given problem (see e.g. \cite{DZ1}).
However, in this paper we will \emph{always} take $w_-=0$, so that 
$v=(I+w_+)$.
Thus $C_w$ always denotes the operator $C_-\bigl(\cdot(v-I) \bigr)$.

In this paper we will not develop the general theory for the solution 
of RHP's, giving conditions under which~\eqref{e-sam4} has a (unique) 
solution, etc. 
Rather, for the convenience of the reader who may not be familiar with 
Riemann-Hilbert theory, we will use the above calculations 
and computations as a guide, 
and verify all the steps directly as they arise.

%%===========

\bigskip
We now consider the RHP for the PII equation 
(\cite{FN}, \cite{JMU} : see also \cite{IN}, \cite{FZ}, \cite{DZ}).
We will consider two equivalent versions of the RHP for PII.
These two RHP's  
will be used in the later sections for the construction of parametrices 
for the solution of~\eqref{e-Y}.

\begin{figure}[ht]
 \centerline{\epsfig{file=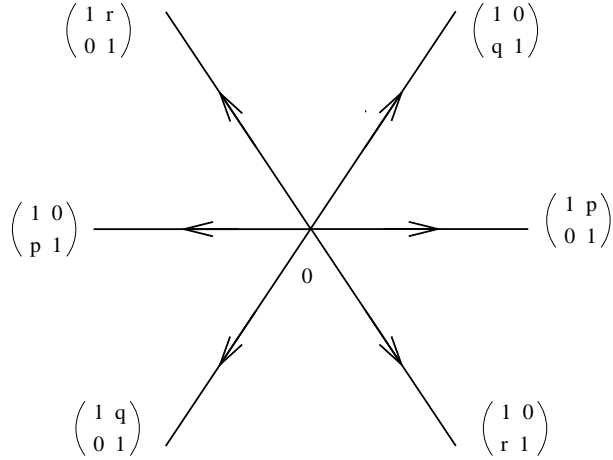, width=8cm}}
 \caption{$v^{PII}$ and $\Sigma^{PII}$}\label{fig-PII}
\end{figure}

Let $\Sigma^{PII}$ denote the oriented contour consisting of 6 rays 
in Figure~\ref{fig-PII}.
Thus $\Sigma^{PII}=\cup_{k=1}^{6}\{ \Sigma^{PII}_{k}=e^{i(k-1)\pi/3} \R_+ \}$,
with associated jump matrix $v^{PII} : \Sigma^{PII} \to M_2(\C)$, 
where the monodromy data $p,q$ and $r$ are 
complex numbers satisfying the relation
\begin{equation}\label{e-mon}
  p+q+r+pqr=0.
\end{equation}
For $x\in\R$ and $z\in\Sigma^{PII}-\{0\}$, set
\begin{equation}
  \begin{split}
     v_x^{PII}(z)&=e^{-i\theta_{PII}\sigma_3}v^{PII}
e^{i\theta_{PII}\sigma_3}\\
     &\equiv e^{-i\theta_{PII}ad\sigma_3}v^{PII}, 
  \end{split}
\end{equation}
where
\begin{equation}
    \theta_{PII}\equiv \frac{4z^3}3+xz.
\end{equation}
The contour $\Sigma^{PII}$ consists precisely of the set $Re(i4z^3/3)=0$.
This implies, in particular, that $v_x^{PII}(z)-I \notin L^2(\Sigma^{PII})$.
For example, as $z\to +\infty$ along the real axis, 
$v_x^{PII}(z)-I$ is oscillatory 
(on the other rays, $\Sigma^{PII}_{k}$, $k=1,2,4,5$, $v_x^{PII}(z)-I$ 
could grow), and so we cannot expect that the RHP 
$(\Sigma^{PII},v_x^{PII})$ has a solution in the sense of~\eqref{e-tec2}.
However, if we rotate $\Sigma^{PII}$ in the clockwise direction by any  
angle $\theta_0$, $0<\theta_0<\pi/3$, $\Sigma^{PII}\to 
\Sigma^{PII}_{\theta_0}\equiv e^{-i\theta_0}\Sigma^{PII}$, then it is 
easy to see that $v_x^{PII}(z)-I\in 
L^2\cap L^{\infty}(\Sigma^{PII}_{\theta_0})$, 
and we may expect that the RHP $(\Sigma^{PII}_{\theta_0},v_x^{PII})$ 
has a solution in the sense of~\eqref{e-tec2}.
Moreover, as $v_x^{PII}(z)$ is analytic, it is clear that if one can 
solve $(\Sigma^{PII}_{\theta_0},v_x^{PII})$ for some $0<\theta_0<\pi/3$, 
then one can solve $(\Sigma^{PII}_{\wt{\theta}_0},v_x^{PII})$ 
for any other $0<\wt{\theta}_0<\pi/3$, and the solution of the 
$\wt{\theta}_0$-problem can be obtained from the $\theta_0$-problem 
by an analytic continuation, and vice versa.
So suppose that for some fixed $0<\theta_0<\pi/3$, and for 
$x\in\R$,  
$m_{\theta_0}^{PII}(z;x)$ is a ($2\times 2$ matrix) solution of the RHP 
$(\Sigma_{\theta_0}^{PII}, v_x^{PII})$, 
\begin{equation}\label{e-mPII}
  \begin{cases}
     m^{PII}_{\theta_0}(z) \qquad  
\text{analytic in} \qquad   \C -\Sigma_{\theta_0}^{PII},\\
     \bigl(m_{\theta_0}^{PII}\bigr)_+(z)=
\bigl(m_{\theta_0}^{PII}\bigr)_-(z) v_x^{PII}(z) , 
\qquad 0\neq z\in \Sigma_{\theta_0}^{PII},\\
     m^{PII}_{\theta_0}(z)\to I \qquad \text{as} \qquad z\to\infty, 
  \end{cases}
\end{equation}
in the sense of~\eqref{e-tec2}.
Let $m_1^{PII}(x)$ denote the residue at $\infty$ of 
$m_{\theta_0}^{PII}(z)$, given by
\begin{equation*}
   m_{\theta_0}^{PII}(z;x)=I+\frac{m_1^{PII}(x)}{z}+O(\frac1{z^2})
\end{equation*}
as $z\to\infty$.
Then 
\begin{equation}\label{e-z4}
   u(x)\equiv 2im_{1,12}^{PII}(x)=-2im_{1,21}^{PII}(x)
\end{equation}
solves PII (see \cite{FN}, \cite{JMU}), 
\begin{equation*}
   u_{xx}=2u^3+xu,\qquad x\in\R,
\end{equation*}
where $m_{1,12}^{PII}(x)$ (resp., $m_{1,21}^{PII}(x)$) denotes 
the $(12)$-entry (resp, $(21)$-entry) of $m_1^{PII}(x)$.
It is easy to see that $m_1^{PII}(x)$, and hence $u(x)$ in~\eqref{e-z4}, 
is independent of the choice of $\theta_0\in (0,\pi/3)$ 

A solution of the RHP $(\Sigma_{\theta_0}^{PII}, v_x^{PII})$ 
for some $\theta_0$, hence for all $\theta_0\in (0,\pi/3)$, may not exist 
for all $p,q,r$ satisfying~\eqref{e-mon} and $x\in\R$.
A sufficient condition (see \cite{FZ})  
for the RHP to have a unique solution (in the sense of~\eqref{e-tec2}) 
for all $x\in\R$, is that
\begin{equation*}
   |q-\bar{p}|<2 \qquad \text{and} \qquad r\in\R.  
\end{equation*}
In this paper, we need the singular case 
\begin{equation}
   p=-q=1 \qquad \text{and} \qquad r=0.
\end{equation}
\begin{figure}[ht]
 \centerline{\epsfig{file=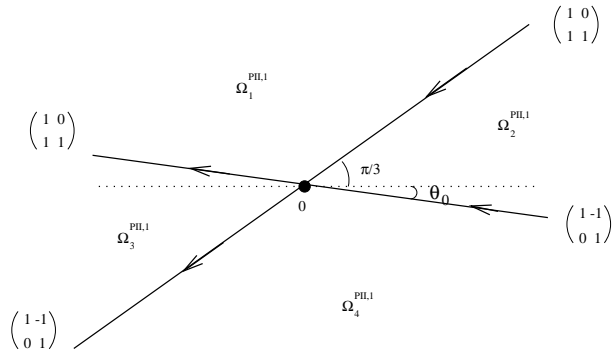, width=8cm}}
 \caption{$p=-q=1,r=0$ case ; $\Sigma_{\theta_0}^{PII,1}$ and 
$v^{PII,1}$}\label{fig-PII1}
\end{figure}
The latter condition $r=0$ implies that 
there is no jump across the rays $\pm e^{i(2\pi/3-\theta_0)}$, 
and we may replace $\Sigma^{PII}_{\theta_0}$ by $\Sigma^{PII,1}_{\theta_0}$ 
as in Figure~\ref{fig-PII1} (note that the orientations across the rays 
$e^{-i\theta_0}, e^{i(2\pi/3-\theta_0)}$ have been reversed). 
As noted in \cite{DZ}, a unique solution in the sense of~\eqref{e-tec2} 
still exists in this singular case 
for all $x\in\R$ : a proof of this fact is not given in \cite{DZ}, 
but can be found in [DKMVZ3 ; nonregular case, Case II]. 
In addition, the solution has the property that 
\begin{equation}\label{e-unif}
  \begin{split}
   &m_{\theta_0}^{PII,1}(z;x) \ \ \text{and its inverse are uniformly bounded}\\
&\text{for}\ \  (z,x)\in(\C-\Sigma_{\theta_0}^{PII,1})\times[-M,M],
  \end{split}
\end{equation}
for any fixed $M>0$.
As $m^{PII,1}_{\theta_0}(z;x)$ solves~\eqref{e-mPII} 
in the sense of~\eqref{e-tec2}, we see 
in particular that~\eqref{e-unif} holds up to the boundary in each sector.

The asymptotics of $u(x)=2im_{1,12}^{PII,1}(x)$  
given in~\eqref{e-asyPII}, is computed in ~\cite{DZ} via the above RHP and 
from the proof in \cite{DZ}, one learns that 
\begin{equation}\label{e-aympm_1}
  \begin{split}
   &m^{PII,1}_{1,22}(x)=O\biggl(\frac{e^{-(4/3)x^{3/2}}}{x^{1/4}}\biggr) 
\quad \text{as} \quad x\to\infty,\\
   &m^{PII,1}_{1,22}(x)\sim \frac{i}8 x^2 \quad \text{as} \quad x\to -\infty,
  \end{split}
\end{equation}
where $m^{PII,1}_{1,22}$ denotes the $(22)$-entry of $m^{PII}_{1}$.
Also, using the methods in \cite{DZ}, for example, 
one obtains the relation
\begin{equation}\label{e-um1PII}
   \frac{d}{dx} 2im^{PII,1}_{1,22}(x) = u^2(x).
\end{equation}
and verifies directly that $2im^{PII,1}_{1,22}(x)$ is real-valued.

\bigskip
\begin{figure}[ht]
 \centerline{\epsfig{file=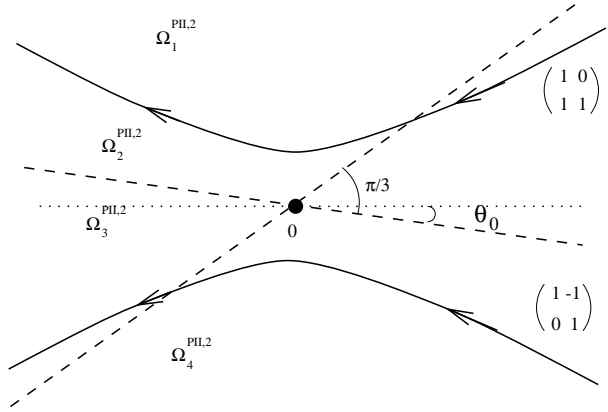, width=8cm}}
 \caption{$\Sigma^{PII,2}$ and $v^{PII,2}$}\label{fig-PII2}
\end{figure}
For the first of the two equivalent RHP's advertised above, 
we consider Figure~\ref{fig-PII2}, which consists of 
the real axis (the dotted line), $\Sigma^{PII,1}_{\theta_0}$ for some 
fixed, small $\theta_0>0$ (the dashed lines), and a contour 
$\Sigma^{PII,2}$ consisting of a pair of curved solid lines. 
The contour $\Sigma^{PII,2}$ is of the general shape indicated 
in the Figure, with one component in $\C_+$ and one component in $\C_-$, 
and we require that $\Sigma^{PII,2}$ is asymptotic to straight lines 
lying \emph{strictly} within 
the region $\{ |arg z|<\pi/3 \}\cup\{ 2\pi/3 <arg z< 4\pi/3 \}$.
Together with the line $\{xe^{-i\theta_0} : x\in\R \}$, 
these contours divide 
the complex plane into 4 open regions, 
$\Omega^{PII,2}_k$, $k=1,2,3,4$, as shown in Figure ~\ref{fig-PII2}.
Let $v^{PII,2}$ be the jump matrix on $\Sigma^{PII,2}$ which is given by 
$\bigl(\begin{smallmatrix} 1&0\\1&1 \end{smallmatrix}\bigr)$ in $\C_+$ 
and by $\bigl(\begin{smallmatrix} 1&-1\\0&1 \end{smallmatrix}\bigr)$ in $\C_-$.
We define
\begin{equation*}
  \begin{cases}
    m^{PII,2}=m_{\theta_0}^{PII,1}e^{-i(\theta_{PII})ad\sigma_3}
\bigl(\begin{smallmatrix} 1&0\\1&1 \end{smallmatrix} \bigr)^{-1}
\quad\text{in}\quad{\Omega^{PII,2}_1\cap\Omega^{PII,1}_2},\\
    m^{PII,2}=m_{\theta_0}^{PII,1}e^{-i(\theta_{PII})ad\sigma_3}
\bigl(\begin{smallmatrix} 1&0\\1&1 \end{smallmatrix} \bigr)
\qquad\text{in}\quad{\Omega^{PII,2}_2\cap\Omega^{PII,1}_1},\\
    m^{PII,2}=m_{\theta_0}^{PII,1}e^{-i(\theta_{PII})ad\sigma_3}
\bigl(\begin{smallmatrix} 1&-1\\0&1 \end{smallmatrix} \bigr)^{-1}
\quad\text{in}\quad{\Omega^{PII,2}_3\cap\Omega^{PII,1}_4},\\
    m^{PII,2}=m_{\theta_0}^{PII,1}e^{-i(\theta_{PII})ad\sigma_3}
\bigl(\begin{smallmatrix} 1&-1\\0&1 \end{smallmatrix} \bigr)
\qquad\text{in}\quad{\Omega^{PII,2}_4\cap\Omega^{PII,1}_3},\\
    m^{PII,2}=m_{\theta_0}^{PII,1}\quad\text{otherwise},
  \end{cases}
\end{equation*}
where the regions $\Omega^{PII,1}_{k}$,$k=1,2,3,4$ are defined  
in Figure ~\ref{fig-PII1}.
A straightforward calculation with the jump relations for 
$m_{\theta_0}^{PII,1}$, 
shows that 
$m^{PII,2}$ solves the new RHP
\begin{equation}\label{e-RHPm2}
  \begin{cases}
    m^{PII,2}\quad\text{is analytic in}\quad\C-\Sigma^{PII,2},\\
    m^{PII,2}_+=m^{PII,2}_-v^{PII,2}_x\quad\text{on}
\quad\Sigma^{PII,2},\\
    m^{PII,2}\to I\quad\text{as}\quad z\to\infty,
  \end{cases}
\end{equation}
where $v^{PII,2}_{x}=e^{-i(\theta_{PII})ad\sigma_3}v^{PII,2}$ 
and $v^{PII,2}$ is given in Figure~\ref{fig-PII2}.
This deformed RHP is clearly equivalent to the original RHP 
for $m_{\theta_0}^{PII,1}$
in the sense that a solution of the one RHP yields a solution 
of the other RHP, and vice versa.
Also we have 
\begin{equation}\label{e-RHPm2'}
   (m_{\theta_0}^{PII,1})_1=m_1^{PII,2},
\end{equation}
for the residues of $m_{\theta_0}^{PII,1}$ (resp, $m^{PII,2}$) at $\infty$.
From~\eqref{e-unif}, 
we see that for any fixed $M>0$, 
\begin{equation}\label{e-unif2}
 \begin{split}
   &m^{PII,2}(z;x) \quad\text{and its inverse are uniformly bounded}\\
   &\text{for}\quad  (z,x)\in(\C-\Sigma^{PII,2})\times[-M,M].
 \end{split}
\end{equation}
A particular choice of contour $\Sigma^{PII}$ will be made in 
Section~\ref{s-K1} (see below).

\bigskip
\begin{figure}[ht]
 \centerline{\epsfig{file=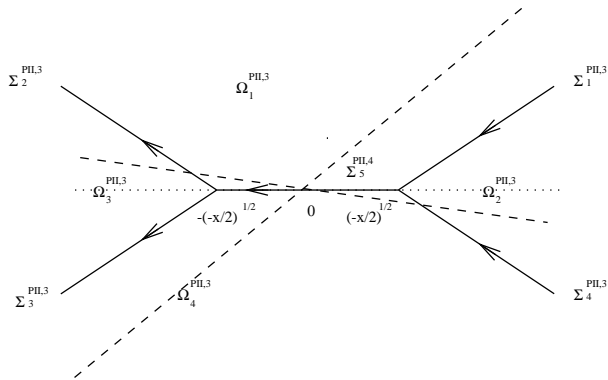, width=8cm}}
 \caption{$\Sigma^{PII,3}$}\label{fig-PII3}
\end{figure}
The second of the equivalent RHP's is restricted to the case $x<0$, 
and we consider 
Figure~\ref{fig-PII3}, which consists of the real axis (the dotted line), 
$\Sigma^{PII,1}_{\theta_0}$ for some fixed small $\theta_0>0$ 
(the dashed lines) and a contour $\Sigma^{PII,3}=
\cup_{k=1}^{5}\Sigma^{PII,3}_{k}$ consisting of 5 straight lines, one 
finite and four infinite.
The regions $\Omega^{PII,3}_k$, $1\le k\le 4$, are the components 
of $\C-\Sigma^{PII,3}$.

The infinite lines make an angle strictly between $0$ and $\pi/3$ 
with the real axis.
Set 
\begin{equation}\label{e-z5}
   g^{PII}(z)=\frac43 \bigl(z^2+\frac{x}2 \bigr)^{3/2}
\end{equation}
which is defined to be analytic in $\C-[-\sqrt{\frac{-x}2}, 
\sqrt{\frac{-x}2}]$, and behaves like 
$\frac43 z^3+xz+\frac{x^2}{8z}+O(\frac1{z^3}) = \theta_{PII}(z)+O(\frac1{z})$ 
as $z\to\infty$.
Therefore for any $M>0$, 
\begin{equation}\label{e-z6}
    e^{i(g^{PII}(z)-\theta_{PII}(z))} \quad \text{is bounded for}
\quad (z,x)\in (\C-[-\sqrt{\frac{-x}2},\sqrt{\frac{-x}2}]) \times [-M,0]
\end{equation}
and
\begin{equation}\label{e-z7}
    e^{i(g^{PII}(z)-\theta_{PII}(z))} \to 1 \quad \text{as $z\to\infty$ 
uniformly for $-M\le x\le 0.$}
\end{equation}
We define $m^{PII,3}$ by 
\begin{equation*}
  \begin{cases}
    m_{\theta_0}^{PII,1}[e^{-i(\theta_{PII})ad\sigma_3}
\bigl(\begin{smallmatrix} 1&0\\1&1 \end{smallmatrix} \bigr)^{-1}]
e^{i(g^{PII}(z)-\theta_{PII}(z))\sigma_3}
\quad\text{in}\quad 
\Omega^{PII,3}_1\cap(\Omega^{PII,1}_2\cup\Omega^{PII,1}_3),\\
    m_{\theta_0}^{PII,1}[e^{-i(\theta_{PII})ad\sigma_3}
\bigl(\begin{smallmatrix} 1&0\\1&1 \end{smallmatrix} \bigr)]
e^{i(g^{PII}(z)-\theta_{PII}(z))\sigma_3}
\qquad\text{in}\quad{\Omega^{PII,3}_3\cap\Omega^{PII,1}_1},\\
    m_{\theta_0}^{PII,1}[e^{-i(\theta_{PII})ad\sigma_3}
\bigl(\begin{smallmatrix} 1&-1\\0&1 \end{smallmatrix} \bigr)^{-1}]
e^{i(g^{PII}(z)-\theta_{PII}(z))\sigma_3}
\quad\text{in}\quad{\Omega^{PII,3}_2\cap\Omega^{PII,1}_4},\\
    m_{\theta_0}^{PII,1}[e^{-i(\theta_{PII})ad\sigma_3}
\bigl(\begin{smallmatrix} 1&-1\\0&1 \end{smallmatrix} \bigr)]
e^{i(g^{PII}(z)-\theta_{PII}(z))\sigma_3}
\qquad\text{in}\quad
\Omega^{PII,3}_4\cap(\Omega^{PII,1}_2\cup\Omega^{PII,1}_3),\\
    m_{\theta_0}^{PII,1}\quad\text{otherwise}.
  \end{cases}
\end{equation*}
Then from the jump relations for $m_{\theta_0}^{PII,1}$,  
we see that $m^{PII,3}$ solves the new RHP 
$(\Sigma^{PII,3}, v_{x}^{PII,3})$ in the sense of~\eqref{e-tec2}, 
\begin{equation}\label{e-z17}
  \begin{cases}
    m^{PII,3}\quad\text{is analytic in}\quad\C-\Sigma^{PII,3},\\
    m^{PII,3}_+=m^{PII,3}_-v^{PII,3}_x\quad\text{on}
\quad\Sigma^{PII,3},\\
    m^{PII,3}\to I\quad\text{as}\quad z\to\infty,
  \end{cases}
\end{equation}
where $v^{PII,3}_{x}$ is given by 
\begin{equation}\label{e-vPII4}
  \begin{cases}
     \begin{pmatrix}1&0\\e^{2ig^{PII}}&1 \end{pmatrix}
\qquad\text{on}\quad \Sigma_1^{PII,3},\Sigma_2^{PII,3}\\
     \begin{pmatrix}1&-e^{-2ig^{PII}}\\0&1 \end{pmatrix}
\quad\text{on}\quad \Sigma_3^{PII,3},\Sigma_4^{PII,3}\\
     \begin{pmatrix}e^{-2ig_{-}^{PII}}&-1\\1&0 \end{pmatrix}
\quad\text{on}\quad \Sigma_5^{PII,3}.\\
  \end{cases}
\end{equation}
Also we have
\begin{equation}\label{e-mPmP4}
   m_1^{PII}=m_1^{PII,3}-\bigl( \frac{ix^2}{8} \bigr)\sigma_3,
\end{equation}
for the respective residues of $m_{\theta_0}^{PII,1}$ and $m^{PII,3}$ 
at $\infty$. 
Finally, from~\eqref{e-unif} and~\eqref{e-z6}, 
we see that, for any fixed $M\in\R$
\begin{equation}\label{e-unif4}
  \begin{split}
   &m^{PII,3}(z;x) \quad\text{and its inverse are uniformly bounded}\\
&\text{for}\quad  (z,x)\in(\C-\Sigma^{PII,3})\times[-M,0].
  \end{split}
\end{equation}

%%%=======================secTL

\section{{\bf Connection to the Toda Lattice and the Painlev\'e III Equation}}
\label{s-TL}

In this Section, we discuss the connection of the RHP~\eqref{e-Y} 
for $\kappa^2_k$ and the RHP for the Toda lattice 
and the Painlev\'e III equation. 
In the RH context, the connection results from 
the specific form of the weight, $e^{\sqrt{\lambda}(z+z^{-1})}$. 
Connections can also be seen from the Toeplitz determinant/orthogonal 
polynomial point of view as in \cite{PS}, \cite{Hi} and \cite{Wi}.
The purpose of this short Section is purely to establish the 
various connection, but we do not use the results in the sequel.

\bigskip
Write $q=k+1$ in~\eqref{e-Y}.
We define 
\begin{equation}\label{e-sTL1}
   m^{TL}(z;q)=
   \begin{cases}
\biggl( \begin{smallmatrix} 0&-1\\1&0 \end{smallmatrix} \biggr) Y(z;q)
\biggl( \begin{smallmatrix} z^{-q}e^{\sqrt{\lambda}z^{-1}}&0\\
0&z^{q}e^{-\sqrt{\lambda}z^{-1}} \end{smallmatrix} \biggr)
\biggl( \begin{smallmatrix} 0&1\\-1&0 \end{smallmatrix} \biggr) \quad |z|>1,\\
\biggl( \begin{smallmatrix} 0&-1\\1&0 \end{smallmatrix} \biggr)Y(z;q)
\biggl( \begin{smallmatrix} e^{\sqrt{\lambda}z}&0\\
0&e^{-\sqrt{\lambda}z} \end{smallmatrix} \biggr) \quad |z|<1.
 \end{cases}
\end{equation}
A simple calculation shows that $m^{TL}$ solves the following RHP,
\begin{equation}\label{e-TL3}
  \begin{cases}
   &m^{TL}(z) \quad\text{is analytic in} \quad \C-\Sigma,\\
   &m^{TL}_+(z)=m^{TL}_-(z)\begin{pmatrix}
0&-z^{-q}e^{-\sqrt{\lambda}(z-z^{-1})}\\
z^{q}e^{\sqrt{\lambda}(z-z^{-1})}&1 \end{pmatrix}\\
   &m^{TL}(z)\to I \quad\text{as} \quad z\to\infty.
  \end{cases}
\end{equation}
Once again, the RHP for $Y$ is equivalent to the RHP for $m^{TL}$ 
in the sense that a solution of one problem yields a solution
of the other problem.

Recall that the RHP related to the Toda Lattice problem, for 
$-\infty<m<\infty$, 
\begin{equation}
  \begin{split}
     \frac{da_m}{dt}=2(b^2_{m}-b^2_{m-1})\\
     \frac{db_m}{dt}=b_{m}(a_{m+1}-a_{m}), 
  \end{split}
\end{equation}
under initial data decaying at infinity is the following (see,e.g.,\cite{Ka}).
Suppose that there are no solitons and denote the   
reflection coefficient by $r(z)$, $z\in\Sigma$.  
Then we find $Q(z)$ such that 
\begin{equation}
  \begin{cases}
     Q(z) \qquad\text{is analytic in $\C-\Sigma$}\\
     Q_+(z)=Q_-(z)\begin{pmatrix} 1-|r(z)|^2&-\bar{r}(z)z^{2m}e^{-t(z-z^{-1})}\\
r(z)z^{-2m}e^{t(z-z^{-1})}&1 \end{pmatrix}\\
     Q(z)\to I \qquad\text{as $z\to\infty$}.
  \end{cases}
\end{equation}
When $q$ is even, if we set $\sqrt{\lambda}=t$ and $q=-2m$ in~\eqref{e-TL3}, 
then the RHP is identical with the above RHP 
with $r(z)\equiv 1$.

For the connection to the Painlev\'e III equation, 
define 
\begin{equation}
  m^{PIII}(z)=
  \begin{cases}
     (-1)^qm^{TL}(z) &\quad |z|<1,\\
     m^{TL}(z) &\quad |z|>1.
  \end{cases}
\end{equation}
Note in~\eqref{e-TL3}, 
\begin{equation}
   (-1)^q\negthickspace\begin{pmatrix}
0\negthickspace&\negthickspace-z^{-q}e^{-\sqrt{\lambda}(z-z^{-1})}\\
z^{q}e^{\sqrt{\lambda}(z-z^{-1})}&1 \end{pmatrix}\negthickspace
= \negthickspace(-1)^qe^{-\frac{\sqrt{\lambda}}2(z-z^{-1})ad\sigma_3}z^{-\frac{q}2ad\sigma_3}  
\negthickspace\begin{pmatrix} 1&1\\-1&0 \end{pmatrix} ^{-1}
\end{equation}
where $z^{\frac{q}2}$ is analytic in 
$\C-(-\infty,0]$ and real-valued for real $z$.
If we set $\sqrt{\lambda}=-ix$, then
this is the same RHP for the particular  
Painlev\'e III equation (see \cite{FMZ} for results and notations)
\begin{equation*}
   u_{xx}=\frac{u_x^2}{u}-\frac1{x}u_x
+\frac1{x}\bigl( -4qu^2+4(1-q)\bigr) +4u^3+\frac{-4}{u}
\end{equation*}
with monodromy data 
\begin{equation*}
  \begin{split}
    &\theta_{\infty}=-\theta_0=q,\\
    &a_0=b_0=a_{\infty}=b_{\infty}=0,\\
    &E=\bigl( \begin{smallmatrix} 1&1\\-1&0 \end{smallmatrix} \bigr).
  \end{split}
\end{equation*}

In the RHP~\eqref{e-Y}, we are interested directly in the quantity 
$-Y_{21}(0;k+1,\lambda)$, or by~\eqref{e-sTL1}, $m^{TL}_{11}(0;q)$. 
On the other hand, for the Toda lattice and the PIII equation, one is 
interested in quantities other than $m^{TL}_{11}(0;q)$ which are 
related to the respective RHP's. 
For example, the solution $u(x)$ of PIII equation is given by 
$u(x)=-ix(m^{PIII}_1)_{12}$ where $m^{PIII}=I+\frac{m^{PIII}_1}{z}
+O(\frac1{z})$, which is clearly different from $(-1)^qm^{PIII}_{11}(0;q)$.
However, the importance of the connection of~\eqref{e-Y} to the RHP's 
for Toda lattice and the PIII equation lies precisely in the fact that 
$(a_m,b_m)$ (resp., $u(x)$) solve differential-difference 
(resp., differential) equations which in turn imply that the coefficients  
of the generating function $\phi_n(\lambda)$, $2\sqrt{\lambda}=-ix$, 
must satisfy a certain class of identities. 
We plan to investigate these relations in a later publication.

Finally, note that for PIII, the interesting asymptotic question is 
to evaluate the limit $x=i\sqrt{\lambda}\to\infty$, with $q$ fixed. 
In this paper, as in the Toda lattice, we are interested in the double 
limit when $\lambda\to\infty$ and $q$ is allowed to vary 
(note that in \cite{Ka}, the singular case $r(z)\equiv 1$ is not 
considered). 
When $\lambda\to\infty$, $\frac{2\sqrt{\lambda}}{q}\sim 1$, 
we are in a region where the solution of PIII equation degenerates 
to a solution of PII equation, and this explains the appearance of 
PII in the parametrix for the solution of $Y$ of the RHP~\eqref{e-Y}.

%%%===================sec1.tex

\section{{\bf Equilibrium Measure and $g$-function}}
\label{s- Begin}

In this Section, the equilibrium measure is explicitly calculated 
for each $\gamma>0$ (Lemma~\ref{lem-ourg}) and, using this equilibrium measure, 
the $g$-function~\eqref{4.8} is introduced in order to 
convert the RHP~\eqref{e-Y} 
into a RHP which is normalized to be $I$ at $\infty$.

\bigskip
Let $\Sigma$ denote the unit circle oriented counterclockwise 
and $f(e^{i\theta})=f(z)$ be a non-negative, periodic, smooth  
function on $\Sigma$.
Let $p_q(z)=\kappa_qz^q+\cdots$ be the $q$-th normalized orthogonal 
polynomial with respect to the weight 
$f(e^{i\theta})\frac{d\theta}{2\pi}$ on the unit circle.
Define the polynomial $p_q^{*}(z)\equiv z^q\bar{p}_q(1/z)
=z^q\overline{p_q(1/\bar{z})}$ (see \cite{Sz}).
We consider the following RHP : 
Let $Y(z)$ be the $2\times 2$ matrix-valued function satisfying 
\begin{equation}\label{e-beg1}
  \begin{cases}
    Y(z) \qquad  \text{is analytic in} \quad  \C -\Sigma,\\
    Y_+(z)=Y_-(z)\begin{pmatrix}1&\frac{1}{z^{k+1}}f(z)\\0&1
\end{pmatrix} \quad\text{on} \quad  \Sigma,\\
    Y(z)z^{-(k+1)\sigma_3}=I+O(\frac1z) \quad\text{as} \quad z\rightarrow \infty.
  \end{cases}
\end{equation}
The following Lemma is the starting point of our calculations.
\begin{lem}\label{lem-exi}
  (cf. \cite{FIK}, \cite{DKMVZ}) The RHP~\eqref{e-beg1} has a unique solution
  \begin{equation*}
    Y(z)=\begin{pmatrix}
    \frac1{\kappa_{k+1}}p_{k+1}(z)& \frac1{\kappa_{k+1}}
\int_{\Sigma}\frac{p_{k+1}(s)}{s-z}\frac{f(s)ds}{2\pi i s^{k+1}}\\
    -\kappa_{k}p^{*}_{k}(z)& -\kappa_{k}
\int_{\Sigma}\frac{p^{*}_{k}(s)}{s-z}\frac{f(s)ds}{2\pi i s^{k+1}}
       \end{pmatrix}.
  \end{equation*}
\end{lem}
\begin{proof}
  Existence : Using the property of Cauchy operator $C_+-C_-=I$, where
$Ch(z)\equiv\int_{\Sigma} \frac{h(s)}{s-z} ds$,
it is a straightforward calculation
to show that the above expression for $Y$ satisfies the jump condition.
The asymptotics at $\infty$ codes in precisely the fact that the $p_k's$ 
are the normalized orthogonal polynomials for the weight 
$f(e^{i\theta})d\theta$. 

  Uniqueness : Suppose that there is another solution $\tilde{Y}$ of RHP. Noting

$det \biggl( \begin{smallmatrix} 1&\frac{1}{z^{k+1}}f(z)\\0&1
\end{smallmatrix} \biggr) =1$, we have that $\det \tilde{Y}$ is entire, 
and $\to 1$
as $z\to\infty$. Therefore by Liouville's theorem, $\det \tilde{Y}\equiv 1$.
In particular, $\tilde{Y}$ is invertible.
Now set $Z=Y\tilde{Y}^{-1}$. Then it has no jump on $\Sigma$ 
hence is entire. 
Also, $Z\to I$ as $z\to\infty$, and therefore $Z\equiv I$.
\end{proof}
From this Lemma,  we have
\begin{equation}
  \kappa^2_{k}=-Y_{21}(0)
\end{equation}
Therefore the RHP~\eqref{e-Y} has a unique solution 
and~\eqref{e-Y21} is verified.

\bigskip
Again set $q=k+1$ and 
\begin{equation}
  \gamma=\frac{2\sqrt{\lambda}}{q}.
\end{equation}
We are interested in the case when $q$ and $2\sqrt{\lambda}$ 
are of the same order
, or more precisely, $\gamma\to 1$.
In this Section, and also in Sections~\ref{s-K1} and~\ref{s-K2}, we consider 
the RHP~\eqref{e-Y} with parameter $\gamma$ and $q$, 
\begin{equation}\label{RHP-Y}
  \begin{cases}
     Y(z;q) \quad \text{analytic in} \quad  \C -\Sigma,\\
     Y_+(z;q)=Y_-(z;q)
\begin{pmatrix}1&\frac{1}{z^q}e^{\frac{q\gamma}{2}(z+z^{-1})}\\0&1
\end{pmatrix} \quad\text{on}\quad \Sigma,\\
     Y(z;q)=\bigl( I+O(\frac1z) \bigr)z^{q\sigma_3} 
\quad\text{as}\quad z\rightarrow \infty,
  \end{cases}
\end{equation}
rather than $\lambda$ and $q$.
With $\gamma$ fixed, the RHP~\eqref{RHP-Y} is of the Plancherel-Rotach type  
with \emph{varying} exponential weight 
$e^{\frac{q\gamma}{2}(z+z^{-1})}$ on the unit circle 
(see \cite{Sz}, \cite{DKMVZ}).
Similar problem on the real line is analyzed in \cite{DKMVZ} 
without double scaling limit ($\gamma$ is kept fixed).
Our goal is to find the large $q$ behavior of $Y_{21}(0;q)$ for all $\gamma > 0$.
\begin{figure}[ht]
 \centerline{\epsfig{file=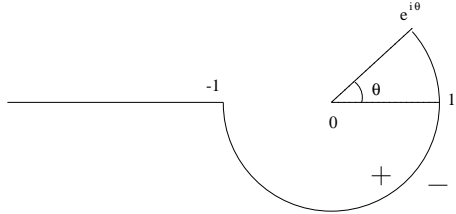, width=6cm}}
 \caption{branch cut of $\log (z-e^{i\theta})$}\label{fig-branch}
\end{figure}

Let $d\mu(s)$ be a probability measure on the unit circle.
Define 
\begin{equation}
  g(z) \equiv \int_{\Sigma}\log (z-s) d\mu(s)
\end{equation}
where for each $\theta$, the branch is chosen such that $\log(z-e^{i\theta})$
is analytic in $\C-(-\infty,-1]\cup \{ e^{it}:-\pi\le t\le\theta \}$
 (see Figure ~\ref{fig-branch})
and $\log(z-e^{i\theta}) \sim \log{z}$ for real $z\to\infty$.
The following Lemma is based on related calculations in \cite{DKM}.
\begin{lem}\label{lem-gfun}
Suppose $d\mu(z)=u(\theta) d\theta$ is an absolutely continuous probability 
measure on the unit circle and $u(\theta)=u(-\theta)$. 
Then $g(z)$ has the following properties : 
\begin{enumerate}
\item $g$ is analytic in $\C -\Sigma\cup (-\infty ,-1).$
\item On $(-\infty ,-1), g_+(z)-g_-(z)=2\pi i.$
\item $g(z)=\log z+O(\frac1z)$  as  $z\rightarrow \infty.$
\item $e^{qg(z)}$ is analytic in $\C-\Sigma.$
\item $e^{qg(z)}=z^q( 1+O(\frac1z))$ as $z\rightarrow \infty.$
\item $g(0)= \pi i.$
\item $g_+(z)+g_-(z)=2\int_{-\pi}^{\pi}\log |z-s| d\mu(s)
      + i(\phi +\pi)$ on $z\in\Sigma$  where $\phi=arg(z).$
\item $g_+(z)-g_-(z)=2\pi i\int_{\phi}^{\pi} d\mu(s)$ on z$\in\Sigma.$
\end{enumerate}
\end{lem} 

\begin{proof}
  (i)-(v) are trivial. 
For (vi), 
$$g(0)=\int_{-\pi}^{\pi} \log{(0-e^{i\theta})} u(\theta)d\theta
=\int_{-\pi}^{\pi} i(\theta+\pi) u(\theta)d\theta = \pi i$$
using the evenness of $u(\theta)$.

For (vii), fix $z=e^{i\phi}\in\Sigma$. 
Then $arg(z-e^{i\theta})$ is analytic if $-\pi<\theta<\phi$ and 
  \begin{equation*}
    \begin{split}
      g_+(z)&=\int_{-\pi}^{\pi}\log |z-e^{i\theta}| d\mu(z)
+i\int_{-\pi}^{\phi} arg(z-e^{i\theta}) d\mu(z)
+i\int_{\phi}^{\pi} arg_+(z-e^{i\theta}) d\mu(z),\\
      g_-(z)&=\int_{-\pi}^{\pi}\log |z-e^{i\theta}| d\mu(z)
+i\int_{-\pi}^{\phi} arg(z-e^{i\theta}) d\mu(z)
+i\int_{\phi}^{\pi} arg_-(z-e^{i\theta}) d\mu(z).
    \end{split}
  \end{equation*}
  Note that for $\phi <\theta <\pi$, 
$$arg_+(e^{i\phi}-e^{i\theta}) - arg_-(e^{i\phi}-e^{i\theta}) = 2\pi.$$
  This yields
$$g_+(z)+g_-(z)=2\int_{-\pi}^{\pi}\log |z-e^{i\theta}| d\mu(z)
+2i\int_{-\pi}^{\pi}arg_+(e^{i\phi}-e^{i\theta}) d\mu(z)
-i\int_{\phi}^{\pi} 2\pi d\mu(z).$$
  Set 
$$F(\phi)=2\int_{-\pi}^{\phi} arg(e^{i\phi}-e^{i\theta}) d\mu(z)
+2\int_{\phi}^{\pi} arg_+(e^{i\phi}-e^{i\theta}) d\mu(z)
-2\pi\int_{\phi}^{\pi} d\mu(z) -\phi$$
  If we show $F(\phi)\equiv\pi$, then (vii) is proved.
  Note that (a) $arg(e^{i\phi}-e^{i\phi_-})=\phi+\frac{\pi}2$, 
(b) $arg_+(e^{i\phi}-e^{i\phi_+})=\phi+\frac{3\pi}2$ and 
(c) $\frac{d}{d\phi} arg(e^{i\phi}-e^{i\theta}) \equiv \frac12$.
  This gives us $F'(\phi)=2arg(e^{i\phi}-e^{i\phi_-})u(\phi)-2arg_+(e^{i\phi}-e^{i\phi_+})u(\phi)+2\pi u(\phi) \equiv 0$. But $F(\pi)=\pi$. Therefore $F(\phi)\equiv\pi$

  For (viii), 
  \begin{equation*}
    \begin{split}
      g_+(z)-g_-(z)&=i\int_{\phi}^{\pi}
[arg_+(e^{i\phi}-e^{i\theta}) - arg_-(e^{i\phi}-e^{i\theta})] d\mu(z)\\
      &=i\int_{\phi}^{\pi} 2\pi d\mu(z).
    \end{split}
  \end{equation*}
\end{proof}

\bigskip
Let $M$ be the set of probability measures on $\Sigma$.
The equilibrium measure $d\mu_V(z)$ for potential $V(z)=-\frac{\gamma}2(z+z^{-1})$
on the unit circle is defined by the following minimization problem,
\begin{equation}
   \inf_{\mu\in M} \iint_{\Sigma\times\Sigma} \log|z-w|^{-1} d\mu(z) d\mu(w)
+ \int_{\Sigma} V(z) d\mu(z).
\end{equation}
The infimum is achieved uniquely (see, e.g. \cite{ST}) at the 
equilibrium measure.
Let $J$ denote the support of $d\mu_V$.
The equilibrium measure and its support are uniquely determined 
by the following Euler-Lagrange variational conditions :
\begin{equation}\label{e-varicon}
  \begin{split}
  &\text{there exits a real constant $l$ such that},\\
  &\quad 2\int_{\Sigma} \log{|z-s|}d\mu_V(s) - V(z)+l=0 \ \ \text{for $z\in \bar{J}$},\\
  &\quad 2\int_{\Sigma} \log{|z-s|}d\mu_V(s) - V(z)+l \le 0 \ \
\text{for $z\in\Sigma -\bar{J}$}.
  \end{split}
\end{equation}
In Lemma ~\ref{lem-ourg} below, we find $d\mu_V$,
its support and $l$ explicitly from this variational condition 
with the aid of Lemma ~\ref{lem-gfun}.
Let 
\begin{equation}\label{4.8}
  g(z)=g_V(z) \equiv \int_{\Sigma}\log (z-s) d\mu_V(s)
\end{equation}
where $d\mu_V$ is the equilibrium measure.
Following \cite{DKMVZ}, we define 
$$\m1 \equiv 
e^{\frac{ql}{2}\sigma_3}Y(z)e^{-qg(z)\sigma_3} e^{-\frac{ql}{2}\sigma_3}.$$
Then $m^{(1)}$ solves the following new RHP,
\begin{equation}\label{e-m1orig}
  \begin{cases} 
    \m1 \qquad\text{is analytic in}\quad \C-\Sigma,\\
    \mpl1 = \mni1 v^{(1)} \quad\text{on}\quad\Sigma,\\
    \m1=I+O(\frac1z) \quad\text{as}\quad z\rightarrow \infty
  \end{cases}
\end{equation}
where $v^{(1)}=\begin{pmatrix} e^{q(g_--g_+)}& 
\frac{1}{z^q}e^{q(g_++g_--V+l)}\\
0& e^{q(g_+-g_-)} \end{pmatrix}$, 
and 
\begin{equation}\label{e-Korig}
   \kappa^2_{q-1}=-Y_{21}(0;q)=-m_{21}^{(1)}(0)e^{ql}e^{qg(0)}
=-(-1)^qm^{(1)}_{21}(0)e^{ql},
\end{equation}
from Lemma~\ref{lem-gfun} (vi).

Once again we note that this RHP for $m^{(1)}$ is equivalent to 
the RHP for $Y$ in the sense that 
a solution of one RHP yields a solution of the other RHP, and vice versa.
Using Lemma ~\ref{lem-gfun}, the jump matrix $v^{(1)}$ is given by 
\begin{equation}\label{e-Omega}
  \begin{cases}
    \text{inside the support of $d\mu_V$},\\
    \qquad\begin{pmatrix} e^{-2q\pi i\int_{\phi}^{\pi} d\mu_V(\theta)}& (-1)^q\\
          0& e^{2q\pi i\int_{\phi}^{\pi} d\mu_V(\theta)} \end{pmatrix}.\\
    \text{outside the support of $d\mu_V$}, \\
    \qquad\begin{pmatrix} e^{-2q\pi i\int_{\phi}^{\pi} d\mu_V(\theta)}& 
(-1)^qe^{q\bigl[ 2\int_{-\pi}^{\pi}\log |z-e^{i\theta}| d\mu_V(\theta) 
-V(z)+l\bigr] }\\
0&e^{2q\pi i\int_{\phi}^{\pi} d\mu_V(\theta)} \end{pmatrix}.
  \end{cases}
\end{equation}
As indicated in the Introduction, the purpose of the $g$-function is to 
turn exponentially growing terms in the jump matrix for the RHP,  
into oscillatory or exponentially decaying terms : this can be seen 
explicitly in~\eqref{e-Omega}, using~\eqref{e-varicon}.

\bigskip
We have explicit formulae for the equilibrium measure and $l$.
For $0<\gamma\leq 1$, the equilibrium measure has the whole circle 
as its support but for $\gamma > 1$, a gap opens up. 
See, also \cite{GW} and \cite{Johan}.

Notation : $\chi_B(\theta)$ denotes the indicator function of the set 
$B\subset\Sigma$.

\begin{lem}\label{lem-ourg}
For the weight $V(z)=-\frac{\gamma}{2}(z+z^{-1})$, the equilibrium measure and $l$ are given as follows : 
\begin{enumerate}
\item If $0\leq \gamma \leq 1$, then
\begin{equation}\label{e-lem4.3.1} 
  d\mu_V(\theta)=\frac{1}{2\pi}(1+\gamma \cos{\theta}) d\theta
\end{equation}
and $l=0$.
\item If $\gamma > 1$, then
\begin{equation}\label{e-lem4.3.2}
  d\mu_V(\theta)=\frac{\gamma}{\pi} \cos(\frac{\theta}{2})
      \sqrt{\frac{1}{\gamma} - \sin^2(\frac{\theta}{2})} 
\chi_{[-\theta_c, \theta_c]}(\theta) d\theta
\end{equation}
and 
\begin{equation}\label{e-lem4.3.3}
  l=-\gamma +\log\gamma +1,
\end{equation}
      where $\sin^2\frac{\theta_c}{2}=\frac{1}{\gamma}, 0< \theta_c<\pi$. 
In this case, the inequality in the variational 
condition ~\eqref{e-varicon} is strict.
\end{enumerate}
\end{lem}

\begin{proof}
  (i) First, it is easy to check that $d\mu_V(\theta)$ defined above 
in~\eqref{e-lem4.3.1} is a positive probability measure. We set 
    $$g(z)=\int_{-\pi}^{\pi}\log (z-e^{i\theta})\frac{1}{2\pi}(1+\gamma \cos{\theta}) d\theta.$$
  Then $$g'(z)=\frac1{2\pi i} \int_{\Sigma} \frac1{z-s}(1+\frac{\gamma}2(s+s^{-1}))\frac{ds}{s}.$$
  Using a residue calculation with $g(z)=\log{z}+O(\frac1{z})$ as 
$z\to\infty$ and $g(0)=\frac1{2\pi} \int_{-\pi}^{\pi} 
\log{e^{i(\theta+\pi)}}(1+\gamma\cos{\theta})d\theta=\pi i$, we have
  \begin{equation}\label{e-g1} g(z)=
    \begin{cases}
      \log{z}-\frac{\gamma}{2z}& |z|>1, z\notin (-\infty,-1),\\
      -\frac{\gamma}2 z+\pi i& |z|<1.
    \end{cases}
  \end{equation}
  Therefore we have $$g_+(z)+g_-(z)=\log{z}-\frac{\gamma}2(z+z^{-1})+\pi i.$$
  From Lemma ~\ref{lem-gfun} (vii), we have
  $$2\int_{-\pi}^{\pi}\log |z-e^{i\theta}| \frac1{2\pi} 
(1+\gamma\cos{\theta}) d\theta +\frac{\gamma}2(z+z^{-1}) = 0$$
  for any $z=e^{i\phi}$ with $l=0$ as $\log{z}=i\phi$.

  (ii) It is straightforward to check that the above 
measure~\eqref{e-lem4.3.2} is 
a positive probability measure.
For $g(z)$ defined in as before, we have
  \begin{equation*}
    \begin{split}
      g'(z)&= \int_{-\theta_c}^{\theta_c} \frac1{z-e^{i\theta}} 
\frac{\gamma}{\pi} \cos(\frac{\theta}{2}) 
\sqrt{\frac1{\gamma}-\sin^2\frac{\theta}{2}} d\theta \\
      &= \frac{\gamma}{4\pi i} \int_{-\theta_c}^{\theta_c} 
\frac1{z-s}
\frac{s+1}{s^2} \sqrt{(s-\xi)(s-\bar{\xi})} ds
    \end{split}
  \end{equation*}
  where $\xi=e^{i\theta_c}$, and the branch is chosen to be analytic in 
$\C-\{ e^{i\theta} : \theta_c\le |\theta| \le \pi \}$ 
and $\sqrt{(s-\xi)(s-\bar{\xi})}>0$ for real $s>>0$. 
From a residue calculation, we obtain
  $$g'(z)= \frac1{2z}-\frac{\gamma}4(1-z^{-2}) +\frac{\gamma}4\frac{z+1}{z^2}\sqrt{(z-\xi)(z-\bar{\xi})}.$$
  Integrating, we have for $|z|>1$, $z\notin (-\infty,-1)$, 
  \begin{equation*}
    g(z)=\frac12\log{z} -\frac{\gamma}{4} (z+z^{-1}) +\frac{\gamma}2
    +\frac{\gamma i}4 \int_{1_{+0}}^{z} \frac{s+1}{s}\sqrt{(s-\xi)(s-\bar{\xi})}
    \frac{ds}{si} +g_-(1)
  \end{equation*}
and for $|z|<1$, $z\notin (-1,0]$,
  \begin{equation*}
    g(z)=\frac12\log{z} -\frac{\gamma}{4} (z+z^{-1}) +\frac{\gamma}2
    +\frac{\gamma i}4 \int_{1_{-0}}^{z} \frac{s+1}{s}\sqrt{(s-\xi)(s-\bar{\xi})}
    \frac{ds}{si} +g_+(1)
  \end{equation*}
 where $g_+$, $g_-$ denote the limit from inside and outside each and 
$1_{+0}$, $1_{-0}$ denote the outside and inside limits.

  (a) For $|\phi| \le \theta_c$,
  \begin{equation*}
    g_+(z)+g_-(z)=\log{z} -\frac{\gamma}2(z+z^{-1}) +\gamma +g_+(1)+g_-(1)
  \end{equation*}
  From Lemma ~\ref{lem-gfun} (vii), we obtain 
  \begin{equation*} 
    \begin{split}
      g_+(1)+g_-(1)&= 2\int_{-\theta_c}^{\theta_c} \log|1-e^{i\theta}| 
\frac{\gamma}{\pi} \cos{\frac{\theta}2} \sqrt{\frac1{\gamma} 
-\sin^2 \frac{\theta}2} d\theta +i(0+\pi)\\
      &= \frac{4\gamma}{\pi} \int_{0}^{\theta_c} \log(2|\sin{\frac{\theta}2}|)
\cos{\frac{\theta}2} \sqrt{\frac1{\gamma} -\sin^2 \frac{\theta}2} d\theta +\pi i\\
      &= 2\log{2} -\log{\gamma} +\frac8{\pi} \int_{0}^{\frac{\pi}2}
\log(\sin{\theta}) \cos^2 \theta d\theta +\pi i\\
      &= 2\log{2} -\log{\gamma} -1 +\frac4{\pi} \int_{0}^{\frac{\pi}2}
\log(\sin{\theta}) d\theta +\pi i,
    \end{split}
  \end{equation*}
 after simple change of variables and integration by parts.
  But we have $\int_{0}^{\frac{\pi}2} \log(\sin{\theta}) d\theta 
= -\frac{\pi}2 \log{2}$ from
   \begin{equation*}  
    \begin{split}
      2\int_{0}^{\frac{\pi}2} \log(\sin{\theta}) d\theta 
      &= \int_{0}^{\frac{\pi}2} \log(\sin{\theta}) d\theta 
+ \int_{0}^{\frac{\pi}2} \log(\cos{\theta}) d\theta \\
      &= \int_{0}^{\frac{\pi}2} \log(\frac12 \sin{2\theta}) d\theta \\
      &=-\frac{\pi}2 \log{2} 
+ \int_{0}^{\frac{\pi}4} \log(\sin{2\theta}) d\theta 
+ \int_{\frac{\pi}4}^{\frac{\pi}2} \log(\sin{2\theta}) d\theta \\
      &=-\frac{\pi}2 \log{2} +\int_{0}^{\frac{\pi}2} \log(\sin{\theta}) d\theta.
    \end{split}
  \end{equation*}
  Therefore  
$$g_+(z)+g_-(z)=\log{z} -\frac{\gamma}2(z+z^{-1})+\gamma-\log{\gamma}-1+\pi i.$$
  From Lemma ~\ref{lem-gfun} (vii), we obtain the desired result 
for $|\phi| \le \theta_c$, $z=e^{i\phi}$,
  \begin{equation*}
    2\int_{-\pi}^{\pi}\log |z-e^{i\theta}| d\mu_V(\theta)
+ \frac{\gamma}2(z+z^{-1}) -\gamma +\log\gamma +1 =0.
  \end{equation*}

  (b) for $\theta_c <\phi <\pi$ ($-\pi <\phi <-\theta_c$ case is similar),
  \begin{equation*}
    \begin{split}
      g_+(z)+g_-(z)&= \log{z} -\frac{\gamma}2(z+z^{-1}) +g_+(1)+g_-(1)\\
      &\qquad +\frac{\gamma}2 \int_{\theta_c}^{\phi} \frac{s+1}{s^2}\sqrt{(s-\xi)(s-\bar{\xi})} ds.
    \end{split}
  \end{equation*}
  But 
$$\frac{\gamma}2 \int_{\theta_c}^{\phi} \frac{s+1}{s^2}\sqrt{(s-\xi)(s-\bar{\xi})} ds = -\frac{\gamma}2 \int_{\theta_c}^{\phi} \cos{\frac{\theta}2} \sqrt{\sin^2 \frac{\theta}2 -\frac1{\gamma}} d\theta < 0.$$
  Therefore, using Lemma ~\ref{lem-gfun} (vii) and 
calculations in (a), we obtain for $|\phi|>\theta_c$,
  \begin{equation*}
    2\int_{-\pi}^{\pi}\log |z-e^{i\theta}| d\mu_V(\theta)
+ \frac{\gamma}2(z+z^{-1}) -\gamma +\log\gamma +1 <0.
  \end{equation*}
\end{proof}

In the following Sections, we distinguish the two cases, 
$\gamma\le1$ and $\gamma>1$, 
due to the difference of the supports of their equilibrium measures.

%%%%====================sec2.tex

\section{$\mathbf{ 0\leq \gamma \le 1}$}
\label{s-K1}

From ~\eqref{e-g1}, we have the explicit formula for the $g$-function :
\begin{equation*}  g(z)=
    \begin{cases}
      \log{z}-\frac{\gamma}{2z}& |z|>1, z\notin (-\infty,-1)\\
      -\frac{\gamma}2 z+\pi i& |z|<1.
    \end{cases}
\end{equation*}
With this $g$, $l=0$ from Lemma~\ref{lem-ourg} (i), and our 
RHP ~\eqref{e-m1orig},  
or equivalently ~\eqref{e-Omega}, becomes
\begin{equation}
\begin{cases}
m^{(1)} \ \ \text{analytic in} \ \ \C-\Sigma,\\
m_+^{(1)}= m_-^{(1)} \begin{pmatrix} (-1)^qz^qe^{\frac{q\gamma}{2}(z-z^{-1})}& (-1)^q\\
                              0& (-1)^qz^{-q}e^{\frac{-q\gamma}{2}(z-z^{-1})}
              \end{pmatrix} \ \ \text{on} \ \  \Sigma,\\
m^{(1)}=I+O(\frac1z) \ \ \text{as} \ \  z\rightarrow \infty\\
\end{cases}
\end{equation}
and $\kappa^2_{q-1}=-(-1)^qm_{21}^{(1)}(0)$ from ~\eqref{e-Korig}.\\
We define $m^{(2)}$ in terms of $m^{(1)}$ as follows :
\begin{equation}\label{e-conju}
 \begin{split}
  &\text{for even $q$,} \\
  &\quad \begin{cases}
    m^{(2)} \equiv m^{(1)}& \quad |z|>1,\\ 
    m^{(2)} \equiv m^{(1)} \bigl( \begin{smallmatrix} 0& -1\\ 1&0 
\end{smallmatrix} \bigr)& \quad |z|<1. \end{cases}\\
  &\text{for odd $q$,}\\
  &\quad\begin{cases}
    m^{(2)} \equiv \bigl( \begin{smallmatrix} 1&0\\0&-1 
\end{smallmatrix} \bigr)  m^{(1)} \bigl( 
\begin{smallmatrix} 1&0\\0&-1 \end{smallmatrix} \bigr)& \quad |z|>1,\\
    m^{(2)} \equiv \bigl( \begin{smallmatrix} 1&0\\0&-1 
\end{smallmatrix} \bigr) m^{(1)} \bigl( 
\begin{smallmatrix} 0&-1\\-1&0 \end{smallmatrix} \bigr)& 
\quad |z|<1. \end{cases}
 \end{split}
\end{equation}
Then we have a new equivalent RHP
\begin{equation}
\begin{cases}
m^{(2)}_+=m^{(2)}_- v^{(2)}\quad\text{on} \quad\Sigma,\\
m^{(2)}=I+O(\frac1z) \ \ \text{as}\ \ z\rightarrow \infty\\
\end{cases}
\end{equation}
where $v^{(2)}= \begin{pmatrix} 
1& -(-1)^qz^qe^{\frac{q\gamma}{2}(z-z^{-1})} \\
(-1)^qz^{-q}e^{\frac{-q\gamma}{2}(z-z^{-1})}&0 \end{pmatrix}$ 
and $\kappa^2_{q-1}=m_{22}^{(2)}(0)$. \\

\begin{figure}[ht]
 \centerline{\epsfig{file=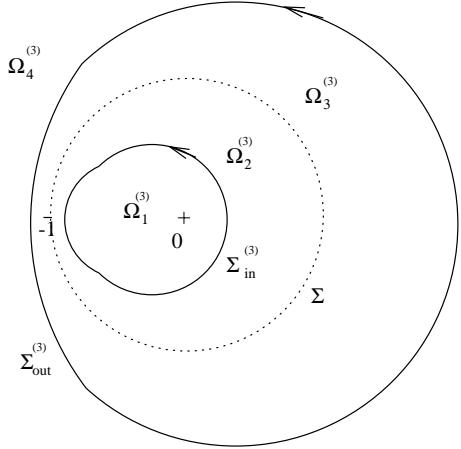, width=6cm}}
 \caption{$\Sigma^{(3)}$ and $\Omega^{(3)}$}\label{fig-ga<1}
\end{figure}
Introduce $\Sigma^{(3)}=\Sigma^{(3)}(\gamma)=\Sigma^{(3)}_{in}\cup\Sigma^
{(3)}_{out}$ (see Figure~\ref{fig-ga<1}) as follows.
For fixed $\pi/2 <|\theta|\le\pi$, 
\begin{equation}\label{e-z8}
  F(\rho)\equiv F(\rho,\theta)=
  Re\bigl( \frac{\gamma}2(z-z^{-1})+\log{z} \bigr)
  =\frac{\gamma}2 (\rho-\rho^{-1})\cos{\theta}+\log{\rho}, 
\end{equation}
where $z=\rho e^{i\theta}$, $0<\rho\le 1$, has the minimum at  
\begin{equation}\label{e-z9}
  \rho=\rho_{\theta}\equiv \frac{1-\sqrt{1-\gamma^2\cos^2 \theta}}
  {-\gamma\cos{\theta}}
\end{equation}
and 
$F(\rho_{\theta})<0$. 
( Note that $\rho_{\theta}<0$ for $|\theta|<\pi/2$. )
For $\frac12\le\gamma\le1$, we take 
\begin{equation}\label{e-Sig}
  \begin{split}
   &\Sigma^{(3)}_{in}=\{\rho_{\theta}e^{i\theta} : 
3\pi/4 \le |\theta|\le\pi \} \cup\{\rho_{3\pi/4}e^{i\theta} : 
|\theta|\le 3\pi/4 \},\\
   &\Sigma^{(3)}_{out}=\{\rho_{\theta}^{-1}e^{i\theta} :
3\pi/4 \le |\theta|\le\pi \} \cup\{\rho_{3\pi/4}^{-1}e^{i\theta} : 
|\theta|\le 3\pi/4 \}.
  \end{split}
\end{equation}
Orient $\Sigma^{(3)}$ as in Figure~\ref{fig-ga<1}. 
And finally, for $0\le\gamma\le\frac12$, set 
$\Sigma^{(3)}(\gamma)=\Sigma^{(3)}(\frac12)$.

Of course, $\Sigma^{(3)}$ varies with $\gamma\in [0,1]$. 
However, using estimates from \cite{GK}, it is not difficult to 
show that the Cauchy operators $C_\pm$ on $L^2(\Sigma^{(3)})$ 
are uniformly bounded, 
\begin{equation}\label{e-Cau}
   \|C_\pm \|_{L^2(\Sigma^{(3)})\to L^2(\Sigma^{(3)})} \le C <\infty
\end{equation}
for all $0\le\gamma\le 1$.
Observe also that in the limit $\gamma\to 1$, $\Sigma^{(3)}$ takes the 
form of the cross 
\begin{equation}\label{e-z10}
  y=\pm |x+1| \quad \text{for $z=x+iy$ near $-1$}.
\end{equation}

Apart from the neighborhood of $z=-1$, there is considerable freedom 
in the choice of $\Sigma^{(3)}$. 
For example, $3\pi/4$ could be replaced by any angle between 
$\pi/2$ and $\pi$. 
Also the form of the contour for $|\theta|<3\pi/4$ is not critical, 
as long as it has the general shape drawn in Figure~\ref{fig-ga<1} : 
all that we really need is that the jump matrix $v^{(3)}$ below 
has the property $\sup_{\{z\in\Sigma^{(3)} : |\arg(z)|<3\pi/4\}} 
|v^{(3)}-I|\to 0$ 
exponentially as $q\to\infty$.

Using the factorization 
\begin{equation*}
    v^{(2)}= \begin{pmatrix} 1&0\\
(-1)^qz^{-q}e^{\frac{-q\gamma}{2}(z-z^{-1})}&1\end{pmatrix}
    \begin{pmatrix} 1&-(-1)^qz^qe^{\frac{q\gamma}{2}(z-z^{-1})}\\0&1\end{pmatrix}
\equiv (b_-^{(2)})^{-1}b_+^{(2)},
\end{equation*}
we define 
\begin{equation*}
  \begin{cases}
     m^{(3)}=m^{(2)} (b_+^{(2)})^{-1}
\qquad\text{in}\quad \Omega^{(3)}_2,\\
     m^{(3)}=m^{(2)} (b_-^{(2)})^{-1}
\qquad\text{in}\quad \Omega^{(3)}_3,\\
     m^{(3)}=m^{(2)} 
\qquad\text{in}\quad \Omega^{(3)}_1, \Omega^{(3)}_4.
  \end{cases}
\end{equation*}
Then $m^{(3)}$ solves the RHP $(\Sigma^{(3)},v^{(3)})$ where 
\begin{equation}\label{e-sec2.5}
  \begin{cases}
    v^{(3)} = \begin{pmatrix} 1&-(-1)^qz^qe^{\frac{q\gamma}{2}(z-z^{-1})}\\0&1
              \end{pmatrix} \quad\text{on}\quad\Sigma_{in}^{(3)},\\
    v^{(3)} = \begin{pmatrix} 1&0\\(-1)^qz^{-q}e^{\frac{-q\gamma}{2}(z-z^{-1})}&1
              \end{pmatrix} \quad\text{on}\quad\Sigma_{out}^{(3)},
  \end{cases}
\end{equation}\\
and 
\begin{equation}\label{e-z1}
  \kappa^2_{q-1}= m_{22}^{(3)}(0). 
\end{equation}

As $q\to\infty$, $v^{(3)}(z)\to I$.
Set $\Sigma^{\infty}=\Sigma^{(3)}$.
The RHP 
\begin{equation}
\begin{cases}
m^{\infty}_+=m^{\infty}_-\text{I} \quad\text{on}\quad \Sigma^{\infty}, \\
m^{\infty}=I+O(\frac1z) \quad\text{as}\quad z\rightarrow \infty\\
\end{cases}
\end{equation}
has, of course, the unique solution, $m^{\infty}(z)\equiv I$. 

\bigskip
Let $0\le\gamma\le 1-\delta_1$ for some $0< \delta_1< 1$.
From the choice of $\Sigma^{(3)}$, 
\begin{equation}\label{e-sec2.13}
 \begin{split}
   \|v^{(3)}-I\|_{L^\infty(\Sigma^{(3)})} 
&=\sup_{3\pi/4 \le\theta\le 5\pi/4} |e^{F(\rho_\theta,\theta)}|
\le \sup_{3\pi/4 \le\theta\le 5\pi/4} |e^{F(\rho_{\pi},\theta)}|\\
&\le e^{F(\rho_{\pi},\pi)}
=e^{q(\sqrt{1-\gamma^2} +\log \frac{1-\sqrt{1-\gamma^2}}{\gamma})}.
 \end{split}
\end{equation}
But, for $0<\gamma \le 1$, a straightforward estimate shows that 
\begin{equation}\label{e-sec2.14}
   \sqrt{1-\gamma^2} +\log (1-\sqrt{1-\gamma^2}) -\log\gamma 
\le -\frac{2\sqrt{2}}{3} (1-\gamma)^{3/2}
\end{equation}
so that 
\begin{equation}\label{e-sec2.15}
   \|w^{(3)}\|_\infty =\|v^{(3)}-I\|_\infty
\le e^{-\frac{2\sqrt{2}}{3}\delta_1^{3/2}q} \to 0 
\quad\text{as}\quad q\to\infty.
\end{equation} 
Since $\|C_{w^{(3)}}\|_{L^2\to L^2} \le C\|w^{(3)}\|_\infty$, for 
some constant $C$ independent of $\gamma$ (see~\eqref{e-Cau}), 
$(I-C_{w^{(3)}})^{-1}$ is invertible for large $q$ and 
the solution for the RHP $(\Sigma^{(3)},v^{(3)})$ is given by 
(see~\eqref{e-sam5}) 
\begin{equation}\label{e-z2}
      m^{(3)}(z)=I+\int_{\Sigma^{(3)}} 
\frac{((I-C_{w^{(3)}})^{-1}I)(s)(v^{(3)}(s)-I)}{s-z}
\frac{ds}{2\pi i}, \quad z\notin\Sigma^{(3)},
\end{equation}
and (see~\eqref{e-z1})
\begin{equation}
   \kappa_{q-1}^2=m^{(3)}_{22}(0)
=1+\biggl( \int_{\Sigma^{(3)}} 
\frac{((I-C_{w^{(3)}})^{-1})I)(s)(v^{(3)}(s)-I)}{s}
\frac{ds}{2\pi i} \biggr)_{22}.
\end{equation}
Now from the fact that the length of $\Sigma^{(3)}$ is uniformly 
bounded and $dist(0, \Sigma) \geq c>0$ for all $\gamma\in [0,1]$, we obtain, 
\begin{equation}\label{e-sec2.18}
   |  \kappa_{q-1}^2 -1 |
\le C\|v^{(3)}-I\|_\infty
\le Ce^{-\frac{2\sqrt{2}}{3}\delta_1^{3/2}q}.
\end{equation}

%%%%==============================sec3.tex

\bigskip
The above calculation also applies to the case when $\gamma\to 1$ slowly.
Indeed, suppose $\frac12 \le \gamma \le 1-\frac{M_1}{2^{1/3}q^{2/3}}$, 
where $M_1>0$ is a fixed, sufficiently large number.
(The lower bound $\frac12$ is chosen for convenience. Any fixed number 
between 0 and 1 would work.) 
From~\eqref{e-sec2.14},~\eqref{e-sec2.15}, for some constant $C$
which is independent of $\gamma$,
\begin{equation}\label{e-sec3.6}
  \|C_{w^{(3)}}\|_{L^2\to L^2}
\le Ce^{-\frac{2\sqrt{2}}3q(1-\gamma)^{3/2}}
\le Ce^{-\frac{2\sqrt{2}}3M_1^{3/2}} \le \frac12 <1,
\end{equation}
if $M_1$ is sufficiently large. 
For convenience, we only consider $M_1\geq 1$.
From~\eqref{e-z2},
\begin{equation*}
  \begin{split}
    m^{(3)}&=I+\frac1{2\pi i} \int_{\Sigma^{(3)}} \frac{[(I-C_{w^{(3)}})^{-1
}I](s)w^{(3)}(s)}{s-z}\,ds\\
   &=I+\frac1{2\pi i} \int_{\Sigma^{(3)}} \frac{w^{(3)}(s)+[(I-C_{w^{(3)}})^
{-1}C_{w^{(3)}}I](s)w^{(3)}(s)}{s-z} \,ds
  \end{split}
\end{equation*}
and, as $diag(w^{(3)})=0$,
\begin{equation}
   \kappa^2_{q-1}=m^{(3)}_{22}(0)
    =1+\bigl( \frac1{2\pi i} \int_{\Sigma^{(3)}}
[(I-C_{w^{(3)}})^{-1}C_{w^{(3)}}I](s)w^{(3)}(s)
    \frac{ds}{s} \bigr)_{22}.
\end{equation}
Hence, we have
\begin{equation}
  \begin{split}
    |\kappa^2_{q-1}-1|&\le
\|(I-C_{w^{(3)}})^{-1}C_{w^{(3)}}I\|_{L^2}
\|\frac{w^{(3)}(s)}{2\pi is}\|_{L^2}\\
    &\le \|(I-C_{w^{(3)}})^{-1}\|_{L^2\rightarrow L^2}
\|C_{w^{(3)}}I\|_{L^2} \|\frac{w^{(3)}(s)}{2\pi is}\|_{L^2}\\
    &\le C\|w^{(3)}\|_{L^2}^2\\
    &\le C\|w^{(3)}\|_{L^\infty} \|w^{(3)}\|_{L^1}\\
    &\le Ce^{-\frac{2\sqrt{2}}3q(1-\gamma)^{3/2}} \|w^{(3)}\|_{L^1},
  \end{split}
\end{equation}
where the (final) constant $C$ is independent of $\gamma,q$ and $M_1$ 
(sufficiently large), provided that 
$0\le\gamma\le 1-\frac{M_1}{2^{1/3}q^{2/3}}$.

Since the length of $\Sigma^{(3)}$ is bounded, we have $\|w^{(3)}\|_{L^1} 
\le Ce^{-\frac{2\sqrt{2}}3q(1-\gamma)^{3/2}}$, which is 
the same estimation~\eqref{e-sec2.18} as in the case $\gamma<1-\delta_1$. 
But for future calculations (see~\eqref{e-z3} below), we need a sharper result.
We estimate $\|w^{(3)}\|_{L^1}$ as follows :
  Focus on $\Sigma^{(3)}_{in}$.
For $\Sigma^{(3)}_{out}$, similar computations apply.
Only the $12$-component of $w^{(3)}$ is non-zero.
  Set $\wt{\theta}=\frac1{q^{1/3}}\log{q}$.
  $$\int_{\Sigma^{(3)}_{in}}|z^qe^{\frac{q\gamma}{2}(z-z^{-1})}||dz|= (1)+(2
)$$
  where $(1)$ is an integration over $|\theta|\le \pi-\wt{\theta}$
and $(2)$ covers the remainder.
  Note from~\eqref{e-z10} that $|dz|\le Cd\theta$. 
Substituting $\rho_\theta$~\eqref{e-z9} into $F(\rho,\theta)$~\eqref{e-z8}, 
we obtain on $\Sigma^{(3)}_{in}$, 
\begin{equation}\label{e-z11}
    |z^qe^{\frac{q\gamma}2(z-z^{-1})}|\le e^{
q\bigl(\sqrt{1-\gamma^2\cos^2\theta}+\log(1-\sqrt{1-\gamma^2\cos^2\theta})
-\log(-\gamma\cos\theta)\bigr)}.
\end{equation}
Setting, $\gamma\to -\gamma\cos\theta$ in~\eqref{e-sec2.14}, 
we obtain for $z\in\Sigma^{(3)}_{in}$, 
\begin{equation}\label{e-z12}
    |z^qe^{\frac{q\gamma}2(z-z^{-1})}|\le
e^{-\frac{2\sqrt{2}q}{3}(1+\gamma\cos\theta)^{3/2}}
\le e^{-Cq(\pi-|\theta|)^3}
\end{equation}
Hence, adjusting the constants $C$ if necessary, 
we have 
  $$(1)\le Ce^{-Cq\wt{\theta}^3} \le \frac{C}{q^{1/3}}$$
and $$(2)\le C\int_0^{\wt{\theta}} e^{-Cq\theta^3} d\theta
\le C\int_0^{\log{q}} e^{-Ct^3} \frac{dt}{q^{1/3}} \le \frac{C}{q^{1/3}}.$$
Therefore, 
\begin{equation}\label{6.9}
   \|w^{(3)}\|_{L^1}\le \frac{C}{q^{1/3}}
\end{equation}
and we obtain
\begin{equation}\label{e-sec3f}
   |\kappa^2_{q-1}-1|\le
\frac{C}{q^{1/3}}e^{-\frac{2\sqrt{2}}3q(1-\gamma)^{3/2}}.
\end{equation}

\bigskip
Let $M_2>0$ be a fixed number 
and consider $1-\frac{M_2}{2^{1/3}q^{2/3}} \le \gamma \le 1$.
For this case, as $q\to\infty$, $\gamma\to 1$ and 
$\rho_{\theta=\pi} \to 1$. 
We need to devote special attention to the neighborhood of $z=-1$, 
where we will introduce a parametrix for the RHP,  
which is related to the special solution of the Painlev\'e II (PII)
equation~\eqref{e-int1} given in Section~\ref{s-PII}.
For a discussion of parametrices in RHP's, see e.g. \cite{DZ}, \cite{DKMVZ}.
Set $\gamma=1-\frac{t}{2^{1/3}q^{2/3}}$. 
The region above corresponds to $0\le t\leq M_2$.
Let $\mathcal{O}$ be a small neighborhood of size $\epsilon$ around
$z=-1$, where $\epsilon>0$ is a fixed number which is small enough  
so that first, 
\begin{equation}\label{e-z18}
  \text{the map $u$ defined below is a bijection from $\mathcal{O}$,} 
\end{equation}
and second, 
\begin{equation}\label{e-z19}
  \text{the inequality~\eqref{e-z20} below is satisfied.}
\end{equation} 
The goal is to solve 
the RHP for $m^{(3)}$ explicitly in this small region.

Let $u=\frac12(z-z^{-1})$ in $\mathcal{O}$.
\begin{figure}[ht]
 \centerline{\epsfig{file=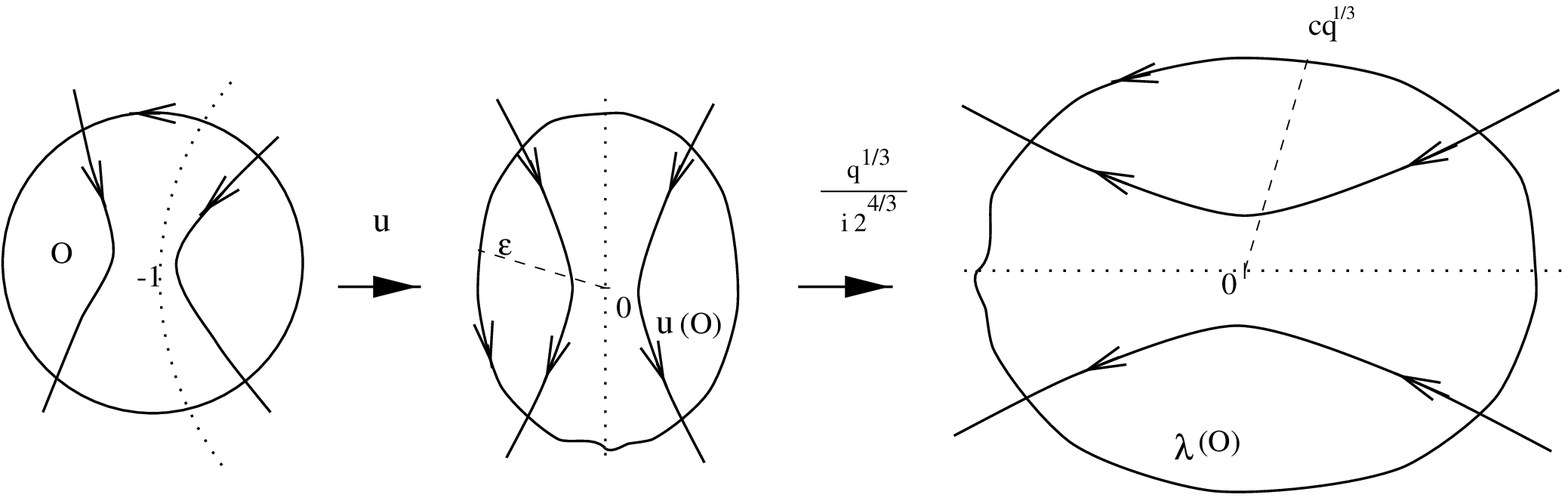, width=11cm}}
 \caption{}\label{fig-map}
\end{figure}
As noted above, we choose and fix $\epsilon>0$ sufficiently small 
(in fact, any number $0<\epsilon<1$ would do.) so that $z\to u(z)$ 
is a bijection 
from $\mathcal{O}$ onto some open neighborhood of $0$ in the $u$-plane : 
under the bijection, 
$\Sigma\cap\mathcal{O}$ becomes a part of
the imaginary axis.
Set $$\lambda(z)=\frac{q^{1/3}u(z)}{i2^{4/3}}=
\frac{q^{1/3}}{i2^{4/3}}\frac12(z-z^{-1}).$$
Note that with $\epsilon$ fixed, there are constants $c_1,c_2>0$ 
such that  
\begin{equation}\label{e-z13}
  c_1 q^{1/3}\le |\lambda(z)| \le c_2 q^{1/3},  
\end{equation}
for all $z\in\partial\mathcal{O}$.
Under the map $z\to\lambda(z)$, $\Sigma\cap\mathcal{O}$ now becomes 
part of the real axis  and 
\begin{equation}
   \lambda(\Sigma^{(3)}\cap\mathcal{O}) = \{ x+iy : 
y^2=\frac{\frac{q^{2/3}}{2^{8/3}}(1-\gamma^2)+x^2}
{1+\frac{2^{8/3}\gamma^2 x^2}{q^{2/3}}}, |x|\le cq^{1/3} \},
\end{equation}
where $c$ is a fixed small number.
As $\frac{q^{2/3}}{2^{8/3}}(1-\gamma^2)=\frac{t}4\bigl(
1-\frac{t}{2^{4/3}q^{2/3}} \bigr)$, $0\le t\le M_2$, we see that 
the contour $\lambda(\Sigma^{(3)}\cap\mathcal{O})$ makes an angle 
$\le \pi/4$ and uniformly bounded away from zero as $q\to\infty$, 
hence has the general shape of the contour in 
Figure~\ref{fig-PII2}, Section~\ref{s-PII}, 
within the ball $\lambda(\mathcal{O})$.
We \emph{define} 
\begin{equation}
   \Sigma^{PII,2}\cap\lambda(\mathcal{O}) \equiv 
\lambda(\Sigma^{(3)}\cap\mathcal{O})
\end{equation}
and extend $\Sigma^{PII,2}$ smoothly outside $\lambda(\mathcal{O})$ 
in such a way that it is asymptotic to straight lines making angles 
between 0 and $\pi/3$ with the real axis.
It is clear from the estimation in Section~\ref{s-PII}, and the 
preceding calculations, that for such a contour $\Sigma^{PII,2}$, 
the bound~\eqref{e-unif2} for the solution $m^{PII,2}(z,t)$ 
of $(\Sigma^{PII,2},v_t^{PII,2})$
\begin{equation}\label{e-z15}
   \text{is uniform for $\gamma, q$ satisfying the relation 
$1-\frac{M_2}{2^{1/3}q^{2/3}}\le \gamma\le 1$.}
\end{equation}

Introduce the parametrix around $z=-1$ as follows.
Define 
\begin{equation*}
  \begin{cases}
    m_p(z)=m^{PII,2}(\lambda(z),t) \quad
\text{in}\quad\mathcal{O}-\Sigma^{(3)}, \\
    m_p(z)=I \quad\text{in}\quad\mathcal{\bar{O}}^c-\Sigma^{(3)}.
  \end{cases}
\end{equation*}
As $q\to\infty$, $|\lambda(z)|\to\infty$ for $z\in\partial\mathcal{O}$,  
and we have for $v_p(z)\equiv v_t^{PII,2}(\lambda(z))$, 
\begin{equation}
 \begin{cases}
   m_p(z) \ \ \text{is analytic in} \ \ 
\C-(\Sigma^{(3)}\cup\partial\mathcal{O})\\
   m_{p+}(z)=m_{p-}(z)v_p(z) \ \ 
\text{on} \ \ \mathcal{O}\cap\Sigma^{(3)}\\
   m_{p+}(z)=m_{p-}(z)I \ \ \text{on} \ \ \mathcal{O}^c\cap\Sigma^{(3)}\\
   m_{p+}(z)=I+\frac{m^{PII,2}_{1}(t)}{\lambda(z)}+O(\frac{1}{\lambda(z)^2}) 
\ \ \text{on} \ \ \partial\mathcal{O} \ \  
\text{as} \ \ q\rightarrow\infty.
 \end{cases}
\end{equation}

The key fact is that $v_p$ is an approximation to  
$v^{(3)}$ with error of order $\frac{1}{q^{2/3}}$. 
We compare, for example, the 12-components of $v^{(3)}$ and 
$v_p$ on $\Sigma^{(3)}$.
We focus on $\mathcal{O}\cap\Sigma_{in}^{(3)}$.
Using the $u$ variable, the $12$-entries of $v^{(3)}$ and $v_p$ are 
$$-\exp(q[\gamma u+\log(\sqrt{1+u^2}-u)])$$ 
and $$-\exp(q[\gamma u-u+\frac16u^3])$$
respectively.
By~\eqref{e-z11} and~\eqref{e-z12}, we have for 
$z\in\mathcal{O}\cap\Sigma_{in}^{(3)}$, 
\begin{equation*}
  |e^{q[\gamma u+\log(\sqrt{1+u^2}-u)]}|
=|z^qe^{\frac{q\gamma}2(z-z^{-1})}|\le
e^{-\frac{2\sqrt{2}q}{3}(1+\gamma\cos\theta)^{3/2}}.
\end{equation*}
From the Taylor expansion of the odd function $\log(\sqrt{1+u^2}-u)$, 
$$\log(\sqrt{1+u^2}-u)=-u+\frac16u^3+u^5r(u),$$
where $r(0)=-\frac3{40}$ and $r(u)$ is bounded for small $u$, 
say $|u|\le\frac12$.
Set $\hat{c}=\sup_{|u|\le\frac12} |r(u)|$.
Note that for $z=\rho_{\theta}e^{i\theta}
\in\mathcal{O}\cap\Sigma^{(3)}_{in}$ (see~\eqref{e-z9}), 
as $q\to\infty$, we have 
$|u|\leq \wt{c}(1+\gamma \cos \theta)^{1/2}\le \wt{c}(c'\epsilon)^{1/2}$.
Therefore, if we have chosen $\epsilon>0$ small enough 
so that 
\begin{equation}\label{e-z20}
  -\frac{2\sqrt{2}}3+\hat{c}(\wt{c})^5c'\epsilon\le -\frac12, 
\end{equation}
we obtain, as $q\to\infty$,
\begin{equation}\label{e-z14}
 \begin{split}
   &|e^{q[\gamma u+\log(\sqrt{1+u^2}-u)]}-e^{q[\gamma u-u+\frac16u^3]}| \\
   &\qquad = |e^{q[\gamma u+\log(\sqrt{1+u^2}-u)]}||1-e^{-qu^5r(u)}| \\
   &\qquad \leq e^{-\frac{2\sqrt{2}q}3(1+\gamma \cos \theta)^{3/2}}\times 
q|u|^5\|r\|_{L^\infty(\{|u|\le\frac12\})}
e^{q\hat{c}(\wt{c})^5(1+\gamma \cos \theta)^{5/2}} \\
   &\qquad \leq Cq(1+\gamma \cos \theta)^{5/2}
e^{q(-\frac{2\sqrt{2}}3+\hat{c}(\wt{c})^5c'\epsilon)(1+\gamma \cos \theta)^{3/2}} \\
   &\qquad \leq Cq(1+\gamma \cos \theta)^{5/2} 
e^{-\frac{q}2(1+\gamma \cos \theta)^{3/2}} \\
   &\qquad \leq \frac{C}{q^{2/3}},
 \end{split}
\end{equation}
where we have used the basic inequality 
$|1-e^z|\le |z|e^{|z|}$ and the fact that 
$\|x^{5/3}e^{-x^{3/2}}\|_{L^{\infty}_{(0,\infty)}} \le C$.
Since 
\begin{equation*}
   v^{(3)}v_p^{-1}= 
   \begin{pmatrix} 1& -e^{q[\gamma u+\log(\sqrt{1+u^2}-u)]}+
e^{q[\gamma u-u+\frac16u^3]}\\ 0&1 
   \end{pmatrix} \ \ \text{on} \ \ \mathcal{O}\cap\Sigma_{in}^{(3)},
\end{equation*}
we have $$\|v^{(3)}v_p^{-1}-I\|_{L^{\infty}
(\mathcal{O}\cap\Sigma^{(3)}_{in})}=O(\frac1{q^{2/3}}).$$
For $\mathcal{O}\cap\Sigma_{out}^{(3)}$, we have a similar estimation.
On the other hand, for 
$\mathcal{O}^c$, the error is exponentially small ; 
$\|v^{(3)}v_p^{-1}-I\|_{L^\infty(\mathcal{O}^c\cap\Sigma^{(3)})}
=\|v^{(3)}-I\|_{L^\infty(\mathcal{O}^c\cap\Sigma^{(3)})}=O(e^{-cq})$.

Now define $R(z)=m^{(3)}m_p^{-1}$.
The ratio is analytic in $\C-(\Sigma^{(3)}\cup\partial\mathcal{O})$ 
and the above calculations show that the jump matrix 
$v_R=m_{p-}v^{(3)}v_p^{-1}m_{p-}^{-1}$ satisfies 
\begin{equation}\label{e-z16}
  \begin{cases}
    \|v_R-I\|_\infty\le\frac{C(M_2)}{q^{2/3}}
&\text{on $\mathcal{O} \cap \Sigma^{(3)}$},\\
    \|v_R-I\|_\infty\le Ce^{-cq} 
&\text{on $\mathcal{O}^c \cap \Sigma^{(3)}$},\\
    v_R=v_{p}^{-1}=m_{p+}^{-1}=I-\frac{m^{PII,2}_1(t)}{\lambda(z)}+
O_{M_2}(\frac{1}{\lambda(z)^2}) 
&\text{on $\partial \mathcal{O}$, as $q\to\infty$.}
  \end{cases}
\end{equation}
In~\eqref{e-z16}, we have used the fact that $m^{PII,2}(z,t)$, 
and hence $m_{p}$, is invertible 
and bounded for $(z,t)\in\C\times [0,M_2]$. (see,~\eqref{e-z15}).

From~\eqref{e-z16} and~\eqref{e-z13}, we see that 
$\|v_R-I\|_\infty=\|w_R\|_\infty\le\frac{C(M_2)}{q^{1/3}}$.
In particular, $(I-C_{w_R})$ is invertible for large $q$ and 
by~\eqref{e-sam5}, $R$ is given by
$$R(z)=I+\frac1{2\pi i} \int_{\Sigma^{(3)}\cup\partial\mathcal{O}} 
\frac{\mu(s)(v_R-I)}{s-z} \,ds$$
where $\mu$ solves $(I-C_{w_R})\mu=I.$
As $\|v_R-I\|_\infty\le\frac{C(M_2)}{q^{1/3}}$, we have 
$\|\mu-I\|_{L^2}=O_{M_2}(\frac1{q^{1/3}})$, and also   
$$R_{22}(0)=1+\frac1{2\pi i} \int_{\Sigma^{(3)}\cup\partial\mathcal{O}} 
(v_R-1)_{22}(s) \frac{ds}{s} 
+O_{M_2}(\frac1{q^{2/3}}).$$
Thus, using $\|v_R-I\|_\infty\le\frac{C(M_2)}{q^{2/3}}$ in~\eqref{e-z16} 
for the second equality, and  $m^{PII,2}_1=m^{PII,1}_1$ 
(see ~\eqref{e-RHPm2'}) for the last equality, we obtain 
\begin{equation}\label{e-R22K1}
  \begin{split}
    \kappa^2_{q-1}&=R_{22}(0)\\
    &= 1+\frac1{2\pi i} \int_{\partial\mathcal{O}} (v_R-1)_{22}(z) \frac{dz}{z} 
+O_{M_2}(\frac1{q^{2/3}})\\
    &=1-\frac1{2\pi i} \int_{\partial\mathcal{O}} 
\frac{m^{PII,2}_{1,22}(t)}{\lambda(z)} \frac{dz}{z} +
O_{M_2}(\frac1{q^{2/3}})\\
    &=1-\frac{m^{PII,2}_{1,22}(t)}{2\pi i} \int_{u(\partial\mathcal{O})} 
\frac{1}{-i\frac{q^{1/3}}{2^{4/3}}u} \,\frac{du}
{(-\sqrt{u^2+1})} +O_{M_2}(\frac1{q^{2/3}})\\
    &=1+\frac{i2^{4/3}}{q^{1/3}}m^{PII,2}_{1,22}(t)+O_{M_2}(\frac1{q^{2/3}})\\
    &=1+\frac{i2^{4/3}}{q^{1/3}}m^{PII,1}_{1,22}(t)+O_{M_2}(\frac1{q^{2/3}}).
  \end{split}
\end{equation}
Note that error in~\eqref{e-R22K1} is uniform for $0\le t\le M_2.$

We summarize as follows.
\begin{lem}\label{lem-K1}
  Let $M_1>0$ be a fixed number which is sufficiently large so 
that~\eqref{e-sec3.6} is satisfied. 
Also let $M_2>0$ and  
$0<\delta_1 <1$ be fixed numbers.
As $q\to \infty$, we have the following results.
\begin{enumerate}
    \item If $0\le \gamma \le 1-\delta_1$,  
then, for some constants $C,c$ which may depend on $\delta_1$,
\begin{equation*}
  |\kappa^2_{q-1}- 1 |\le Ce^{-cq}.
\end{equation*}
    \item If $\frac12 \le \gamma \le 1-\frac{M_1}{2^{1/3}q^{2/3}}$, 
then, 
for some constant $C$ which is independent of $M_1$ satisfying~~\eqref{e-sec3.6}, 
\begin{equation*}
  |\kappa^2_{q-1}-1| \le \frac{C}{q^{1/3}}
e^{-\frac{2\sqrt{2}}3q(1-\gamma)^{3/2}}.
\end{equation*}
    \item If $1-\frac{M_2}{2^{1/3}q^{2/3}} \le \gamma \le 1$, 
\begin{equation*}
  \bigl|\kappa^2_{q-1}-1-\frac{i2^{4/3}}{q^{1/3}}m^{PII}_{1,22}(t)\bigr|
\le \frac{C(M_2)}{q^{2/3}}, 
\end{equation*}
where $t$ is defined by $\gamma=1-\frac{t}{2^{1/3}q^{2/3}}$.
\end{enumerate}
\end{lem}

%%%%=================sec4.tex

\bigskip
\section{$\mathbf{\gamma > 1}$}
\label{s-K2}

Let $\theta_c$ be given in Lemma ~\ref{lem-ourg}, 
$\sin^2\frac{\theta_c}{2}=\frac{1}{\gamma}, 0< \theta_c<\pi$.
Decompose $\Sigma=\overline{C_1}\cup C_2$ where 
$C_1=\{ e^{i\theta} : \theta_c < |\theta | \le \pi \}$
and $C_2=\Sigma-\overline{C_1}$.
Note that on the support of the measure $d\mu_V$ in~\eqref{e-lem4.3.2}, 
$d\mu_V(\theta)=\frac{\gamma}{\pi} \cos(\frac{\theta }{2}) 
\sqrt{\frac{1}{\gamma} - \sin^2(\frac{\theta}{2})} 
\chi_{[-\theta_c, \theta_c]}(\theta) d\theta
=\frac\gamma{4\pi i} \frac{s+1}{s^2}\sqrt{(s-\xi)(s - \xi^{-1})}ds$ 
for $s=e^{i\theta}$.
\begin{figure}[ht]
 \centerline{\epsfig{file=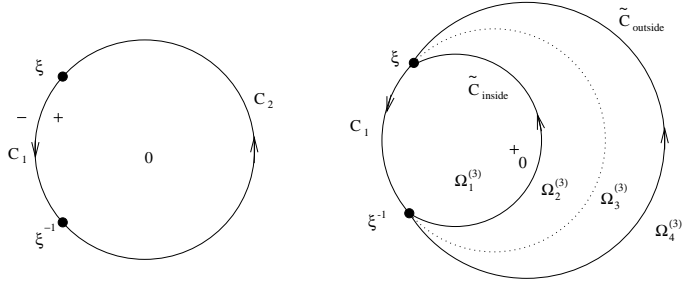, width=9cm}}
 \caption{$\Sigma$ and $\Sigma^{(3)}$}\label{fig-Cinout}
\end{figure}

\begin{lem}\label{lem-alpha}
Define 
$\alpha(z) = -\frac{\gamma}4 \int^z_\xi  
\frac{s+1}{s^2}\sqrt{(s-\xi)(s - \xi^{-1})}ds$ 
where $\xi=e^{i\theta_c}$ and the branch is chosen to 
be analytic in $\C-\overline{C_1}$
and $\sqrt{(s-\xi)(s - \xi^{-1})}>0$ for real $s>0$. Then
\begin{enumerate}
\item $e^{2\alpha}$ is independent of the path in $\C-(\bar{C_1}\cup\{0\})$,.
\item $\exp\bigl( -2\pi i 
\int_{\theta_c}^\phi d\mu_V(\theta)\bigr) 
= \exp\bigl( 2\alpha(z)\bigr) \quad \text{for} \quad z=e^{i\phi}, 
|\phi|<\theta_c$.
\item $\exp\bigl( 2\int_{-\pi}^{\pi}\log |z-e^{i\theta}| d\mu_V(\theta) 
-V(z)+l\bigr)  = \exp\bigl( -2\alpha_-(z)\bigr) \quad \text{for} 
\quad z=e^{i\phi}, |\phi|>\theta_c$.
\end{enumerate}
\end{lem}
\begin{proof}
   Property (i) follows from a standard reside calculation 
: the change in $\alpha(z)$ around
the point at $0$ is $-\pi i$, and the change in $\alpha(z)$
around $C_1$ is $0$. Property (ii) follows from the definition of $\alpha(z)$.
For (iii), set 
\begin{equation*}
   F(\phi)=2\int_{-\pi}^{\pi}\log |e^{i\phi}-e^{i\theta}| d\mu_V(\theta)
+\gamma\cos{\phi} +l +2\alpha_-(e^{i\phi})
\end{equation*}
for $z=e^{i\phi}$, $|\phi|>\theta_c$.
From the variational condition ~\eqref{e-varicon}, we have $F(\theta_c)=0$.
Differentiating,
\begin{equation*}
  \begin{split}
   F'(\phi)&=\int_{-\theta_c}^{\theta_c} i\bigl[ 
\frac{2e^{i\phi}}{e^{i\phi}-e^{i\theta}}-1 \bigr] d\mu_V(\theta)
-\frac{\gamma}{2i} \bigl( e^{i\phi}-e^{-i\phi} \bigr)\\
   &\qquad -\frac{i\gamma}2\frac{e^{i\phi}+1}{e^{i\phi}}
\sqrt{(e^{i\phi}-\xi)(e^{i\phi} - \xi^{-1})}_-\\
  &= \frac{\gamma}{2\pi} \int_{C_2} \frac{z}{z-s} \frac{s+1}{s^2} 
\sqrt{(s-\xi)(s - \xi^{-1})} ds\\
   &\qquad -i-\frac{\gamma}{2i}(z-z^{-1}) -\frac{i\gamma}2\frac{z+1}{z}
\sqrt{(z-\xi)(z - \xi^{-1})}_-.
  \end{split}
\end{equation*}
A residue calculation similar to that in (i), now shows that $F'(\phi)=0$.
Therefore we have $F(\phi)\equiv 0$.
\end{proof}

Note that $e^{2q\pi i\int_{\phi}^{\pi} d\mu_V(\theta)}=1$ for $\phi$ 
outside the support of $d\mu_V$ , i.e.~for $|\phi|>\theta_c$. 
By~\eqref{e-Omega} and above Lemma, our RHP becomes 
\begin{equation}
\begin{cases}
m^{(1)}(z) \ \ \text{is analytic in} \ \ \C-\Sigma, \\
m^{(1)}_+=m^{(1)}_- \begin{pmatrix} e^{-2q\alpha}& (-1)^q\\ 0& e^{2q\alpha} \end{pmatrix} \ \ \text{on} \ \ C_2,\\
m^{(1)}_+=m^{(1)}_- \begin{pmatrix} 1& (-1)^qe^{-2q\alpha_-}\\ 0&1 \end{pmatrix} \ \ \text{on} \ \ C_1,\\
m^{(1)}=I+O(\frac1z) \ \ \text{as}\ \  z\rightarrow \infty
\end{cases}
\end{equation}
and $\kappa^2_{q-1}=-(-1)^qm^{(1)}_{21}(0)e^{ql}=
-(-1)^qe^{q(-\gamma+\log \gamma+1)}m^{(1)}_{21}(0)$
by~\eqref{e-Korig} and~\eqref{e-lem4.3.3}.

We use the same conjugation ~\eqref{e-conju} 
for $m^{(1)}$ as in the case $\gamma\le 1$.
Then our new jump matrices for $m^{(2)}$ are 
\begin{equation}\label{e-m2tilde}
\begin{cases}
v^{(2)}= \begin{pmatrix} 1& -e^{-2q\alpha}\\ e^{2q\alpha}&0 \end{pmatrix} \ \ \text{on} \ \ C_2\\
v^{(2)}= \begin{pmatrix} e^{-2q\alpha_-}& -1\\ 1&0 \end{pmatrix} \ \ \text{on} \ \ C_1
\end{cases}
\end{equation}
and $\kappa^2_{q-1}=e^{q(-\gamma+\log \gamma+1)}m^{(2)}_{22}(0)$.

Set $\Sigma^{(3)}=\overline{C_1}\cup \wt{C}_{inside} \cup \wt{C}_{outside}$ 
where $\wt{C}_{inside}$ and $\wt{C}_{outside}$ are open arcs as  chosen below.
Note the factorization $v^{(2)}=
\begin{pmatrix} 1&0\\e^{2q\alpha}&1 \end{pmatrix}
\begin{pmatrix} 1& -e^{-2q\alpha}\\0&1 \end{pmatrix}$on $C_2$. 
Set $Re \alpha =R$, $Im \alpha =I$ so that $\alpha =R+iI$.
Recall the Cauchy-Riemann equations in polar coordinates $(r,\theta)$,
\begin{equation*}
   r\frac{\partial R}{\partial r} =\frac{\partial I}{\partial\theta},
\quad r\frac{\partial I}{\partial r} =-\frac{\partial R}{\partial\theta}. 
\end{equation*}
For $z=e^{i\theta}\in C_2$, 
$\alpha(z)=-\pi i \int_{\theta_c}^{\theta} d\mu_V(\theta')$ 
is pure imaginary and 
\begin{equation*}
\frac{\partial I}{\partial\theta}
=\frac{\partial}{\partial\theta} (-i\alpha)=-\pi\frac{\gamma}{\pi}
\cos(\frac{\theta }{2})
\sqrt{\frac{1}{\gamma} - \sin^2(\frac{\theta}{2})} <0.
\end{equation*}
Hence 
\begin{equation*}
   R=0  \quad\text{and}\quad
   \frac{\partial R}{\partial r} =\frac{\partial I}{\partial\theta} <0
\quad\text{on}\quad C_2. 
\end{equation*}
Therefore for fixed $\theta$, $e^{i\theta}\in C_2$, there is 
$\epsilon_1=\epsilon_1(\theta) >0$ 
such that $R=Re \alpha >0$ (resp. $<0$) for $z=re^{i\theta}$ with 
$1-\epsilon_1 < r< 1$ (resp. $1<r<1+\epsilon_1$).
We take $\wt{C}_{inside}$ (resp, $\wt{C}_{outside}$) 
such that $|e^{-2\alpha}|<1$ (resp, $|e^{2\alpha}|<1$) 
on $\wt{C}_{inside}$ (resp, $\wt{C}_{outside}$).
Clearly there exist $0<\rho_1, \rho_2 <1$ such that $|e^{-2\alpha}|<\rho_1$ 
(resp, $|e^{2\alpha}|<\rho_2$), for all $z\in\wt{C}_{inside}$ 
(resp, $z\in\wt{C}_{outside}$), apart from a small neighborhood 
of the endpoints.
Introduce the regions $\Omega^{(3)}_{k}$, $k=1,2,3,4$ 
as in Figure~\ref{fig-Cinout}.
Define $m^{(3)}$ as follows, 
\begin{equation*}
  \begin{cases}
     &m^{(3)}=m^{(2)} \begin{pmatrix} 1&-e^{-2q\alpha}\\0&1 \end{pmatrix}^{-1}
\quad\text{in}\quad \Omega^{(3)}_2,\\
     &m^{(3)}=m^{(2)} \begin{pmatrix} 1&0\\e^{2q\alpha}&1 \end{pmatrix}
\qquad\text{in}\quad \Omega^{(3)}_3,\\
     &m^{(3)}=m^{(2)} \qquad\qquad\text{in}\quad \Omega^{(3)}_1,\Omega^{(3)}_4.
  \end{cases} 
\end{equation*}
Then $v^{(3)}$ is given by
\begin{equation*}
  \begin{cases}
     &\begin{pmatrix} 1&-e^{-2q\alpha}\\0&1 \end{pmatrix}
\quad\text{on}\quad \wt{C}_{inside},\\
     &\begin{pmatrix} 1&0\\e^{2q\alpha}&1 \end{pmatrix}
\qquad\text{on}\quad \wt{C}_{outside},\\
     &\begin{pmatrix} e^{-2q\alpha_-}&-1\\1&0 \end{pmatrix}
\quad\text{on}\quad C_1,
  \end{cases}
\end{equation*}
and 
\begin{equation}\label{e-Km3K2}
   \kappa^2_{q-1}=e^{q(-\gamma+\log \gamma+1)}m^{(3)}_{22}(0).
\end{equation}

From Lemma ~\ref{lem-alpha} (iii) and the second variational 
condition in~\eqref{e-varicon}, we have, for any $z\in C_1$,
\begin{equation*}
   e^{-2q\alpha_-}\to 0 \quad\text{as} \quad q\to\infty.
\end{equation*}
(Recall that the inequality in ~\eqref{e-varicon} is strict
from Lemma ~\ref{lem-ourg} (ii)).
Also from the choice of $\wt{C}_{inside}$ and $\wt{C}_{outside}$, 
\begin{equation*}
   e^{-2q\alpha}\to 0, \quad e^{2q\alpha}\to 0
\quad\text{as} \quad q\to\infty
\end{equation*}
on $\wt{C}_{inside}$, $\wt{C}_{outside}$, respectively.
Therefore $v^{(3)} \rightarrow v^\infty$ as $q\rightarrow \infty$,
\begin{equation}
\begin{cases}
v^\infty =
 \begin{pmatrix} 1&0\\0&1 \end{pmatrix} \ \ \text{on} \ \ 
\wt{C}_{inside}\cup\wt{C}_{outside},\\
v^\infty =
 \begin{pmatrix} 0&-1\\1&0 \end{pmatrix} \ \ \text{on} \ \ 
C_1.
\end{cases}
\end{equation}

The following result can be verified by direct calculation.
\begin{lem}\label{lem-minfty}
RHP $(C_1, v^\infty)$ can be solved explicitly.
\begin{equation}
m^\infty = \begin{pmatrix} \frac12(\beta+\beta^{-1})& \frac1{2i}(\beta-\beta^{-1})\\
                           -\frac1{2i}(\beta-\beta^{-1})& \frac12(\beta+\beta^{-1}) \end{pmatrix}
\end{equation}
where $\beta(z)\equiv (\frac{z-\xi}{z-\xi^{-1}})^{1/4}, \ \ 
\text{is analytic in} \ \ \C-\bar{C}_1 \ \ \text{such that} \ \ 
\beta\sim +1 \ \ \text{as} \ \ z\rightarrow \infty$
\end{lem}
Hence we expect from~\eqref{e-Km3K2}, 
$\kappa^2_{q-1}\sim e^{q(-\gamma+\log \gamma+1)}\frac1{\sqrt{\gamma}},$ 
because $m^\infty_{22}(0)=\frac1{\sqrt{\gamma}}$.

Our goal now is to show that indeed $m^{(3)}\to m^{\infty}$ 
as $q\to\infty$. As in Section ~\ref{s-K1}, we must control 
the behavior of the solution of the RHP for $m^{(3)}$ near the 
endpoints, where the rate of exponential convergence 
$v^{(3)}\to v^{\infty}$, becomes smaller and smaller.

%%%%%===================================sec5.tex

\bigskip
Let $\delta_3, M_4>0$ be fixed numbers, let $0< \delta_4<1$ 
be a fixed, sufficiently small number satisfying~\eqref{e-z30} below, 
and let $M_3 >0$ be a fixed, sufficiently large number 
satisfying~\eqref{e-z29} and~\eqref{e-z28} below.
We consider 3 cases for $\gamma$ : 
\begin{enumerate}
    \item $1+\delta_3 \le \gamma$
    \item $1+\frac{M_3}{2^{1/3}q^{2/3}} \le \gamma \le 1+\delta_4$
    \item $1\le \gamma \le 1+\frac{M_4}{2^{1/3}q^{2/3}}$
\end{enumerate}

\bigskip
Calculations similar to those that are needed for the asymptotics 
of the orthogonal polynomial on the real line (see, \cite{DKMVZ}), 
show that for $\gamma>1+\delta_3$,
\begin{equation}\label{8.1}
  \kappa^2_{q-1}= e^{q(-\gamma+\log \gamma+1)}\frac1{\sqrt{\gamma}}(1+O(\frac1{q})).
\end{equation}
The error is uniform for $1+\delta_3 \le \gamma \le L$ for any fixed $L<\infty$.
However, we will not use this result, utilizing instead (stronger) estimates 
from \cite{Johan} (see next section).

\bigskip
\begin{figure}[ht]
 \centerline{\epsfig{file=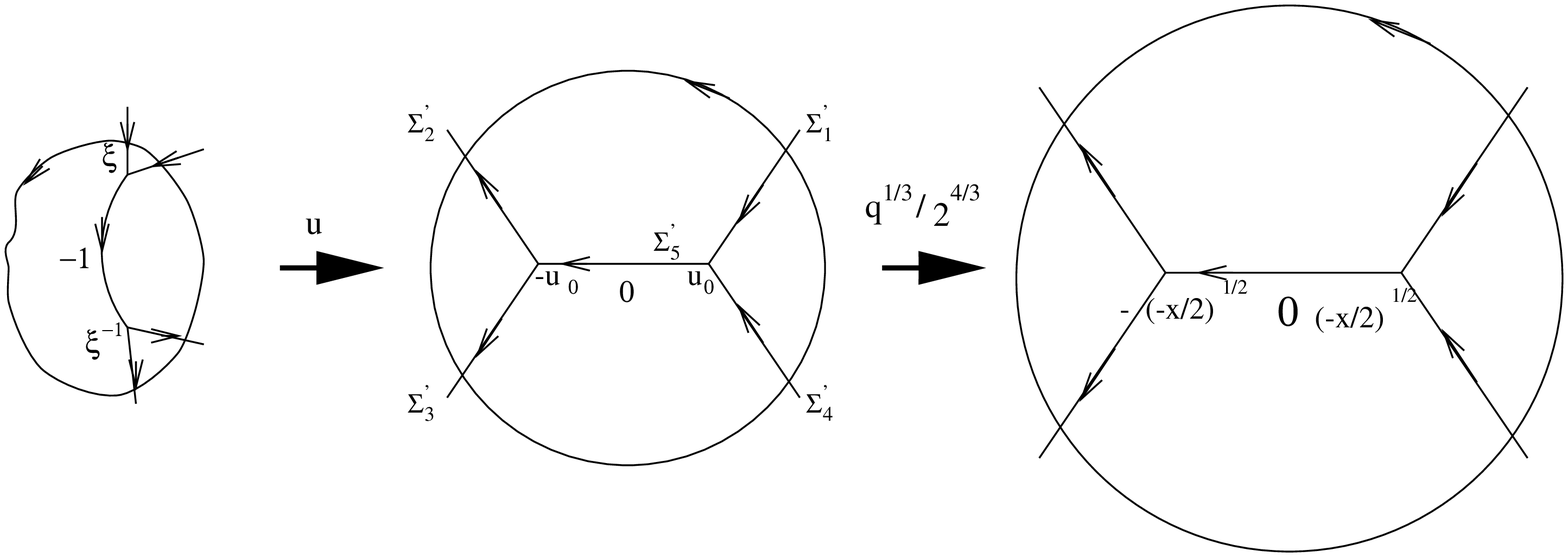, height=4cm}}
 \caption{}\label{fig-umap2}
\end{figure}
We consider case (iii).
Set $\gamma =1+\frac{t}{2^{1/3}q^{2/3}}$ with $0\le t \le M_4$ and $u_0=\sin{\theta_c}=\frac2{\gamma} \sqrt{\gamma-1}$.
In defining $\Sigma^{(3)}$ above, there is some freedom in the choice 
of $\wt{C}_{inside}$ and $\wt{C}_{outside}$.
We make the following choice (see~\eqref{e-z23} below). 
Set $x=-\frac{q^{2/3}u_0^2}{2^{5/3}}=
-t(1+\frac{t}{2^{1/3}q^{2/3}})^{-2} \sim -t<0$ as $q\to\infty$,  
and let $\Sigma^{PII,3}$ be the contour defined in Figure~\ref{fig-PII3} 
for this specific $x$.
Let $\Sigma'=\{ u=\frac{2^{4/3}}{q^{1/3}}\lambda : \lambda\in\Sigma^{PII,3} \}
=\cup_{k=1}^5 \Sigma'_k$, 
and let $\epsilon' >0$ be small and fixed (see~\eqref{e-z21},~\eqref{e-z22} below).
For definiteness, we can, and do assume that the rays 
$\Sigma'_1,\dots,\Sigma'_4$ make an angle of $\pi/6$ with the real axis.
Consider $\mathcal{O}'=\{ u : |u|<\epsilon' \}$.
If $q$ is large enough, then $u_0\in\mathcal{O}'$.
Set $u=u(z)=\frac1{2i}(z-z^{-1})$. 
We choose $\epsilon'$ such that (cf.~\eqref{e-z18}) 
\begin{equation}\label{e-z21}
   \text{$u$ is a bijection from an open neighborhood of $z=-1$ onto $\mathcal{O}'$.}
\end{equation}
Clearly there are constants $c_1,c_2>0$ such that 
$c_1\le |z(u)|\le c_2$ for all $u\in\partial\mathcal{O}'$.

Under $u^{-1}$, the points $u_0, -u_0$ are mapped into $\xi, \bar{\xi}$ respectively,  
and $u^{-1}(\Sigma'_5)=C_1$. 
Consider a point $z\in u^{-1}(\Sigma'_4 \cap\mathcal{O}')$, the inverse image of 
a point $u\in\Sigma'_4 \cap\mathcal{O}$.
Changing variables twice, 
$v=\frac1{2i}(s-s^{-1})$ and $w=v^2$, 
\begin{equation}\label{8.2}
  \begin{split}
    -2\alpha(z)
    &= \frac{i\gamma}{2} \int^z_\xi \bigl(\sqrt{s}+\sqrt{s}^{-1}\bigr)
\sqrt{(s+s^{-1})-(\xi+\xi^{-1})} \frac{ds}{is}\\
    &= -i\gamma \int_{u_0}^u \biggl( \sqrt{1-u_0^2}-\sqrt{1-v^2}\biggr)^{1/2}
\frac{\sqrt{1-\sqrt{1-v^2}}}{\sqrt{1-v^2}} dv\\
    &= \frac{-i\gamma}2 \int_{u_0}^u \biggl( \frac{v^2-u_0^2}
{1+k(u_0^2)+k(v^2)}\biggr)^{1/2} v\bigl( 1+h(v^2)\bigr) dv\\
    &= -\frac{i\gamma}4 \int_{u_0^2}^{u^2} (w-u_0^2)^{1/2}
\bigl[ 1+(h(w)-h(u_0^2))+h(u_0^2)\bigr] dw\\
    &\qquad -\frac{i\gamma}4 \int_{u_0^2}^{u^2} (w-u_0^2)^{1/2}
\frac{\bigl[ (k(u_0^2)-k(w))-2k(u_0^2)\bigr] (1+h(w))dw}
{(1+k(u_0^2)+k(w))^{1/2}
\bigl[ 1+(1+k(u_0^2)+k(w))^{1/2}\bigr]}\\
    &= \frac{-i\gamma}{6} (u^2-u_0^2)^{3/2}
+O\biggl( |u^2-u_0^2|^{5/2}\biggr)
+O\biggl( |u^2-u_0^2|^{3/2}u_0^2\biggr),
  \end{split}
\end{equation}
where $\sqrt{u^2-u_0^2}$ is defined to be analytic in $\C-[-u_0,u_0]$ and 
positive for real $u>u_0$ ; 
$h(w)=\frac{\sqrt{2}\sqrt{1-\sqrt{1-w}}}{\sqrt{w}\sqrt{1-w}}-1$,
which is analytic in $|w|\le\epsilon'$ and $h(0)=0$ ; 
$k(w)=\frac12(\sqrt{1-w}-1)$, which is also analytic
in $|w|\le\epsilon'$ and $k(0)=0$.
Since $\Sigma'_4$ is a straight ray of angle $\frac{-\pi}{6}$ 
at $u_0$, $Re(-i(u^2-u_0^2)^{3/2}) \le -\frac12|u^2-u_0^2|^{3/2}$, 
which yields 
\begin{equation}\label{e-z31}
   |\exp(-2\alpha(z))|\le 
\exp\bigl(-\frac1{24}|u^2-u_0^2|^{3/2}\bigr) < 1, 
\end{equation} 
provided  $\epsilon'$ is sufficiently small so that 
\begin{equation}\label{e-z22}
   O(|u^2-u_0^2|^{5/2})+O(|u^2-u_0^2|^{3/2}u_0^2)
\le c|u^2-u_0^2|^{3/2}(|u^2-u_0^2|+u_0^2) \le\frac1{24} |u^2-u_0^2|^{3/2}
\end{equation}
for $u\in\Sigma'_4\cap\mathcal{O}'$, where the terms on the LHS 
are given in~\eqref{8.2}.
The same choice of $\epsilon'$ gives rise the same 
result for $z\in u^{-1}(\Sigma'_3 \cap\mathcal{O}')$, and 
also $|\exp(2\alpha(z))|\le \exp(-\frac1{24}|u^2-u_0^2|^{3/2}) < 1$ 
for $z\in u^{-1}((\Sigma'_1\cup\Sigma'_2)\cap\mathcal{O}')$.
We thus fix $\Sigma^{(3)}$ by choosing 
\begin{equation}\label{e-z23}
   \Sigma^{(3)}\equiv u^{-1}(\Sigma'\cap\mathcal{O}') 
\quad\text{inside $\mathcal{O}'$}, 
\end{equation}
and extending it to a contour of the general shape 
$\bar{C}_1\cup\wt{C}_{inside}\cup\wt{C}_{outside}$ 
as in Figure~\ref{fig-Cinout}.

Define 
\begin{equation}
  \begin{cases}
    m_p(z)=m^{PII,3}\bigl(\frac{q^{1/3}}{2^{4/3}}u(z),x\bigr)& 
\quad\text{in} \quad\mathcal{O}-\Sigma^{(3)},\\
    m_p(z)= I& \quad\text{in}\quad\bar{\mathcal{O}}^c-\Sigma^{(3)},
  \end{cases}
\end{equation}
where $m^{PII,3}(z,x)$ solves the RHP of Painlev\'e II equation 
given by~\eqref{e-z17} and~\eqref{e-vPII4}.
Then $m_p$ solves the RHP on $\Sigma^{(3)}\cup\partial\mathcal{O}$ 
in which the jump matrix $v_p(z)$ is given by 
\begin{equation}
  \begin{cases}
    v^{PII,3}(\frac{q^{1/3}}{2^{4/3}}u(z))&\quad 
z\in\Sigma^{(3)}\cap\mathcal{O},\\
    I& \quad z\in\Sigma^{(3)}\cap\mathcal{O}^c,\\
    m_{p+}(z)& \quad z\in\partial\mathcal{O},
  \end{cases}
\end{equation}
where $v^{PII,3}$ is given in~\eqref{e-vPII4}.

We compare 
$v^{(3)}$ and $v_p$.
First, let $z\in \Sigma^{(3)}\cap\mathcal{O}$
such that $u(z)\in\Sigma'_4\cap\mathcal{O}'$. 
The $12$-entries of $v^{(3)}$ and $v_p$ are 
$-\exp(-2q\alpha(z))$ and $-\exp(-2ig^{PII}(\frac{q^{1/3}}{2^{4/3}}u(z)))
=-\exp(-\frac{iq}6(u^2-u_0^2)^{3/2})$, respectively.
With $\epsilon'$ chosen small as above, using 
\begin{align*}
  &Re(-i(u^2-u_0^2)^{3/2})\le -\frac12|u^2-u_0^2|^{3/2},\\
  &u_0^2\le \frac{4M_4}{q^{2/3}},\\
  &\|x^{5/2}e^{-x^{3/2}}\|_{L^{\infty}[0,\infty)} \le C,\\
  &\gamma-1 \le \frac{M_4}{2^{1/3}q^{2/3}},
\end{align*}
we obtain from~\eqref{8.2},
\begin{equation}
  \begin{split}
    &|e^{-2q\alpha(z)} - e^{-2ig^{PII}(\frac{q^{1/3}}{2^{4/3}}u(z))}|\\ 
    &\qquad \le e^{Re(-\frac{iq}6(u^2-u_0^2)^{3/2})} 
|e^{-2q\alpha(z)+2ig^{PII}(\frac{q^{1/3}}{2^{4/3}}u(z))} -1|\\
    &\qquad \le Ce^{-\frac{q}{24}|u^2-u_0^2|^{3/2}} 
\bigl[ q|u^2-u_0^2|^{3/2}(|u^2-u_0^2|+u_0^2+(\gamma-1)) \bigr]\\
    &\qquad \le \frac{C(M_4)}{q^{2/3}}.
  \end{split}
\end{equation}
In a similar manner, for $z$ such that $u(z)\in\Sigma'_3\cap\mathcal{O}'$, 
same result holds and for $z\in\wt{C}_{outside}\cap\mathcal{O}$, 
the difference of the $21$-entries of $v^{(3)}$ and $v_p$ satisfies 
$|e^{2q\alpha(z)} - e^{2ig^{PII}(\frac{q^{1/3}}{2^{4/3}}u(z))}|
\le\frac{C(M_4)}{q^{2/3}}$.
For $z\in C_1$,  
$Re(-i(u^2-u_0^2)_-^{3/2})=-|u^2-u_0^2|^{3/2}$.  
Again by~\eqref{8.2}, 
the difference of the $11$-entries of $v^{(3)}$ and $v_p$ satisfies
$|e^{-2q\alpha_-(z)} - e^{-2ig_-^{PII}(\frac{q^{1/3}}{2^{4/3}}u(z))}|
\le\frac{C(M_4)}{q^{2/3}}$.
Therefore, we have
\begin{equation}\label{8.7}
   \|v^{(3)}v_p^{-1} -I\|_{L^{\infty}(\Sigma^{(3)}\cap\mathcal{O})}
\le\frac{C(M_4)}{q^{2/3}}.
\end{equation}

Secondly, for $z\in\Sigma^{(3)}\cap\mathcal{O}^c$, 
$|z-\xi|,|z-\xi^{-1}|\geq c>0$ implies 
exponential decay for $e^{-2q\alpha(z)}$ and $e^{2q\alpha(z)}$ 
for $z\in\wt{C}_{inside}\cap\mathcal{O}^c$ and 
$z\in\wt{C}_{outside}\cap\mathcal{O}^c$, respectively.
Therefore we have 
\begin{equation}\label{8.8}
   \|v^{(3)}v_p^{-1} -I\|_{L^{\infty}(\Sigma^{(3)}\cap\mathcal{O}^c)} 
\le Ce^{-cq}.
\end{equation}

Finally, for $z\in\partial\mathcal{O}$, 
$|u(z)|=\epsilon'$ and 
\begin{equation}\label{8.9}
   m_{p+}(z)=m^{PII,3}\bigl(\frac{q^{1/3}}{2^{4/3}}u(z),x\bigr)
=I+\frac{m^{PII,3}_1(x)}{\frac{q^{1/3}}{2^{4/3}}u(z)}+ 
O_{M_4}(\frac1{q^{2/3}}),
\end{equation}
by~\eqref{e-z17}.
Here the error is uniformly for $0\le x\le 2M_4$.

Now as in the case $\gamma\le 1$, 
define $R(z)=m^{(3)}m_p^{-1}$. Then the jump matrix for $R$ is given by 
$v_R=m_{p-}v^{(3)}v_p^{-1}m_{p-}^{-1}$. From~\eqref{8.7},~\eqref{8.8} and ~\eqref{8.9}, 
and also~\eqref{e-unif4}, 
we have 
\begin{equation}
  \begin{cases}
    \|v_R-I\|_{L^{\infty}(\Sigma^{(3)}\cap\mathcal{O})} 
\le\frac{C(M_4)}{q^{2/3}}& \text{on $\mathcal{O} \cap \Sigma^{(3)}$},\\
    \|v_R-I\|_{L^{\infty}(\Sigma^{(3)}\cap\mathcal{O}^c)}
\le Ce^{-cq}& \text{on $\mathcal{O}^c \cap \Sigma^{(3)}$},\\
    v_R=m_{p+}^{-1}=I-\frac{2^{4/3}m^{PII,3}_1(x)}{q^{1/3}u(z)}
+O_{M_4}(\frac{1}{q^{2/3}})& 
\text{on $\partial \mathcal{O}$}.
  \end{cases}
\end{equation}

As in ~\eqref{e-R22K1} for the case $\gamma \le 1$, 
using $m_{1,22}^{PII}=m_{1,22}^{PII,3}+(ix^2/8)$
from ~\eqref{e-mPmP4},
\begin{equation}\label{e-z34}
  \begin{split}
    m_{22}^{(3)}(0)&= R_{22}(0)\\
    &= 1+\frac{i2^{4/3}}{q^{1/3}}m^{PII,3}_{1,22}(x)+O_{M_4}(\frac1{q^{2/3}})\\
    &= 1+\frac{i2^{4/3}}{q^{1/3}} \biggl[m^{PII}_{1,22}(x) 
-\frac{it^2}{8}\bigl( 1+\frac{t}{2^{1/3}q^{2/3}}\bigr)^{-4}\biggr] 
+O_{M_4}(\frac1{q^{2/3}})\\
    &= 1+\frac{i2^{4/3}}{q^{1/3}}m^{PII}_{1,22}(x) +
\frac{t^2}{2^{5/3}q^{1/3}} +O_{M_4}(\frac1{q^{2/3}}).
  \end{split}
\end{equation}
Therefore, from ~\eqref{e-Km3K2} 
using  $x=-t(1+\frac{t}{2^{1/3}q^{2/3}})^{-2}=-t+O_{M_4}(\frac1{q^{2/3}})$
and the fact that $\frac{d}{dt}m^{PII}_{1,22}(t)$ is bounded for $0\le t\le M_4$,
(this follows, for example from~\eqref{e-um1PII} and the boundedness of 
$u(x)=2im_{1,12}^{PII}$ ; alternatively statements like~\eqref{e-unif4} are true also 
for all the $x$-derivatives of $m^{PII,3}(z;x)$, etc.), 
\begin{equation}
  \begin{split}
  \kappa^2_{q-1}&=e^{q(-\gamma+\log \gamma+1)}m^{(3)}_{22}(0)\\
  & =e^{q(-\gamma+\log \gamma+1)} \biggl( 1 +\frac{i2^{4/3}}{q^{1/3}}m^{PII}_{1,22}(x) +\frac{t^2}{2^{5/3}q^{1/3}} +O_{M_4}(\frac1{q^{2/3}}) \biggr)\\
  & = \biggl( 1-\frac{t^2}{2^{5/3}q^{1/3}} +O_{M_4}(\frac1{q}) \biggr) 
\biggl( 1 +\frac{i2^{4/3}}{q^{1/3}}m^{PII}_{1,22}(x) 
+\frac{t^2}{2^{5/3}q^{1/3}} +O_{M_4}(\frac1{q^{2/3}}) \biggr) \\
  & =1 +\frac{i2^{4/3}}{q^{1/3}}m^{PII}_{1,22}(x) +O_{M_4}(\frac1{q^{2/3}})\\
  & =1 +\frac{i2^{4/3}}{q^{1/3}}m^{PII}_{1,22}(-t) +O_{M_4}(\frac1{q^{2/3}}).
  \end{split}
\end{equation}

\bigskip
Finally we consider the case (ii), $1+\frac{M_3}{2^{1/3}q^{2/3}}\le 
\gamma\le 1+\delta_4$.
We conjugate $m^{(2)}$ with jump matrix $v^{(2)}$ given 
by~\eqref{e-m2tilde}, as follows.\\
\begin{equation}
 \begin{cases}
    m^{(4)} \equiv m^{(2)} \qquad |z|>1,\\
    m^{(4)} \equiv m^{(2)} \bigl(\begin{smallmatrix} 0& 1\\-1&0 \end{smallmatrix} \bigr)
\quad |z|<1.
 \end{cases}
\end{equation}
Define $\wt{\alpha}(z)= -\pi i \int^z_\xi \frac\gamma{4\pi i} 
\frac{s+1}{s^2}\sqrt{(s-\xi)(s - \xi^{-1})}ds$, 
where $\wt{\alpha}(z)$ is the same as $\alpha(z)$ in Lemma~\ref{lem-alpha}, but 
now we choose the branch so that 
$\sqrt{(s-\xi)(s - \xi^{-1})}$ is analytic in $\C-\bar{C}_2$, and 
$\sqrt{(s-\xi)(s - \xi^{-1})} \sim +s$ as $s\to\infty$.
Then the jump matrix $v^{(4)}$ for $m^{(4)}$ is given by
\begin{equation}
 \begin{cases}
  v^{(4)}= \begin{pmatrix} e^{2q\wt{\alpha}_+} &1\\ 0&e^{2q\wt{\alpha}_-} 
\end{pmatrix} \ \ \text{on} \ \ C_2,\\
  v^{(4)}= \begin{pmatrix} 1&e^{-2q\wt{\alpha}}\\ 0&1 \end{pmatrix} \ \ \text{on} 
\ \ C_1
 \end{cases}
\end{equation}
and $\kappa^2_{q-1}=-e^{q(-\gamma+\log \gamma+1)}m^{(4)}_{21}(0)$.
\begin{figure}[ht]
 \centerline{\epsfig{file=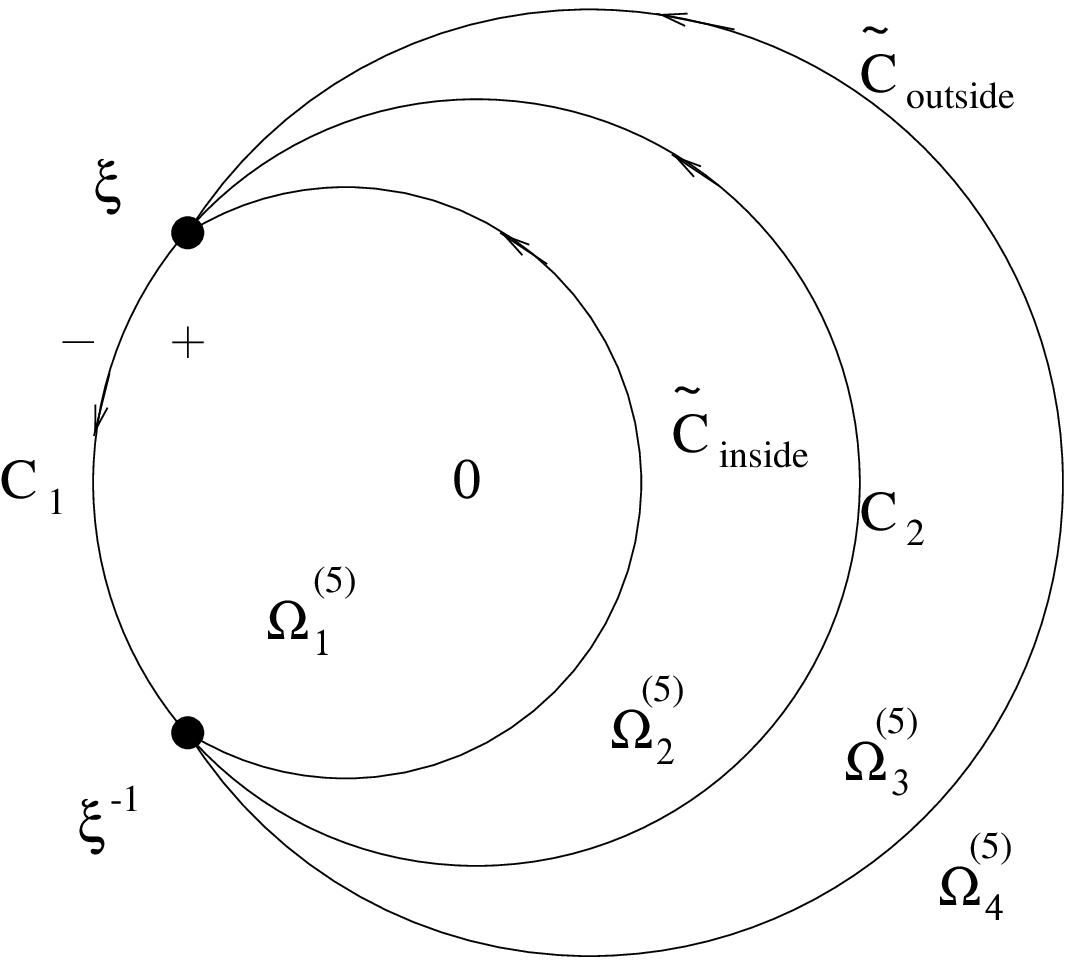, width=6cm}}
 \caption{}\label{fig-Ctilde}
\end{figure}
Noting the factorization $v^{(4)}=\begin{pmatrix} 
1&0\\ e^{2q\wt{\alpha}_-}&1 \end{pmatrix}$
$\begin{pmatrix} 0&1\\-1&0 \end{pmatrix}$
$\begin{pmatrix} 1&0\\e^{2q\wt{\alpha}_+}&1 \end{pmatrix}$ on $C_2$,
we define (see Figure~\ref{fig-Ctilde}) 
\begin{equation*}
  \begin{cases}
     m^{(5)}=m^{(4)} \begin{pmatrix} 1&0\\ e^{2q\wt{\alpha}}&1 \end{pmatrix}^{-1}
\quad\text{in}\quad \Omega^{(5)}_2,\\
     m^{(5)}=m^{(4)} \begin{pmatrix} 1&0\\ e^{2q\wt{\alpha}}&1 \end{pmatrix}
\qquad\text{in}\quad \Omega^{(5)}_3,\\
     m^{(5)}=m^{(4)} 
\qquad\text{in}\quad \Omega^{(5)}_1\cup\Omega^{(5)}_4,
  \end{cases}
\end{equation*}
so that 
\begin{equation*}
  \begin{cases}
    v^{(5)}=\begin{pmatrix} 0&1\\-1&0 \end{pmatrix} \quad\text{on}\quad C_2,\\
    v^{(5)}=\begin{pmatrix} 1&e^{-2q\wt{\alpha}}\\ 0&1 \end{pmatrix} \quad\text{on}\quad C_1,\\
    v^{(5)}=\begin{pmatrix} 1&0\\ e^{2q\wt{\alpha}}&1 \end{pmatrix}
\quad\text{on}\quad \wt{C}_{inside}\cup \wt{C}_{outside}.
  \end{cases}
\end{equation*}

As in the case of $\alpha$, we have $|e^{-\wt{\alpha}(z)}|<1$ for $z\in C_1$ and 
$|e^{\wt{\alpha}(z)}|<1$ for $z\in \wt{C}_{inside}\cup \wt{C}_{outside}$.
Therefore taking $q\to\infty$, we have 
\begin{equation}
\begin{cases}
v^{(5,\infty)} =
 \begin{pmatrix} 0&1\\-1&0 \end{pmatrix} \ \ \text{on} \ \ C_2,\\
v^{(5,\infty)} =
 \begin{pmatrix} 1&0\\0&1 \end{pmatrix} \ \ \text{on} \ \ C_1.
\end{cases}
\end{equation}
This RHP can be solved explicitly as in Lemma ~\ref{lem-minfty}, 
and we find
\begin{equation}\label{e-z24}
m^{(5,\infty)} = \begin{pmatrix} \frac12(\wt{\beta}+\wt{\beta}^{-1})& 
\frac1{2i}(\wt{\beta}-\wt{\beta}^{-1})\\
      -\frac1{2i}(\wt{\beta}-\wt{\beta}^{-1})& 
\frac12(\wt{\beta}+\wt{\beta}^{-1}) \end{pmatrix}
\end{equation}
where $\wt{\beta}(z)\equiv (\frac{z-\xi}{z-\xi^{-1}})^{1/4}$ 
is now analytic in $\C-\bar{C}_2$ and 
$\wt{\beta}\sim +1$ as $z\to\infty$.
From~\eqref{e-z24}, we have $m^{(5,\infty)}_{21}(0)=-\frac1{\sqrt{\gamma}}$ and 
$\kappa^2_{q-1}\sim e^{q(-\gamma+\log \gamma+1)}\frac1{\sqrt{\gamma}}$ 
as $q\to\infty$.
Again, we need to construct parametrices around $\xi$ and $\xi^{-1}$ 
in order to prove 
that indeed $m^{(5)}\to m^{(5,\infty)}$.
Note that $\det m^{(5,\infty)}=1$ : this follows either by direct calculation 
or by a general argument as $\det v^{(5,\infty)}=1$.

\begin{figure}[ht]
 \centerline{\epsfig{file=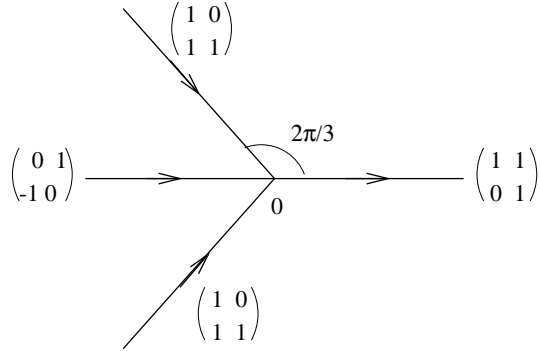, width=7cm}}
 \caption{$\Sigma''$ and $v_{\Psi}$}\label{fig-Sigma2}
\end{figure}

Let $\Sigma''$ be the contour $\Sigma''=\R\cup\R_+e^{2\pi i/3}\cup\R_+e^{4\pi i/3}$ 
shown in Figure~\ref{fig-Sigma2}.
Let $\omega=e^{2\pi i/3}$ and set (see \cite{DZ})
\begin{equation}
  \begin{cases}
        \Psi(s)=\begin{pmatrix}Ai(s)&Ai(\omega^2s)\\
Ai'(s)&\omega^2Ai'(\omega^2s) \end{pmatrix}
e^{-\frac{i\pi}{6}}\sigma_3, &\quad 0< arg s < \frac{2\pi}{3},\\
        \Psi(s)=\begin{pmatrix}Ai(s)&Ai(\omega^2s)\\
Ai'(s)&\omega^2Ai'(\omega^2s) \end{pmatrix}
e^{-\frac{i\pi}{6}\sigma_3}\begin{pmatrix}1&0\\-1&1\end{pmatrix}, 
&\quad \frac{2\pi}{3}< arg s <\pi, \\
        \Psi(s)=\begin{pmatrix}Ai(s)&-\omega^2Ai(\omega s)\\
Ai'(s)&-Ai'(\omega s) \end{pmatrix}
e^{-\frac{i\pi}{6}\sigma_3}\begin{pmatrix}1&0\\1&1\end{pmatrix}, 
&\quad \pi < arg s < \frac{4\pi}{3}, \\
        \Psi(s)=\begin{pmatrix}Ai(s)&-\omega^2Ai(\omega s)\\
Ai'(s)&\omega Ai'(\omega s) \end{pmatrix}
e^{-\frac{i\pi}{6}}\sigma_3, &\quad \frac{4\pi}{3}<arg s <2\pi,\\
\end{cases}
\end{equation}
where $Ai(s)$ is the Airy function.
Then $\Psi$ satisfies the jump conditions 
\begin{equation}
  \begin{cases}
        \Psi_+=\Psi_-\begin{pmatrix}1&1\\0&1\end{pmatrix}& \quad z\in \R_+,\\
        \Psi_+=\Psi_-\begin{pmatrix}0&1\\-1&0\end{pmatrix}& \quad z\in \R_-,\\
        \Psi_+=\Psi_-\begin{pmatrix}1&0\\1&1\end{pmatrix}& 
\quad z\in \R_+e^{2\pi i/3}, \R_+e^{4\pi i/3}.\\
  \end{cases}
\end{equation}

Let $\mathcal{O}_{\xi}$ and ${\mathcal{O}}_{\bar{\xi}}$ be neighborhoods 
around $\xi$ and $\bar{\xi}$ of size 
$\epsilon''\sqrt{\gamma-1}$, respectively, where $\epsilon''>0$ is a small,  
fixed number chosen to satisfy~\eqref{e-z25},~\eqref{e-z26} below.
Since $\frac12 |\xi-\bar{\xi}|=\frac2{\gamma}\sqrt{\gamma-1} > 
\epsilon''\sqrt{\gamma-1}$, 
$\mathcal{O}_{\xi}$ and ${\mathcal{O}}_{\bar{\xi}}$ have no intersection, 
provided 
\begin{equation}\label{e-z25}
   0<\epsilon''<1.
\end{equation}
For definiteness we assume that $\partial\mathcal{O}_{\xi}$ and 
$\partial{\mathcal{O}}_{\bar{\xi}}$ are oriented counterclockwise.
In $\mathcal{O}_{\xi}$, a simple substitution shows that 
$\wt{\alpha}(z)=\frac23(z-\xi)^{3/2}G(z)$,
where $G$ is analytic and $G(\xi)=(\gamma-1)^{3/4}
e^{-i(\frac32\theta_c+\frac34\pi)}$.
Here $(z-\xi)^{3/2}=|z-\xi|^{3/2}e^{\frac32 i arg(z-\xi)}$ and
$\theta_c-\pi/2 < arg(z-\xi) <\theta_c+ 3\pi /2$.
Define
\begin{equation}
\lambda(z)\equiv (z-\xi)(G(z))^{2/3},
\end{equation}
where $(G(z))^{2/3}$ is analytic in $\mathcal{O}_{\xi}$ and 
$(G(z))^{2/3}\to (\gamma-1)^{1/2}e^{-i(\theta_c+\pi /2)}$ as $z\to\xi$.
Of course, $\lambda^{3/2}=\frac32 \wt{\alpha}$.
\begin{figure}[ht]
 \centerline{\epsfig{file=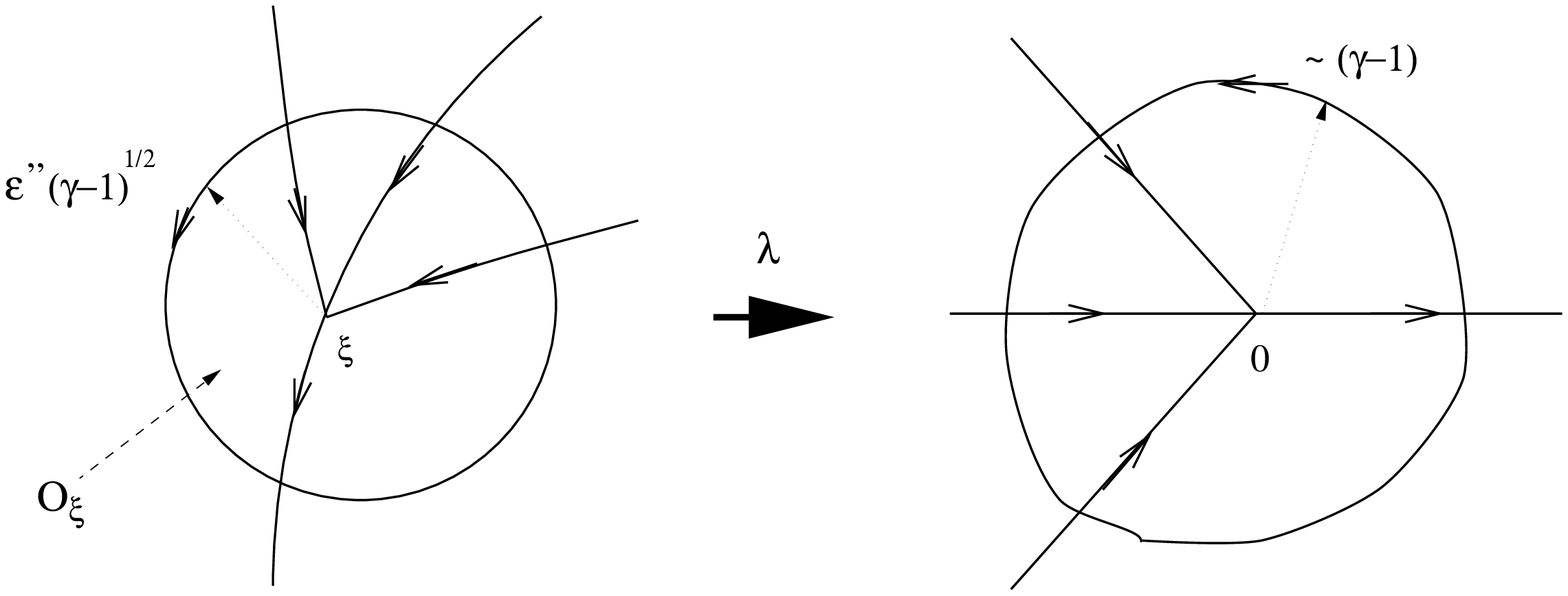, width=10cm}}
 \caption{}\label{fig-maptil}
\end{figure}
It is a simple calculus question to verify that we may choose $\epsilon''$ 
sufficiently small so that 
\begin{equation}\label{e-z26}
  \begin{split}
   &\text{$z\to\lambda(z)$ is a bijection from $\mathcal{O}_{\xi}$ 
onto an open neighborhood}\\
   &\text{of $0$ in the $\lambda$-plane, 
of radius $\sim (\gamma-1)$.}
  \end{split}
\end{equation}

Define $\Sigma^{(5)}\cap\mathcal{O}_{\xi}\equiv 
\{ z\in\mathcal{O}_{\xi} : \lambda(z)\in\Sigma''\}$.
As in the construction in \cite{DZ}, set (cf. (4.34) in \cite{DZ}) 
\begin{equation}\label{8.20}
   E(z) =\begin{pmatrix}1&-1\\-1&-1\end{pmatrix}
\sqrt\pi e^{i\pi/6}q^{\frac{\sigma_3}{6}}
\biggl((z-\bar\xi)(G(z))^{2/3}\biggr)^{\frac{\sigma_3}{4}}, 
\end{equation}
and for $z\in\mathcal{O}_{\xi}-\Sigma^{(5)}$, define the parametrix 
for $m^{(5)}$ by  
\begin{equation}
    m_p(z)=E(z)\Psi(q^{2/3}\lambda(z))e^{q\wt{\alpha}(z)\sigma_3}.
\end{equation}
Then $m_p$ satisfies the same jump conditions 
on $\Sigma^{(5)}\cap\mathcal{O}_{\xi}$ as $m^{(5)}$ ; $m_{p+}=m_{p-}v^{(5)}$.
And if $q$ becomes large, then 
for $z\in\partial\mathcal{O}_{\xi}$,
$|q^{2/3}\lambda(z)|\geq cq^{2/3}(\gamma-1) 
\geq cM_3/2^{1/3}$.
Therefore, 
\begin{equation}\label{e-z29}
 \begin{split}
  &\text{if $M_3$ is sufficiently large so that the leading terms dominate}\\ 
  &\text{in the asymptotics for the Airy functions in $\Psi(q^{2/3}\lambda(z))$ 
(see e.g.\cite{AS}),}
 \end{split}
\end{equation}
then by the explicit choice of $E(z)$ in~\eqref{8.20}, we find 
for $z\in\partial\mathcal{O}_{\xi}$, 
\begin{equation}\label{e-z27}
  m_{p+}=m^{(5,\infty)}(z)\biggl(I+O\bigl(\frac1{q\lambda^{3/2}(z)}\bigr)\biggr)
=m^{(5,\infty)}(z)\biggl(I+O\bigl(\frac1{q(\gamma-1)^{3/2}}\bigr)\biggr). 
\end{equation}

Noting the symmetry $m^{(5)}=\overline{m^{(5)}(\bar{z})}$, 
define $\Sigma^{(5)}\cap\mathcal{O}_{\bar{\xi}}\equiv 
\overline{\Sigma^{(5)}\cap\mathcal{O}_{\xi}}$, and 
for $z\in\mathcal{O}_{\bar{\xi}}-\Sigma^{(5)}$ 
set $m_p(z)=\overline{m_p(\bar{z})}$.
We now extend $\Sigma^{(5)}\cap(\mathcal{O}_{\xi}\cup\mathcal{O}_{\bar{\xi}})$ 
to $\Sigma^{(5)}$ to have the same general shape as in Figure~\ref{fig-Ctilde}. 
Finally, for $z\in\C-({\mathcal{O}}_\xi 
\cup {\mathcal{O}}_{\bar{\xi}} \cup\Sigma^{(5)})$, 
define $m_p(z)=m^{(5,\infty)}$.

Set $\mathcal{O}= {\mathcal{O}}_\xi\cup {\mathcal{O}}_{\bar{\xi}}$.
Then $\wt{R}\equiv m^{(5)}m_p^{-1}$ solves a RHP on 
$\Sigma^{(5)}\cup\partial\mathcal{O}$ with 
the jump matrix $v_{\wt{R}}=m_{p-}v^{(5)}v_p^{-1}m_{p-}^{-1}$, 
\begin{equation}
  \begin{cases}
    v_{\wt{R}}= I& \text{on} \quad (\Sigma^{(5)}\cap\mathcal{O})\cup C_2,\\
    v_{\wt{R}}= m^{(5,\infty)} \begin{pmatrix} 1&e^{-2q\wt{\alpha}}\\ 0&1 
\end{pmatrix} (m^{(5,\infty)})^{-1}& 
\text{on} \quad  C_1 \cap\bar{\mathcal{O}}^c, \\
    v_{\wt{R}}= m^{(5,\infty)} \begin{pmatrix} 1&0\\ 
e^{2q\wt{\alpha}}&1 \end{pmatrix}
(m^{(5,\infty)})^{-1}& \text{on} \quad 
(\wt{C}_{inside}\cup\wt{C}_{outside})\cap\bar{\mathcal{O}}^c ,\\
    v_{\wt{R}}=I+O(\frac1{q(\gamma-1)^{3/2}})& \text{on} 
\quad \partial\mathcal{O}.\\
  \end{cases}
\end{equation}

Let $\epsilon_1>0$ be a fixed, small number : for example, we may take 
$\epsilon_1=\epsilon'$ satisfying~\eqref{e-z21},~\eqref{e-z22} above.
Choose $\delta_4$ sufficiently small so that 
\begin{equation}\label{e-z30}
   \xi,\bar{\xi} \in \{z : |z+1|\le\epsilon_1\}
\end{equation}
for $1\le\gamma\le 1+\delta_4$.
By calculations similar to~\eqref{8.2} and~\eqref{e-z31}, 
$Re(2\tilde\alpha(z))\leq -c(\gamma-1)^{3/2}$ for
$z\in (\wt{C}_{inside}\cup\wt{C}_{outside})\cap \{|z+1|\le\epsilon_1\}
\cap \bar{\mathcal{O}}^c$ (in fact the estimate is true on the full set 
$(\wt{C}_{inside}\cup\wt{C}_{outside})\cap \{|z+1|\le\epsilon_1\}$) 
and also $|e^{2q\tilde\alpha(z)}|\leq e^{-cq}$ for 
$z\in (\wt{C}_{inside}\cup\wt{C}_{outside})\cap \{|z+1|>\epsilon_1\}$.
Thus, 
$|e^{2q\tilde\alpha(z)}|\leq e^{-cq(\gamma-1)^{3/2}}$ for 
$z\in (\wt{C}_{inside}\cup\wt{C}_{outside})\cap \mathcal{O}^c$. 
Also, by calculations similar to~\eqref{8.2} and~\eqref{e-z31} again,
$Re(-2\tilde\alpha(z))\leq -c(\gamma-1)^{3/2}$ for 
$z\in C_1\cap \mathcal{O}^c$. 
Therefore we have $L^{\infty}$ estimation 
\begin{equation}\label{e-z32}
   \|v_{\wt{R}} -I \|_{L^\infty(\Sigma^{(5)}\cap\mathcal{O}^c)}
\leq Ce^{-cq(\gamma-1)^{3/2}}.
\end{equation}

Furthermore, from calculations similar to~\eqref{6.9},  
on $\wt{C}_{inside}\cap\mathcal{O}^c\cap \{Im(z)\geq0\}$,  
using $|u+u_0|\geq |u-u_0|$ on the integration contour for the second inequality, 
\begin{equation}
  \begin{split}
    \int |e^{-2q\wt{\alpha}(z)}||dz|
    &\le \int_{\{u=u_0+xe^{-i\pi/3} : x\geq c\sqrt{\gamma-1}\}}
Ce^{-qc|u^2-u^2_0|^{3/2}}du + Ce^{-cq}\\
    &\le \int_{c\sqrt{\gamma-1}}^{\infty} 
Ce^{-qcx^{3}} dx + Ce^{-cq}\\
    &\le  \frac{C}{q(\gamma-1)}.
  \end{split}
\end{equation}
The same calculations apply to the other part of 
$\wt{C}_{inside}\cap\mathcal{O}^c$ and also 
to $\wt{C}_{outside}\cap\mathcal{O}^c$, 
so that 
$\|v_{\wt{R}}-I\|_{L^1((\wt{C}_{inside}\cup\wt{C}_{outside})\cap\mathcal{O}^c)}
\le \frac{C}{q(\gamma-1)}$.
On the other hand, $length(\partial\mathcal{O})\le C\sqrt{\gamma-1}$ and 
$length(C_1\cap\mathcal{O}^c)\le C\sqrt{\gamma-1}$, and hence, by the 
above $L^\infty$ estimates,  
$\|v_{\wt{R}}-I\|_{L^1(\partial\mathcal{O})}\le C/(q(\gamma-1))$ and 
$\|v_{\wt{R}}-I\|_{L^1(C_1\cap\mathcal{O}^c)}\le C/(q(\gamma-1))$.
Thus 
\begin{equation}\label{e-z33}
   \|v_{\wt{R}}-I\|_{L^1(\Sigma^{(5)}\cup\partial\mathcal{O})}\le 
\frac{C}{q(\gamma-1)}.
\end{equation}

Using the choice of $\mathcal{O}_{\xi}$ and $\mathcal{O}_{\bar{\xi}}$, direct 
calculation shows that $m^{(5,\infty)}$, hence $(m^{(5,\infty)})^{-1}$ (as 
$\det m^{(5,\infty)}=1$), are uniformly bounded for $\gamma$ in the region 
$1+\frac{M_3}{2^{1/3}q^{2/3}}\le \gamma \le 1+\delta_4$,  
for $z\in\mathcal{O}^c - \Sigma^{(5)}$.
On the other hand, even though the contour $\Sigma^{(5)}\cup\partial\mathcal{O}$ 
varies with $\gamma$ and $q$, the length of $\Sigma^{(5)}\cup\partial\mathcal{O}$ 
is uniformly bounded for $1+\frac{M_3}{2^{1/3}q^{2/3}}\le \gamma \le 1+\delta_4$.
Also a simple scaling argument shows that the Cauchy operators $C_\pm$ on 
$L^2(\Sigma^{(5)}\cup\partial\mathcal{O})$ are uniformly bounded for 
$1+\frac{M_3}{2^{1/3}q^{2/3}}\le \gamma \le 1+\delta_4$.
Therefore,  
\begin{equation}\label{e-z28}
  \begin{split}
    \|C_{w_{\wt{R}}}\|_{L^2(\Sigma^{(5)}\cup\partial\mathcal{O})
\to L^2(\Sigma^{(5)}\cup\partial\mathcal{O})} 
&\le C\|w_{\wt{R}}\|_{L^\infty(\Sigma^{(5)}\cup\partial\mathcal{O})}\\
    &\le \frac{C}{q(\gamma-1)^{3/2}}+Ce^{-Cq(\gamma-1)^{3/2}}\\
    &\le \frac{C}{M_3^{2/3}}+Ce^{-CM_3^{3/2}}\\
    &\le \frac12 <1
  \end{split}
\end{equation}
provided that $M_3$ is sufficiently large.

From~\eqref{e-sam5} and~\eqref{e-z33}, we have 
\begin{equation*}
  \begin{split}
    |\wt{R}_{22}(0) -1| &= \biggl| \frac1{2\pi i} 
\int_{\Sigma^{(5)}\cup\partial\mathcal{O}} \biggl( 
\frac{w_{\wt{R}}-[(I-C_{w_{\wt{R}}})^{-1}C_{w_{\wt{R}}}I](z)
w_{\wt{R}}(z)}{z}\,ds 
\biggr)_{22} \biggr|\\
    &\le C(\|w_{\wt{R}}\|_{L^1(\Sigma^{(5)}\cup\partial\mathcal{O})}+
\|w_{\wt{R}}\|^2_{L^2(\Sigma^{(5)}\cup\partial\mathcal{O})}) \\
    &\le C\|w_{\wt{R}}\|_{L^1(\Sigma^{(5)}\cup\partial\mathcal{O})} 
\qquad\text{, as $\|w_{\wt{R}}\|
_{L^\infty(\Sigma^{(5)}\cup\partial\mathcal{O})}$ is bounded,}\\
    &\le\frac{C}{q(\gamma-1)},
  \end{split}
\end{equation*}
and 
\begin{equation*}
 |\wt{R}_{21}(0)| \le \frac{C}{q(\gamma-1)}.
\end{equation*}
Therefore, from $m^{(5)}_{21}(0)=\wt{R}_{22}(0)m^{(5,\infty)}_{21}(0)+
\wt{R}_{21}(0)m^{(5,\infty)}_{11}(0)$, we obtain 
\begin{equation}
  \kappa^2_{q-1}
=-e^{q(-\gamma+\log \gamma+1)}m^{(5)}_{21}(0)
=e^{q(-\gamma+\log \gamma+1)} \frac1{\sqrt{\gamma}} 
\biggl(1+O\bigl( \frac1{q(\gamma-1)} \bigr)\biggr).
\end{equation}
Note that this is consistent with the result~\eqref{8.1} for case (i) 
where $\gamma-1\geq \delta_3$.

Summarizing, we have proven the following results.
\begin{lem}\label{lem-K2}
  Let $\delta_3, M_4>0$ be fixed numbers.
  Let $\delta_4>0$ be a fixed sufficiently 
small number satisfying~\eqref{e-z30}, 
and $M_3>0$ be a fixed, sufficiently large number 
satisfying~\eqref{e-z29} and~\eqref{e-z28}.
As $q \to \infty$, we have the following asymptotics.
  \begin{enumerate}
    \item If $1+\delta_3 \le \gamma$,
\begin{equation*}
  \kappa^2_{q-1}
= e^{q(-\gamma+\log \gamma+1)}\frac1{\sqrt{\gamma}}
\biggl(1+O\bigl(\frac1{q}\bigr)\biggr),
\end{equation*}
where the error is uniform for $1+\delta_4 \le \gamma\le L$ for 
any fixed $L< \infty$.
    \item If $1+\frac{M_3}{2^{1/3}q^{2/3}} \le \gamma \le 1+\delta_4$,
\begin{equation*}
  \kappa^2_{q-1}=e^{q(-\gamma+\log \gamma+1)} \frac1{\sqrt{\gamma}}
\biggl(1+O\bigl(\frac1{q(\gamma-1)}\bigr) \biggr),
\end{equation*}
where the error is uniform in the region.
    \item If $1< \gamma \le 1+\frac{M_4}{2^{1/3}q^{2/3}}$, 
\begin{equation*}
  \bigl|\kappa^2_{q-1}
- 1 -\frac{i2^{4/3}}{q^{1/3}}m^{PII}_{1,22}(-t)\bigr|\le  \frac{C(M_4)}{q^{2/3}},
\end{equation*}
where $t$ is defined by $\gamma=1+\frac{t}{2^{1/3}q^{2/3}}$, $0\le t\le M_4$.
  \end{enumerate}
\end{lem}

Note that, comparing Lemma ~\ref{lem-K2} (iii) with Lemma ~\ref{lem-K1} (iii), 
we have same result everywhere in the region $1-\frac{M}{2^{1/3}q^{2/3}} \le 
\gamma \le 1+\frac{M}{2^{1/3}q^{2/3}}$, 
\begin{equation}\label{e-z36}
     \bigl|\kappa^2_{q-1}
- 1 -\frac{i2^{4/3}}{q^{1/3}}m^{PII}_{1,22}(t)\bigr|\le  \frac{C(M)}{q^{2/3}},
\end{equation}
where $t$ is defined by $\gamma=1-\frac{t}{2^{1/3}q^{2/3}}$, 
and $M$ is any fixed positive number.

Also note from Lemma~\ref{lem-K2} (ii), that as $q\to\infty$, 
\begin{equation}\label{e-z37}
   \log\kappa^2_{q-1}\le q(-\gamma+\log\gamma +1) + \frac{C^{\#}}{q(\gamma-1)},
\end{equation}
where $C^{\#}$ is independent of $M_3$ and is fixed once $\delta_4$ 
satisfying~\eqref{e-z30} is determined.

%%%%====================================sec6.tex

\bigskip
\section{\bf{Asymptotics of $\mathbf{\phi_n(\lambda)}$ as $\mathbf{ n\to\infty}$ }}
\label{s-phi}

In this Section, using Lemmas~\ref{lem-K1} and~\ref{lem-K2}, 
we obtain the large $n$ behavior of $\phi_n(\lambda)$.
In the following $\delta_5,\delta_6,\delta_7$ are fixed numbers 
between $0$ and $1$, and 
$M_5,M_6,M_7$ are fixed and positive.
These numbers are free apart from the following requirements : 
\begin{align*}
   &\text{(a) $\delta_6$ satisfies~\eqref{e-z30},}\\
   &\text{(b) $M_5\geq 1$ satisfies~\eqref{e-sec3.6},}\\
   &\text{(c) $\frac12 M_6\geq 1$ satisfies~\eqref{e-sec3.6}, and}\\
   &\text{(d) $M_7\geq 1$ satisfies~\eqref{e-z29}, ~\eqref{e-z28} 
and condition~\eqref{e-z39} below.}
\end{align*}
We consider the following five cases for $\lambda>0$ and $n$ :
\begin{enumerate}
   \item $0\le \frac{2\sqrt{\lambda}}{n+1}\le 1-\delta_5$
   \item $\frac12 \le \frac{2\sqrt{\lambda}}{n+1} 
\le 1-\frac{M_5}{2^{1/3}{(n+1)}^{2/3}}$
   \item $1-\frac{M_6}{2^{1/3}{(n+1)}^{2/3}} \le \frac{2\sqrt{\lambda}}{n+1} 
\le 1+\frac{M_6}{2^{1/3}{(n+1)}^{2/3}}$
   \item $1+\frac{M_7}{2^{1/3}{(n+1)}^{2/3}} \le \frac{2\sqrt{\lambda}}{n+1}
\le 1+\delta_6$
   \item $1+\delta_7 \le \frac{2\sqrt{\lambda}}{n+1}$
\end{enumerate}

\bigskip
Consider case (i). 
For any $k\geq n$, $\frac{2\sqrt{\lambda}}{k+1}\le 1-\delta_5$. 
From Lemma~\ref{lem-K1} (i), we have as $n\to\infty$, 
\begin{equation}
   |\log{\phi_n(\lambda)}|
   =\bigl|\sum_{k=n}^{\infty} \log{\kappa^2_k}(\lambda)\bigr|
   \le \sum_{k=n}^{\infty} Ce^{-ck} \le Ce^{-cn}.
\end{equation}

\bigskip
Consider case (ii).
We split the sum into two pieces.
\begin{equation*}
   \begin{split}
      \log{\phi_n(\lambda)}
      &=\sum_{k=n}^{\infty} \log{\kappa^2_k}(\lambda)\\
      &=\sum_{(1)} \log{\kappa^2_{k}}(\lambda)
      + \sum_{(2)} \log{\kappa^2_{k}}(\lambda)
  \end{split}
\end{equation*}
where $(1)$ and $(2)$ represent the regions
\begin{align*}
  (1)\ \  &n+1\le k+1 \le 4\sqrt{\lambda},\\
  (2)\ \  &4\sqrt{\lambda} <k+1. 
\end{align*}
For $(1)$, $\frac12 \le \frac{2\sqrt{\lambda}}{k+1}\le
1-\frac{M_5}{2^{1/3}(k+1)^{2/3}}$.
From Lemma~\ref{lem-K1} (ii), for some constant $C$, independent of $M_5$ 
satisfying~\eqref{e-sec3.6},
\begin{equation*}
  \bigl| \log \kappa^2_k(\lambda) \bigr|
\le C\frac{ {e^{-\frac{2\sqrt{2}}{3}(k+1)
{(1-\frac{2\sqrt\lambda}{k+1})}^{3/2}}}}{{(k+1)}^{1/3}}.
\end{equation*}
Using the fact that 
$f(x)=\frac1{x^{1/3}}e^{-\frac{2\sqrt{2}}{3}x(1-\frac{2\sqrt{\lambda}}{x})^{3/2}}$ 
is monotone decreasing in the second inequality below, 
we have, as $n\to\infty$,  
\begin{equation}\label{e-z3}
  \begin{split}
    \bigl|\sum_{(1)} \log{\kappa^2_{k}}(\lambda) \bigr|
    &\le C\sum_{(1)} \frac{e^{-\frac{2\sqrt{2}}{3}(k+1)
(1-\frac{2\sqrt{\lambda}}{k+1})^{3/2}}}
{{(k+1)}^{1/3}}\\
    &\le C\int_{n+1}^{4\sqrt{\lambda}}
e^{-\frac{2\sqrt{2}}{3}x(1-\frac{2\sqrt{\lambda}}{x})^{3/2}}
\frac{dx}{x^{1/3}}
+C\frac{e^{-\frac{2\sqrt{2}}{3}(n+1)
(1-\frac{2\sqrt{\lambda}}{n+1})^{3/2}}}{{(n+1)}^{1/3}}\\
    &\le C(2\sqrt{\lambda})^{2/3}\int_{\frac{(n+1)}{2\sqrt{\lambda}}-1}^{1}
e^{-\frac{4\sqrt{2\lambda}y^{3/2}}{3(1+y)^{1/2}}}
\frac{dy}{(1+y)^{1/3}}
+Ce^{-\frac12(n+1)(1-\frac{2\sqrt{\lambda}}{n+1})^{3/2}}\\
    &\le C(2\sqrt{\lambda})^{2/3}\int_{\frac{(n+1)}{2\sqrt{\lambda}}-1}^{1}
e^{-\sqrt{\lambda}y^{3/2}} dy
+Ce^{-\frac12(n+1)(1-\frac{2\sqrt{\lambda}}{n+1})^{3/2}}\\
    &\le C\int_{\sqrt{\lambda}(\frac{(n+1)}{2\sqrt{\lambda}}-1)^{3/2}}^{\infty}
e^{-s}\frac{ds}{s^{1/3}}
+Ce^{-\frac12(n+1)(1-\frac{2\sqrt{\lambda}}{n+1})^{3/2}}\\
    &\le Ce^{-\sqrt{\lambda}(\frac{(n+1)}{2\sqrt{\lambda}}-1)^{3/2}}
+Ce^{-\frac12(n+1)(1-\frac{2\sqrt{\lambda}}{n+1})^{3/2}}\\
    &\le C\exp\biggl(-\frac12(n+1)
\bigl(1-\frac{2\sqrt{\lambda}}{n+1}\bigr)^{3/2}\biggr).
  \end{split}
\end{equation}
We use the change of variable $y=\frac{x}{2\sqrt{\lambda}}-1$ for 
the integral in the third line.
The fifth inequality is obtained from the substitution 
$s=\sqrt{\lambda}y^{3/2}$, 
and at the end, we have used $\frac{2\sqrt{\lambda}}{n+1} \le 1$.

For $(2)$, $\frac{2\sqrt{\lambda}}{k+1}\le \frac12$.
Therefore, from Lemma~\ref{lem-K1} (i), we have
\begin{equation*}
   \bigl|\sum_{(2)} \log{\kappa^2_k}(\lambda)\bigr|
   \le \sum_{k+1=[4\sqrt{\lambda}]}^{\infty} Ce^{-ck} \le Ce^{-cn}.
\end{equation*}

Summing up the above two calculations, we have, for case (ii), 
\begin{equation}\label{e-z35}
   |\log{\phi_n(\lambda)}| \le C\exp\biggl(-c(n+1)
\bigl(1-\frac{2\sqrt{\lambda}}{n+1}\bigr)^{3/2}\biggr),
\end{equation}
as $n\to\infty$.
Note again that the constants C,c can be taken independent of $M_5$.

\bigskip
Consider case (iii).
Set 
\begin{equation}\label{e-z38}
   \frac{2\sqrt{\lambda}}{n+1}=1-\frac{t}{2^{1/3}{(n+1)}^{2/3}}
\end{equation}
so that $-M_6\le t\le M_6$. 
We divide the sum into three pieces,  
\begin{equation*}
  \begin{split}
     \log{\phi_n(\lambda)}
     &=\sum_{k=n}^{\infty} \log{\kappa^2_k}(\lambda)\\
     &=\sum_{(1)} \log{\kappa^2_{k}}(\lambda) 
     + \sum_{(2)} \log{\kappa^2_{k}}(\lambda) 
     + \sum_{(3)} \log{\kappa^2_{k}}(\lambda)
  \end{split}
\end{equation*}
where $(1),(2)$ and $(3)$ indicate the following regions :
\begin{align*}
  (1)\ \  &n+1\le k+1 \le (n+1)+\frac{(M_6-t)}{2^{1/3}}{(n+1)}^{1/3}\\
  (2)\ \  &(n+1)+\frac{(M_6-t)}{2^{1/3}}{(n+1)}^{1/3} <k+1 
< \frac32(n+1)-\frac{t}{2^{1/3}}(n+1)^{1/3}\\
  (3)\ \  &\frac32(n+1)-\frac{t}{2^{1/3}}(n+1)^{1/3} \le k+1.
\end{align*}

\medskip
For $(1)$, as $n\to\infty$,  
\begin{equation*}
  1-\frac{6M_6}{2^{1/3}{(k+1)}^{2/3}} \le 
\frac{2\sqrt{\lambda}}{k+1} \le 1+\frac{2M_6}{2^{1/3}{(k+1)}^{2/3}}.
\end{equation*} 
Hence from~\eqref{e-z36}, 
we have as $k\geq n\to\infty$,  
\begin{equation*}
\log \kappa^2_k(\lambda)
= \frac{i2^{4/3}}{{(k+1)}^{1/3}} m^{PII}_{1,22}
\biggl(2^{\frac13}{(k+1)}^{\frac23}(1-\frac{2\sqrt{\lambda}}{k+1})\biggr) 
+O_{M_6}\bigl(\frac1{k^{2/3}}\bigr).
\end{equation*}
This leads to 
\begin{equation*}
  \begin{split}
    &\sum_{(1)} \log{\kappa^2_{k}}(\lambda)\\
    &= \sum_{k+1=n+1}^{[(n+1)+\frac{(M_6-t)}{2^{1/3}}{(n+1)}^{1/3}]} 
\biggl[ \frac{i2^{4/3}}{{(k+1)}^{1/3}} 
m^{PII}_{1,22}(2^{\frac13}{(k+1)}^{\frac23}(1-\frac{2\sqrt{\lambda}}{k+1})) 
+ O_{M_6}\bigl(\frac1{k^{2/3}}\bigr) \biggr]\\
    &= \int_{(n+1)}^{(n+1)+\frac{(M_6-t)}{2^{1/3}}{(n+1)}^{1/3}} 
\frac{i2^{4/3}}{x^{1/3}} 
m^{PII}_{1,22}(2^{1/3}x^{2/3}(1-\frac{2\sqrt{\lambda}}{x})) dx  +
O_{M_6}\bigl(\frac1{n^{1/3}}\bigr)\\
    &= \int_{0}^{\frac{(M_6-t)}{2^{1/3}}} i2^{4/3} 
m^{PII}_{1,22}\biggl(2^{1/3}{(n+1)}^{2/3}
(1+\frac{2s}{3{(n+1)}^{2/3}} +\cdots)\\ 
    &\qquad \bigl(1-\frac{(n+1)-\frac{t}{2^{1/3}}(n+1)^{1/3}}
{(n+1)+s{(n+1)}^{1/3}}\bigr)\biggr)
     \frac{ds}
{(1+\frac{s}{3{(n+1)}^{2/3}}+\cdots)} +
O_{M_6}\bigl(\frac1{n^{1/3}}\bigr)\\
    &=\int_0^{\frac{(M_6-t)}{2^{1/3}}} i2^{4/3} 
m^{PII}_{1,22}\biggl((t+2^{1/3}s)(1-\frac{s}{3{(n+1)}^{2/3}} +\cdots)\biggr)(1-\frac{s}{3{(n+1)}^{2/3}} +\cdots)ds \\
&\qquad + O_{M_6}\bigl(\frac1{n^{1/3}}\bigr)\\
    &= \int_0^{\frac{(M_6-t)}{2^{1/3}}} i2^{4/3} 
m^{PII}_{1,22}(t+2^{1/3}s) ds 
+O_{M_6}\bigl(\frac1{n^{1/3}}\bigr)\\
    &= \int_{t}^{M_6} 2im^{PII}_{1,22}(y) dy 
+O_{M_6}\bigl(\frac1{n^{1/3}}\bigr).
  \end{split}
\end{equation*}
The fourth equation is obtained using the change of variable 
$x=(n+1)+s{(n+1)}^{1/3}$, and for the sixth equation, 
we use the fact that $\frac{d}{dt}m^{PII}_{1,22}(t)$ 
is uniformly bounded for $-M_6\le t\le M_6$ 
(see the remark below~\eqref{e-z34}).
To pass from the second to the third line, note that 
for integers $b>a$, 
\begin{equation}\label{e-z40}
  \begin{split}
     |\sum_{n=a}^{b-1} f(x) - \int_a^b f(x)| &\le 
\sum_{n=a}^{b-1} sup\{ |f(\alpha)-f(\beta)| : n\le \alpha,\beta\le n+1\}\\ 
     & \le \|f'\|_{L^{\infty}(a,b)}(b-a).
  \end{split}
\end{equation}
For the case at hand, a simple calculation shows that 
$$\|f'\|_{L^\infty\bigl(n+1,[(n+1)+\frac{(M_6-t)}{2^{1/3}}{(n+1)}^{1/3}]\bigr)} 
\le C(M_6)/n^{2/3}.$$
Also, the contribution to the integral from the interval 
$\bigl((n+1)+\frac{(M_6-t)}{2^{1/3}}{(n+1)}^{1/3},
[(n+1)+\frac{(M_6-t)}{2^{1/3}}{(n+1)}^{1/3}]+1\bigr)$ is $O_{M_6}(1/n^{1/3})$.

\medskip
For $(2)$, 
\begin{equation*}
   \frac12 \le \frac{2\sqrt\lambda}{k+1} \le 
1-\frac{\frac12M_6}{2^{1/3}(k+1)^{2/3}}\quad 
\text{as}\quad  k\geq n\rightarrow \infty.
\end{equation*}
As $n\to\infty$, by a calculation similar to the case (ii), again using 
the monotonicity of $f(x)=\frac1{x^{1/3}}
e^{-\frac{2\sqrt{2}}{3}x(1-\frac{2\sqrt{\lambda}}{x})^{3/2}}$ 
for the second inequality, 
we have  
\begin{equation*}
  \begin{split}
    \bigl|\sum_{(2)} \log{\kappa^2_{k}}(\lambda) \bigr| 
    &\le C\sum_{(2)} \frac{e^{-\frac{2\sqrt{2}}{3}(k+1)
(1-\frac{2\sqrt{\lambda}}{k+1})^{3/2}}}
{{(k+1)}^{1/3}}\\
    &\le C\int_{(n+1)+\frac{(M_6-t)}{2^{1/3}}{(n+1)}^{1/3}}
^{\frac32(n+1)-\frac{t}{2^{1/3}}(n+1)^{1/3}} 
e^{-\frac{2\sqrt{2}}{3}x(1-\frac{2\sqrt{\lambda}}{x})^{3/2}} 
\frac{dx}{x^{1/3}}
+Ce^{-(\frac{M_6}2)^{3/2}}\\
    &\le C(2\sqrt{\lambda})^{2/3}
\int_{\frac{M_6(n+1)^{1/3}}{2^{4/3}\sqrt{\lambda}}}
^{\frac{(n+1)}{4\sqrt{\lambda}}} 
e^{-\frac{4\sqrt{2\lambda}}{3}(\frac{y^3}{1+y})^{1/2}}
\frac{dy}{(1+y)^{1/3}}
+Ce^{-(\frac{M_6}2)^{3/2}}\\
    &\le C(2\sqrt{\lambda})^{2/3}
\int_{\frac{M_6(n+1)^{1/3}}{2^{4/3}\sqrt{\lambda}}}
^{\frac{(n+1)}{4\sqrt{\lambda}}} 
e^{-\sqrt{\lambda}y^{3/2}} dy
+Ce^{-(\frac{M_6}2)^{3/2}}\\
    &\le C \int_{(\frac{n+1}{2\sqrt{\lambda}})^{\frac12}
(\frac{M_6}{2})^{\frac32}}^{\infty} 
e^{-s} \frac{ds}{s^{1/3}}
+Ce^{-(\frac{M_6}2)^{3/2}}\\
    &\le C\int_{\frac14M_6^{3/2}}^{\infty} \frac{e^{-s}}{s^{1/3}} ds
+Ce^{-(\frac{M_6}2)^{3/2}}\\
    &\le Ce^{-\frac14M_6^{3/2}}
+Ce^{-(\frac{M_6}2)^{3/2}}
\le Ce^{-\frac14M_6^{3/2}}.
  \end{split}
\end{equation*}
The first inequality follows from Lemma~\ref{lem-K1} (ii) (note that, 
by assumption, $\frac12 M_6$ satisfies~\eqref{e-sec3.6}).
For the second line, in order to control the contribution 
to the integral from the interval 
$[[(n+1)+\frac{(M_6-t)}{2^{1/3}}{(n+1)}^{1/3}], 
(n+1)+\frac{(M_6-t)}{2^{1/3}}{(n+1)}^{1/3}]$, we use the inequality 
$1+\frac{M_6-t}{2^{1/3}(n+1)^{2/3}} \le 1+\frac{2^{2/3}M_6}{(n+1)^{2/3}}
\le \frac{32}{9}$ for large enough $n$.
For the third line, we use the change of variable 
$y=\frac{x}{2\sqrt{\lambda}}-1$, and 
for the fourth line, we use the inequality 
$\frac{4\sqrt{\lambda}}{(n+1)\delta} \geq 1-\frac{M_6}{2^{1/3}(n+1)^{2/3}}
\geq \frac9{23}$ for sufficiently large $n$. 
The fifth equation is obtained from the substitution $s=\sqrt{\lambda}y^{3/2}$,
and for the sixth line, we have used the inequality 
$\frac{2\sqrt{\lambda}}{n+1} \le 1+\frac{M_6}{2^{1/3}(n+1)^{2/3}} 
\le 2$ for sufficiently large $n$.

\medskip
For $(3)$, as $n\to\infty$,  $0\le \frac{2\sqrt\lambda}{k+1} \le \frac34$, 
which yields, from Lemma~\ref{lem-K1} (i), 
\begin{equation*}
   \bigl| \sum_{(3)} \log{\kappa^2_{k}(\lambda)} \bigr| \le Ce^{-cn}.
\end{equation*}

Summing up all these calculations, for 
$\frac{2\sqrt{\lambda}}{n+1} = 1- \frac{t}{2^{1/3}{(n+1)}^{2/3}}$ 
with $-M_6\le t\le M_6$, we have, as $n\to\infty$, 
\begin{equation*}
  \bigl| \log{\phi_n(\lambda)} - \int_{t}^{M_6} 2im^{PII}_{1,22}(y) dy \bigr|
\le \frac{C(M_6)}{n^{1/3}}+Ce^{-\frac14M_6^{3/2}}, 
\end{equation*}
for a constant $C(M_6)$ which depends on $M_6$, and for a 
constant $C$ which is independent of $M_6$. 
Using the asymptotics of $m^{PII}_{1,22}(x)$ as 
$x\rightarrow +\infty$ (see~\eqref{e-aympm_1}), we have 
(recall $\frac12 M_6\geq 1$) 
\begin{equation}
  \bigl| \log{\phi_n(\lambda)} - \int_{t}^{\infty} 2im^{PII}_{1,22}(y) dy \bigr|
\le \frac{C(M_6)}{n^{1/3}}+Ce^{-\frac14M_6^{3/2}}. 
\end{equation}

\bigskip
Now we consider case (iv), 
$1+\frac{M_7}{2^{1/3}{(n+1)}^{2/3}} \le 
\frac{2\sqrt{\lambda}}{n+1} \le 1+\delta_6$.
We write 
\begin{equation*}
     \log{\phi_n(\lambda)} 
     =\sum_{(1)} \log{\kappa^2_{k}} + \sum_{(2)} \log{\kappa^2_{k}},
\end{equation*}
where $(1),(2)$ indicate the following regions : 
\begin{align*}
  (1)\quad  &n+1\le k+1 \le 2\sqrt\lambda - \frac{M_7}{2^{1/3}}(n+1)^{1/3}\\ 
  (2)\quad  &2\sqrt\lambda - \frac{M_7}{2^{1/3}}(n+1)^{1/3}\leq k+1.
\end{align*}

\medskip
For $(1)$, we have for $n$ sufficiently large,  
$1+\frac{M_7}{2^{1/3}(k+1)^{2/3}} \leq \frac{2\sqrt\lambda}{k+1} \leq 1+\delta_6$.
Therefore, using~\eqref{e-z37}, we obtain 
\begin{equation}\label{e-z41}
\begin{split}
  &\sum_{(1)} \log{\kappa^2_{k}(\lambda)}\\
  &\le \sum_{k+1=n+1}^{[2\sqrt\lambda -\frac{M_7}{2^{1/3}}(n+1)^{1/3}]}
        (k+1) \biggl( -\frac14(\frac{2\sqrt\lambda}{k+1}-1)^2\biggr) 
      +\frac{C^{\#}}{2\sqrt{\lambda}-(k+1)}\\
  &\leq -\frac14\int_{n+1}^{2\sqrt\lambda -\frac{M_7}{2^{1/3}}(n+1)^{1/3}}
       x\bigl(\frac{2\sqrt\lambda}{x}-1\bigr)^2 dx
       + \int_{n+1}^{2\sqrt\lambda -\frac{M_7}{2^{1/3}}(n+1)^{1/3}}
       \frac{C^{\#}}{2\sqrt{\lambda}-x}dx +C(M_7)\\
  &\leq -\frac14\int_{\frac{n+1}{2\sqrt\lambda}}
   ^{1-\frac{M_7}{2^{1/3}(n+1)^{2/3}}\frac{n+1}{2\sqrt\lambda}}
    (2\sqrt\lambda)^2 \frac{(1-y)^2}{y} dy
    + C^{\#}\log\biggl(\frac{2\sqrt{\lambda}-(n+1)}{(n+1)^{1/3}} \biggr)
    +C(M_7)\\
  &\leq -\frac14\int_{\frac{n+1}{2\sqrt\lambda}}
   ^{1-\frac{M_7}{2^{1/3}(n+1)^{2/3}}\frac{n+1}{2\sqrt\lambda}}
    (2\sqrt\lambda)^2 (1-y)^2 dy
    + C^{\#}\frac{2\sqrt{\lambda}-(n+1)}{(n+1)^{1/3}} 
    +C(M_7)\\
  &\leq \frac14\int_{1-\frac{n+1}{2\sqrt\lambda}}
   ^{\frac{M_7}{2^{1/3}(n+1)^{2/3}}\frac{n+1}{2\sqrt\lambda}}
    (2\sqrt\lambda)^2 z^2 dz
    + C^{\#}(n+1)^{2/3}(\frac{2\sqrt\lambda}{n+1}-1)
    +C(M_7)\\
  &\leq \frac1{12}(2\sqrt\lambda)^2
    \biggl[( \frac{M_7}{2^{1/3}(n+1)^{2/3}}\frac{n+1}{2\sqrt\lambda})^3
    - (1-\frac{n+1}{2\sqrt\lambda})^3 \biggr]\\
  & \qquad + C^{\#}(n+1)^{2/3}(\frac{2\sqrt\lambda}{n+1}-1)
    +C(M_7)\\
  &\leq \frac{M_7^3}{24}(\frac{n+1}{2\sqrt\lambda}) 
    -\frac{n+1}{24\sqrt\lambda}(n+1)^2(\frac{2\sqrt\lambda}{n+1}-1)^3 
    + C^{\#}(n+1)^{2/3}(\frac{2\sqrt\lambda}{n+1}-1)
    +C(M_7)\\
  &\leq - \frac1{48}\biggl[2^{1/3}(n+1)^{2/3}
    (\frac{2\sqrt\lambda}{n+1}-1)\biggr]^3 
    + C^{\#}(n+1)^{2/3}(\frac{2\sqrt\lambda}{n+1}-1)
    +C(M_7)\\
  &\le - \frac1{96}\biggl[2^{1/3}(n+1)^{2/3}
    (\frac{2\sqrt\lambda}{n+1}-1)\biggr]^3 +C(M_7).
\end{split}
\end{equation}
The first line follows from the inequality 
$ -\gamma + \log \gamma +1 \leq -\frac{(\gamma-1)^2}{4}$ 
for $1\le \gamma \le 2$ ( note from~\eqref{e-z37} that $C^{\#}$ is independent of 
$M_7$).
In the second line, we use the monotonicity of 
$x(\frac{2\sqrt{\lambda}}{x}-1)^2$ and of $(2\sqrt{\lambda}-x)^{-1}$ 
in the region $n+1\le x\le 2\sqrt{\lambda}-\frac{M_7}{2^{1/3}}(n+1)^{1/3}$.
In the succeeding lines, we have used the changes of variables 
$x=2\sqrt{\lambda} y$, $1-y=z$ and $2\sqrt{\lambda} z^{3/2}=s$.
For the last line, note that $2^{1/3}(n+1)^{2/3}(\frac{2\sqrt\lambda}{n+1}-1)
\geq M_7$, and we require 
\begin{equation}\label{e-z39}
    M_7\geq \sqrt{96C^{\#}}.
\end{equation}

\textbf{Remark} : In estimating the sum in the second line of~\eqref{e-z41} 
by an integral, the monotonicity of the integrand plays a crucial role : 
we cannot, for example, use an estimate of the form~\eqref{e-z40}, as 
the derivative is not sufficiently small.

For $(2)$, we have 
$\frac{2\sqrt{\lambda}}{k+1} \le 1+\frac{2M_7}{2^{1/3}{(k+1)}^{2/3}}$.
Calculations similar to 
the previous cases (i), (ii) and (iii) show that 
\begin{equation*}
   \sum_{(2)} \log{\kappa^2_{k}(\lambda)} 
   \le  \bigl| \int_{-2M_7}^{\infty} 2im^{PII}_{1,22}(y) dy \bigr|
+ \frac{C(M_7)}{n^{1/3}}+Ce^{-M_7} \le C(M_7).
\end{equation*}
The result follows by splitting the sum $\sum_{(2)}$ into the following 
regions $2\sqrt{\lambda}-\frac{M_7}{2^{1/3}}(n+1)^{1/3}\le k+1 \le
2\sqrt{\lambda}+\frac{M_7}{2^{1/3}}(n+1)^{1/3}$, 
$2\sqrt{\lambda}+\frac{M_7}{2^{1/3}}(n+1)^{1/3} < k+1 <3\sqrt{\lambda}$ and 
$3\sqrt{\lambda} \le k+1$  : we leave the detail to the reader.

Therefore, for $1+\frac{M_7}{2^{1/3}{(n+1)}^{2/3}} \le 
\frac{2\sqrt{\lambda}}{n+1}
\le 1+\delta_6$, we have 
\begin{equation}
   \log \phi_n(\lambda) \le - \frac1{96}\biggl[2^{1/3}(n+1)^{2/3}
    (\frac{2\sqrt\lambda}{n+1}-1)\biggr]^3 + C(M_7) 
\end{equation}

\bigskip
For case (v), we use the estimation of \cite{Johan} given in the 
Lemma~\ref{lem-phi} (v) below.

\bigskip
Summarizing, we have
\begin{lem}\label{lem-phi}
  Let $0<\delta_5,\delta_6,\delta_7<1$ and $M_5,M_6,M_7>0$ be fixed number.
  Suppose that $\delta_6, M_5, M_6$ and $M_7$ satisfy conditions 
(a),(b),(c) and (d) given at the beginning of this Section, respectively.
  Set 
\begin{equation}\label{e-z42}
    t=2^{1/3}(n+1)^{2/3}\bigl(1-\frac{2\sqrt{\lambda}}{n+1}\bigr)
\quad\text{so that}\quad 
\frac{2\sqrt{\lambda}}{n+1}=1-\frac{t}{2^{1/3}(n+1)^{2/3}}.
\end{equation}
  We have the following estimates for the large $n$ 
behavior of $\phi_n(\lambda)$ : 
  \begin{enumerate}
    \item If $0\le \frac{2\sqrt{\lambda}}{n+1} \le 1-\delta_5$,  
\begin{equation*}
  \bigl| \log{\phi_n(\lambda)} \bigr| \le C\exp(-cn), 
\end{equation*}
for some constants $C,c$ which may depend on $\delta_5$.
    \item If $\frac12 \le \frac{2\sqrt{\lambda}}{n+1} 
\le 1-\frac{M_5}{2^{1/3}{(n+1)}^{2/3}}$ ,
\begin{equation*}
  \bigl| \log{\phi_n(\lambda)} \bigr| \le C\exp(-ct^{3/2}),
\end{equation*}
for constants C,c independent of $M_5$.
    \item If $1-\frac{M_6}{2^{1/3}{(n+1)}^{2/3}} \le 
\frac{2\sqrt{\lambda}}{n+1} \le 1+\frac{M_6}{2^{1/3}{(n+1)}^{2/3}}$, 
so that $-M_6\le t\le M_6$, 
there is a constant $C(M_6)$ which depends on $M_6$, 
and a constant $C$ which is independent of $M_6$,  
such that 
\begin{equation*}
  \bigl| \log{\phi_n(\lambda)} - \int_{t}^{\infty} 2im^{PII}_{1,22}(y) dy \bigr|
\le \frac{C(M_6)}{n^{1/3}}+Ce^{-\frac14M_6^{3/2}}.
\end{equation*}
    \item If $1+\frac{M_7}{2^{1/3}{(n+1)}^{2/3}} \le \frac{2\sqrt{\lambda}}{n+1}
\le 1+\delta_6$, 
\begin{equation*}
      \log \phi_n(\lambda) \le \frac1{96}t^3 + C(M_7),
\end{equation*}
for a constant $C(M_7)$.
   \item \cite{Johan} If $1+\delta_7 \le \frac{2\sqrt{\lambda}}{n+1}$,
\begin{equation*}
   \phi_n(\lambda)\le Ce^{-C\lambda} \le Ce^{-cn^2}.
\end{equation*}
  \end{enumerate}
\end{lem}
The really new results in this Lemma are (iii) and (iv).
Indeed, (i) and (ii) can also be obtained from~\eqref{e-J1.13'}, 
and as indicated, (v) is given in \cite{Johan}.

%%%%=====================================secdep.tex

\bigskip
\section{{\bf De-Poissonization Lemmas}}
\label{s- dep}

In this section, we present two Lemmas which show that $\phi_n(N)$ is 
a good approximation of $q_{n,N}=f_{N,n}/N!$.

We need a Lemma showing the monotonicity of $q_{n,N}$ in $N$.
The statement and proof can be found in \cite{Johan}.
\begin{lem}\label{lem-dep1}
   For all $n,N\geq 1$, 
\begin{equation*}
  q_{n,N+1} \le q_{n,N}.
\end{equation*}
\end{lem}

Using this monotonicity result, 
the following Tauberian-like ``de-Poissonization'' Lemma can be proved.
This is a modification of Lemma 2.5 in \cite{Johan} and 
the proof is the same.
\begin{lem}\label{lem-dep2}
   Let $m>0$ be a fixed real number. 
   Set $\mu_N^{(m)}=N+(2\sqrt{m+1}+1)\sqrt{N\log N}$
and $\nu_N^{(m)}=N-(2\sqrt{m+1}+1)\sqrt{N\log N}$.
Then there are constants $C=C(m)$ and $N_0=N_0(m)$ such that 
\begin{equation*}
   \phi_n(\mu_N^{(m)})-\frac{C}{N^m} \le q_{n,N}
\le \phi_n(\nu_N^{(m)})+\frac{C}{N^m}
\end{equation*}
for $N\geq N_0$, $0\le n\le N$.
\end{lem}

The reader will observe that the above Lemma is actually 
enough for all of our future calculations.  
Nevertheless, for convenience and the purpose of illustration, 
we use the following Lemma for the convergence of moments.
\begin{lem}\label{lem-dep3}
   There exists $C>0$ such that 
\begin{equation*}
    q_{n,N} \le C\phi_n(N-\sqrt{N}),\quad
1-q_{n,N} \le C\bigl(1-\phi_n(N+\sqrt{N})\bigr)
\end{equation*}
for all sufficient large $N$, $0\le n\le N$.
\end{lem}
\begin{proof}
   Note that $q_{n,N} \geq 0$. Using Lemma~\ref{lem-dep1} and 
Stirling's formula for sufficiently large $N$, we have from~\eqref{e-int4}, 
\begin{equation*}
  \begin{split}
   &\phi_n(N-\sqrt{N})
= \sum_{N'=0}^{\infty} \frac{e^{-(N-\sqrt{N})}(N-\sqrt{N})^{N'}}{(N')!} q_{n,N'}\\
   &\qquad \geq \sum_{N'\geq N-\sqrt{N}}^{N} 
\frac{e^{-(N-\sqrt{N})}(N-\sqrt{N})^{N'}}{(N')!} q_{n,N'}\\
   &\qquad \geq q_{n,N}\sum_{N'\geq N-\sqrt{N}}^{N} 
\frac{e^{-(N-\sqrt{N})}(N-\sqrt{N})^{N'}}{(N')!}\\
   &\qquad \geq C q_{n,N} \sum_{N'\geq N-\sqrt{N}}^{N} 
\frac{e^{-(N-\sqrt{N})}(N-\sqrt{N})^{N'}}{(N')^{N'+1/2}e^{-N'}}
   = C q_{n,N} \sum_{N'\geq N-\sqrt{N}}^{N} e^{f(N')}, 
  \end{split}
\end{equation*}
where $f(x)=-(N-\sqrt{N})+x\log (N-\sqrt{N}) +x-(x+\frac12)\log x$.
One can easily check that $f(x)$ is a decreasing function 
for $x\geq (N-\sqrt{N})$. Thus 
\begin{equation*}
   \phi_n(N-\sqrt{N})\geq Cq_{n,N} \sqrt{N} e^{f(N)}
=Cq_{n,N} e^{\sqrt{N}+N\log (1-1/\sqrt{N})} \geq Cq_{n,N},
\end{equation*}
for sufficiently large $N$, $0\le n\le N$.

 For the second inequality, note that $q_{n,N} \le 1$ by definition.
 Again, using Lemma~\ref{lem-dep1} and
Stirling's formula for sufficiently large $N$,
\begin{equation*}
  \begin{split}
   &1-\phi_n(N+\sqrt{N})
= \sum_{N'=0}^{\infty} \frac{e^{-(N+\sqrt{N})}(N+\sqrt{N})^{N'}}{(N')!} 
(1-q_{n,N'})\\
   &\qquad \geq C (1-q_{n,N}) \sum_{N'=N}^{N+\sqrt{N}} e^{g(N')}
  \end{split}
\end{equation*}
where $g(x)=-(N+\sqrt{N})+x\log (N+\sqrt{N}) +x-(x+\frac12)\log x$.
One can check that for $N\le x\le N+\sqrt{N}$, 
$g''(x)<0$ so that $\min g(x)=\min \bigl(g(N), g(N+\sqrt{N})\bigr)$.
If $N$ is sufficiently large, $\min \bigl(g(N), g(N+\sqrt{N})\bigr)
=g(N+\sqrt{N})=-\frac12 \log(N+\sqrt{N})$. 
Therefore 
\begin{equation*}
   1-\phi_n(N+\sqrt{N})\geq C(1-q_{n,N}) \sqrt{N} e^{g(N)}
\geq C(1-q_{n,N}),
\end{equation*}
for sufficiently large $N$, $0\le n\le N$.
\end{proof}

%%%%===================================== secpf.tex

\bigskip
\section{{\bf Proofs of main Theorems}}
\label{s-pf}

In this Section, we prove the main Theorems.\\

\noindent{\it\bf Proof of Theorem 1.1} 
  Assume for definiteness that $t<0$. 
For $t\geq 0$, the calculation is similar.
  From the definition of $q_{n,N}\equiv \frac{f_{N,n}}{N!}$, 
\begin{equation}\label{e-ind2}
  \begin{split}
F_N(t)=Prob (\frac{l_N-2\sqrt{N}}{N^{1/6}}\le t)
    =q_{[2\sqrt{N}+tN^{1/6}],N}.
  \end{split}
\end{equation}
Set 
\begin{equation*}
   n=[2\sqrt{N}+tN^{1/6}].
\end{equation*}
As $t$ is fixed, observe that $0\le n\le N$, as $N\to\infty$.
Using Lemma~\ref{lem-dep2} with any fixed value of $m>0$, we have
\begin{equation*}
   \phi_n(\mu_N^{(m)})-\frac{C}{N^m} \le F_N(t)
\le \phi_n(\nu_N^{(m)})+\frac{C}{N^m}.
\end{equation*}
Set 
\begin{equation*}
   t_N=2^{1/3}(n+1)^{2/3}\bigl(1-\frac{2\sqrt{\mu_N^{(m)}}}{n+1}\bigr).
\end{equation*}
(cf. the definition of $t$ in~\eqref{e-z42}.)
Then, for all large $N$, 
\begin{equation*}
   2t\le t_N \le\frac12 t, \quad\text{and}\quad \lim_{N\to\infty}t_N=t.
\end{equation*}
Let $M_6\geq 2|t|$ be any sufficiently large, fixed number satisfying 
condition $(c)$ in Lemma~\ref{lem-phi}. 
Using Lemma~\ref{lem-phi} (iii), we have, for some constant 
$C(M_6)$ which depends on $M_6$, and a constant $C$ which is 
independent of $M_6$, 
\begin{equation*}
   \phi_n(\mu_N^{(m)})=\exp\biggl(\int_{t_N}^{\infty} 
2im^{PII}_{1,22}(y) dy\biggr)
\bigl(1+O_{M_6}(\frac{1}{n^{1/3}})+O(e^{-\frac14M_6^{3/2}})\bigr)
\end{equation*}
Taking $N\to\infty$, and then taking $M_6\to\infty$, we obtain, 
\begin{equation*}
   \lim_{N\to\infty} \phi_n(\mu_N^{(m)})=
\exp\biggl(\int_{t}^{\infty} 2im^{PII}_{1,22}(y) dy\biggr).
\end{equation*}
For $\phi_n(\nu_N^{(m)})$, we obtain the same limit by a similar calculation,
\begin{equation*}
   \lim_{N\to\infty} \phi_n(\nu_N^{(m)})=
\exp\biggl(\int_{t}^{\infty} 2im^{PII}_{1,22}(y) dy\biggr).
\end{equation*}
Thus, recalling $\frac{d}{dx} 2i(m_1^{PII})_{22}(x) = u^2(x)$ in \eqref{e-um1PII}, 
integration by parts yields
\begin{equation*}
   \lim_{N\to\infty} Prob\biggl(\frac{l_N-2\sqrt{n}}{N^{1/6}}\le t\biggr)
=\exp\biggl(\int_{t}^{\infty} 2im^{PII}_{1,22}(y) dy\biggr)
=F(t).
\end{equation*}

\hfill   $\square$ 
\bigskip

\noindent{\it\bf Proof of Theorem 1.2}
   Integrating by parts,
\begin{equation*}\label{e-mon2}
   \E_N(\chi_N^m)=\int_{-\infty}^{\infty}
t^mdF_N(t)
=-\int_{-\infty}^{0} mt^{m-1}F_N(t) dt
+\int_{0}^{\infty} mt^{m-1} (1-F_N(t)) dt
\end{equation*}
where $F_N(t)\equiv Prob\biggl(\frac{l_N-2\sqrt{N}}{N^{1/6}}\le t\biggr)$
as in Theorem~\ref{thm1}.
From Theorem~\ref{thm1}, we have pointwise convergence of $F_N(t)$ to $F(t)$.
We need uniform control of $F_N$ for large $N$.
Let $M>0$ be a sufficiently large, fixed number
and $0< \delta <\frac14$ be a fixed, sufficiently small number.

Set $n=[2\sqrt{N}+tN^{1/6}]$.
First consider the case when $t\le -M$.
If $t<-2N^{1/3}$, then $F_N(t)=Prob (l_N\le 2\sqrt{N}+tN^{1/6})\le Prob (l_N<0)=
0$.
For $-2N^{1/3}\le t\le -M$,
~\eqref{e-ind2} and Lemma~\ref{lem-dep3} yield
\begin{equation}\label{e-mon3}
   F_N(t)=q_{n,N}\le C\phi_n(N-\sqrt{N}).
\end{equation}
If $-2N^{1/3}\le t\le -2\delta N^{1/3}$, when $N$ is sufficiently large,
\begin{equation*}
   \frac{2\sqrt{N-\sqrt{N}}}{n+1} \geq
\frac{2\sqrt{N}(1-\frac1{\sqrt{N}})^{1/2}}{2\sqrt{N}+tN^{1/6}+1}
\geq \frac{2\sqrt{N}(1-\frac{\delta}4)}{2(1-\delta)\sqrt{N}+1}
\geq 1+\frac{\delta}2.
\end{equation*}
Thus, using Lemma~\ref{lem-phi} (v), for large $N$, 
\begin{equation}\label{e-z43}
   \phi_n(N-\sqrt{N})\le Ce^{-cN}\le Ce^{ct^3}.
\end{equation}
If $-2\delta N^{1/3}\le t\le -M$,
\begin{equation*}\label{e-mon4}
   \frac{2\sqrt{N-\sqrt{N}}}{n+1} \le
\frac{2\sqrt{N}}{2\sqrt{N}-2\delta\sqrt{N}}\le 1+2\delta
\end{equation*}
and, using the monotonicity of $(2\sqrt{\lambda}-x)/x^{1/3}$
as a function of $x\le 2\sqrt{\lambda}$,
\begin{equation*}\label{e-mon5}
\begin{split}
    & 2^{1/3}\biggl(\frac{2\sqrt{N-\sqrt{N}}-(n+1)}{(n+1)^{1/3}}\biggr)\\
    &\qquad\geq 2^{1/3}\biggl(\frac{2\sqrt{N}
(1-\frac1{\sqrt{N}})^{1/2}-(2\sqrt{N}-MN^{1/6}+1)}
{(2\sqrt{N}-MN^{1/6}+1)^{1/3}}\biggr) \geq \frac{M}2
\end{split}
\end{equation*}
as $N\to\infty$.
Thus, we have for $-2\delta N^{1/3}\le t\le -M$,
\begin{equation*}\label{e-mon6}
   1+\frac{\frac12M}{2^{1/3}{(n+1)}^{2/3}} \le \frac{2\sqrt{N-\sqrt{N}}}{n+1}
\le 1+2\delta
\end{equation*}
Therefore, from Lemma~\ref{lem-phi} (iv), provided $\frac{M}2$ satisfies 
condition $(d)$ and $2\delta$ satisfies condition $(a)$, 
\begin{equation}\label{e-mon7}
   \phi_n(N-\sqrt{N})\le C_M
\exp\biggl(-\frac1{96}\bigl(2^{1/3}(n+1)^{2/3}
(\frac{2\sqrt{N-\sqrt{N}}}{n+1}-1)\bigr)^3 \biggr).
\end{equation}
On the other hand,
using the monotonicity of $(2\sqrt{\lambda}-x)/x^{1/3}$
as a function of $x\le 2\sqrt{\lambda}$,
\begin{equation*}\label{e-mon8}
 \begin{split}
    &2^{1/3}\biggl(\frac{2\sqrt{N-\sqrt{N}}-(n+1)}{(n+1)^{1/3}}\biggr)\\
    &\qquad \geq 2^{1/3}\biggl(
\frac{2\sqrt{N}(1-\frac1{\sqrt{N}})^{1/2}-(2\sqrt{N}+tN^{1/6}+1)}
{(2\sqrt{N}+tN^{1/6}+1)^{1/3}}\biggr)
\geq -\frac{t}2
 \end{split}
\end{equation*}
for all $-2\delta N^{1/3}\le t\le -M$, as $N\to\infty$.
Therefore ~\eqref{e-mon7} gives us
\begin{equation}\label{e-z44}
    \phi_n(N-\sqrt{N})\le C(M) \exp(\frac1{768}t^3)
\end{equation}
for $-2\delta N^{1/3}\le t\le -M$.

Inserting the above estimates~\eqref{e-z43} and~\eqref{e-z44} 
into~\eqref{e-mon3}, 
we obtain 
\begin{equation}\label{e-z45}
    F_N(t)\le C(M)e^{ct^3}\qquad\text{for $-2N^{1/3}\le t\le -M$}, 
\end{equation}
and as $F_N(t)=0$ for $t< -2N^{1/3}$, it follows by the 
dominated convergence theorem that 
\begin{equation}
   \lim_{N\to\infty} \int_{-\infty}^{0} mt^{m-1}F_N(t) dt
=\int_{-\infty}^{0} mt^{m-1}F(t) dt.
\end{equation}

Now consider the case when $t\geq M$.
If $t> N^{5/6}-2N^{1/3}$, then $1-F_N(t)=
1-Prob (l_N\le 2\sqrt{N}+tN^{1/6}) \le 1-Prob(l_N> N)=0$.
For $M\le t\le N^{5/6}-2N^{1/3}$,
again, \eqref{e-ind2} and Lemma~\ref{lem-dep3} yield
\begin{equation}\label{e-z51}
   1-F_N(t)=1-q_{n,N}\le C(1-\phi_n(N+\sqrt{N})).
\end{equation}
If $2\delta N^{1/3}\le t\le N^{5/6}-2N^{1/3}$, when $N$ is sufficiently large,
\begin{equation*}
   \frac{2\sqrt{N+\sqrt{N}}}{n+1} \le
\frac{2\sqrt{N}(1+\frac{\delta}4)}{2\sqrt{N}+2\delta\sqrt{N}}
\le 1-\frac{\delta}2.
\end{equation*}
Thus, using Lemma~\ref{lem-phi} (i),
\begin{equation}\label{e-z46}
   1-\phi_n(N+\sqrt{N})\le Ce^{-cn}\le Ce^{-c\sqrt{N}}\le Ce^{-ct^{3/5}}.
\end{equation}
If $M\le t\le 2\delta N^{1/3}$, similar calculations to the case
$-2\delta N^{1/3}\le t\le -M$ yield
\begin{equation*}
   \frac12\le 1-2\delta \le \frac{2\sqrt{N+\sqrt{N}}}{n+1} \le
1-\frac{\frac12M}{2^{1/3}{(n+1)}^{2/3}}.
\end{equation*}
Therefore, from Lemma~\ref{lem-phi} (ii), as $N\to\infty$, 
\begin{equation}\label{e-mon15}
   1-\phi_n(N-\sqrt{N})\le
\exp\biggl(-c\bigl(2^{1/3}(n+1)^{2/3}
(1-\frac{2\sqrt{N+\sqrt{N}}}{n+1})\bigr)^{3/2} \biggr), 
\end{equation}
provided $\frac12 M$ satisfies condition $(d)$.
However, as in the case $-2\delta N^{1/3}\le t\le -M$, we have
\begin{equation*}
    2^{1/3}\biggl(\frac{(n+1)-2\sqrt{N+\sqrt{N}}}{(n+1)^{1/3}}\biggr)
\geq 2^{1/3}\biggl(
\frac{(2\sqrt{N}+tN^{1/6})-2\sqrt{N}(1+\frac1{\sqrt{N}})^{1/2}}
{(2\sqrt{N}+tN^{1/6})^{1/3}}\biggr)
\geq \frac{t}2.
\end{equation*}
for all $M\le t\le 2\delta N^{1/3}$, as $N\to\infty$.
Therefore \eqref{e-mon15} gives us
\begin{equation}\label{e-z47}
   1-\phi_n(N-\sqrt{N})\le C\exp(-ct^{3/2})
\end{equation}
for $M\le t\le 2\delta N^{1/3}$.

Inserting the above estimates~\eqref{e-z46} and~\eqref{e-z47} 
into~\eqref{e-z51}, we obtain for $M\le t\le N^{5/6}-2N^{1/3}$ 
\begin{equation}\label{e-z48}
   1-F_N(t)\le Ce^{-ct^{3/5}}
\end{equation}
as $N\to\infty$.
Once again, as $1-F_N(t)=0$ for $t>N^{5/6}-2N^{1/3}$, it follows by 
the dominated convergence theorem that 
\begin{equation}
   \lim_{N\to\infty} \int_{0}^{\infty} mt^{m-1} (1-F_N(t)) dt
=\int_{0}^{\infty} mt^{m-1} (1-F(t)) dt.
\end{equation}

\hfill  $\square$

%%%%%%%%%%================ref.tex

%%\include{appk}

\appendix 
\section{}

As advertised in the Introduction, in 
this Appendix we give a new derivation of the formula
\begin{equation}\label{a1}
\sum_{N=0}^\infty\frac{\lambda^NF_N(n)}{N!}=
\det(d_{j-k})_{0\le j,k\le n-1},
\end{equation}
where $d_j=(2\pi)^{-1}\int_0^{2\pi}\exp(2\sqrt{\lambda}\cos\theta-ij\theta)
d\theta$, and $F_N(n)$ is the distribution function for the length, 
$\ell_N(\pi)$, of the longest increasing subsequence in the random permutation
$\pi$ from $S_N$. We set $F_0(0)=1$. 
 
Let $\mu=(\mu_1,\mu_2,\dots,\mu_r,0,0,\dots)$, $\mu_1\ge\mu_2\ge\dots$,
be a partition of $N$, i.e. $\mu_j$, $1\le j\le r$, are positive integers and
$N=\mu_1+\dots+\mu_r$; we write $\mu\vdash N$. With $\mu$ we can associate 
a Young diagram, also denoted by $\mu$, in the standard way, see for example
[Sa]. In the Young diagram there are $\mu_j$ boxes in the $j:th$ row. If we
insert the numbers $1,\dots, N$ in the boxes in such a way that the numbers 
in every row and column are increasing we get a (standard) Young tableau
$t\,$; $t$ has shape $\mu$, $s(t)=\mu$. Let $r(\mu)$ denote the number of
rows in $\mu$.
 
Schensted, [Sc], has constructed a certain bijection, the Schensted
correspondence, between the permutation group $S_N$ and pairs of
Young tableaux $(t,t')$ with the same shape $s(t)=s(t')=\mu$, where
$\mu\vdash N$. This correspondence has the property that if $S_N\owns
\pi\to (t,t')$, $\mu=s(t)$, then $\ell_N(\pi)$ equals the length, $\mu_1$, 
of the first row in $\mu$, and the length, $\ell'_N(\pi)$, of the longest
{\it decreasing} subsequence in $\pi$ equals $r(\mu)$, the number of rows 
in $\mu$. For details see [Sa].
 
If we put the uniform probability distribution on $S_N$ then clearly
the random variables $\ell_N$ and $\ell'_N$ have the same distribution
(just reverse the permutation). Let $f(\mu)$ denote the number of Young 
tableaux with shape $\mu$. Then, by the Schensted correspondence,
\begin{equation}\label{a2}
F_N(n)=\frac 1{N!}\sum_{\substack{\mu\vdash N \\ r(\mu)\le n}} f(\mu)^2.
\end{equation}
If we set $h_j=\mu_j+r-j$, $r=r(\mu)$, we have the following formula,
due to Frobenius and Young,
\begin{equation}\label{a3}
f(\mu)=N!\prod_{1\le i<j\le r} (h_i-h_j)\prod_{i=1}^r\frac 1{h_i!},
\end{equation}
see for example [Si]. Note that $N=\sum_{j=1}^r\mu_j=\sum_{j=1}^rh_j-r(r-1)/2$
and $h_{j-1}-h_j=\mu_{j-1}-\mu_j+1\ge 1$. Combining the formulas~\eqref{a2} 
and~\eqref{a3} we get
\begin{equation}\label{a4}
F_N(n)=N!\sum_{r=1}^n\frac 1{r!}
\sum_{(\ast)}\Delta(h)^2\prod_{j=1}^r\frac 1{(h_j!)^2},
\end{equation}
where the $(\ast)$ means that we sum over all different integers $h_i\ge 1$ 
such that $\sum h_j=N+r(r-1)/2$, and $\Delta(h)=\prod_{i<j}(h_j-h_i)$ is the
Vandermonde determinant. That we can remove the ordering of the
$h_j\,$'s in~\eqref{a4} follows from symmetry under
permutation of $h_1,\dots,h_r$. The constraint $\sum h_j=N+r(r-1)/2$ is
removed by the Poissonization
\begin{equation}\label{a5}
\phi_n(\lambda)=e^{-\lambda}\sum_{N=0}^\infty\frac{\lambda^N}{N!}F_N(n)=
e^{-\lambda}[1+\sum_{r=1}^n\lambda^{-r(r-1)/2}H_r(\lambda)],
\end{equation}
where
$$
H_r(\lambda)=\frac 1{r!}\sum_{h\in \mathbb Z_+^r}\Delta(h)^2\prod_{j=1}^r
\frac{\lambda^{h_j}}{(h_j!)^2}.
$$
We have used the fact that $\sum h_j\ge 1+\dots +r=r(r-1)/2+r$ and
$N\ge r$, since the $h_j\,$'s are different integers. The condition that
the $h_j\,$'s are different can then be removed since otherwise $\Delta(h)=0$.
Observe that $H_r(\lambda)$ is a Hankel determinant with respect to the
discrete measure
$$
\nu(\{m\})=\frac{\lambda^m}{(m!)^2},\quad m\in\mathbb Z_+,
$$
see [Sz1], i.e.
$$
H_r(\lambda)=\det(\sum_{m=1}^\infty m^{j+k}\frac{\lambda^m}{(m!)^2})_
{0\le j,k\le r-1}.
$$
If $q_j$, $j\ge 0$, are any polynomials with $\deg q_j=j$ and leading 
coefficient $1$, row and column operations on the determinant gives
\begin{equation}\label{a6}
H_r(\lambda)=\det(\sum_{m=1}^\infty q_j(m)q_k(m)\frac{\lambda^m}{(m!)^2})_
{0\le j,k\le r-1}.
\end{equation}
We now make a particular choice of $q_j$, $q_j(x)=x(x-1)\dots (x-(j-1))$,
if $j\ge 1$ and $q_0(x)=1$, so that
\begin{equation}\label{a7}
a^j\frac {d^j}{da^j} a^m=q_j(m)a^m,\quad m,j\ge 0.
\end{equation}
The elements in the Hankel determinants can then be written
\begin{equation}\label{a8}
\left.\sum_{m=1}^\infty q_j(m)q_k(m)\frac{\lambda^m}{(m!)^2}=
a^jb^k\frac{d^j}{da^j}\frac{d^k}{db^k}\sum_{m=0}^\infty
\frac{a^mb^m}{(m!)^2}\right|_{a=b=\sqrt{\lambda}}-\delta_{j0}\delta_{k0}.
\end{equation}
Now,
$$
\sum_{m=0}^\infty\frac{a^mb^m}{(m!)^2}=\frac 1{2\pi}\int_0^{2\pi}
e^{ae^{i\theta}+be^{-i\theta}}d\theta
$$
and hence we can perform the differentiations in~\eqref{a8} and get
$$
\sum_{m=1}^\infty q_j(m)q_k(m)\frac{\lambda^m}{(m!)^2}=\lambda^{(j+k)/2}
d_{j-k}-\delta_{j0}\delta_{k0},
$$
where $d_{j-k}=(2\pi)^{-1}\int_0^{2\pi}\exp(2\sqrt{\lambda}\cos\theta-
i(j-k)\theta)d\theta$. Inserting this identity into 
the formula~\eqref{a6} yields
\begin{equation}\label{a9}
H_r(\lambda)=\lambda^{r(r-1)/2}(D_r-D_{r-1}), \quad r\ge 1,
\end{equation}
where $D_r$ is the Toeplitz determinant $\det(d_{j-k})_{0\le j,k\le r-1}$
and $D_0=1$. Hence, using the formula~\eqref{a5}, we get $\phi_n(\lambda)=
e^{-\lambda}D_n$, which is what we wanted to prove.
 
In the remaining part of this appendix we will give a heuristic argument 
showing why we can expect the random variable $\ell_N(\pi)$ to behave
like the largest eigenvalue of a random hermitian matrix. From our
considerations above we see that
$$
F_N(n)=\frac 1{N!}\sum_{\substack{\mu\vdash N \\ \mu_1\le n}} f(\mu)^2.
$$
By the same computations as above this leads to
\begin{equation}\label{a10}
\phi_n(\lambda)=e^{-\lambda}[1+\sum_{r=1}^\infty\lambda^{-r(r-1)/2}
H_r(\lambda;n)],
\end{equation}
where
$$
H_r(\lambda;n)=\frac 1{r!}\sum_{h\in\{1,\dots,n+r-1\}^r}\Delta(h)^2\prod_{j=1}^r
\frac{\lambda^{h_j}}{(h_j!)^2}.
$$
Note that $H_r(\lambda;n)\nearrow H_r(\lambda)$ as $n\to\infty$. We can think
of
\begin{equation}\label{a11}
\frac 1{r!H_r(\lambda)}\Delta(h)^2\prod_{j=1}^r
\frac{\lambda^{h_j}}{(h_j!)^2}=
\frac 1{r!H_r(\lambda)}e^{-2\sum_{i<j}\log|h_i-h_j|^{-1}+\sum_j
[(\log\lambda)h_j+2\log(h_j!)]}
\end{equation}
as the probability of the configuration 
$h\in\mathbb Z_+^r$. This probability has
the form of a {\it discrete} Coulomb gas 
on $\mathbb Z_+$ at inverse temperature 
$\beta=2$ confined by an external potential. An $N\times N$ random hermitian 
matrix with a probability density of the form 
$Z_N^{-1}\exp(-\text{Tr\,}V(M))$ has an eigenvalue density
$$
\frac 1{Z_N}e^{-2\sum_{i<j}\log|x_i-x_j|^{-1}+\sum_j V(x_j)},
$$
with $x\in\mathbb R^N$; $x_1,\dots ,x_N$ are the eigenvalues of $M$.
Thus we can think of the $h_j\,$;s as some kind of ``eigenvalues''.
 
Let
$$
P_r(\lambda;n)=H_r(\lambda;n)/H_r(\lambda),
$$
i.e. $P_r(\lambda;n)$ is the probability that the largest ``eigenvalue''
is $\le n+r-1$. Then, by~\eqref{a9} and~\eqref{a10},
\begin{equation}\label{a12}
\phi_n(\lambda)=e^{-\lambda}[1+\sum_{r=1}^\infty P_r(\lambda;n)
(D_r-D_{r-1})]=e^{-\lambda}+\sum_{r=1}^\infty P_r(\lambda;n)
(\phi_r(\lambda)-\phi_{r-1}(\lambda)).
\end{equation}
Now, the essential contribution to the right-hand side of~\eqref{a12} 
comes from
$r$ around $2\sqrt{\lambda}$ since otherwise 
$\phi_r(\lambda)-\phi_{r-1}(\lambda)$ is very small. Thus
$$
\phi_n(\lambda)\approx P_{2\sqrt{\lambda}}(\lambda;n),
$$
i.e. $\phi_n(\lambda)$ is like the probability that the largest
``eigenvalue'' in the discrete Coulomb gas~\eqref{a11} 
is $\le n+2\sqrt{\lambda}$.
\newline\newline


\begin{thebibliography}{10}

\bibitem[AS]{AS}
~M.Abramowitz and ~I.A.Stegun, \emph{Handbook of Mathematical Functions}, 
Dover Publications, New York, (1965). 

\bibitem[AD]{AD}
~D.Aldous and ~P.Diaconis, \emph{Hammersley's Interacting Particle
Process and Longest Increasing Subsequences},
Prob. Th. and Rel. Fields, \textbf{103}, 199-213, (1995).

\bibitem[Bo]{Bo}
~A.Borodin, \emph{Longest increasing subsequences of random 
colored permutations}, Electron. J. Combin., \textbf{6 (1)}, 
R13, (1999).

\bibitem[BB]{BB}
~R.M.Baer and ~P.Brock, \emph{Natural sorting over permutation spaces}, 
Math. Comp., \textbf{22}, 385-410, (1968).

\bibitem[BC]{BC}
~R.Beals and ~R.Coifman, \emph{Scattering and inverse scattering for 
first order systems}, Comm. Pure Appl. Math., \textbf{37}, 
39-90, (1984).

\bibitem[BDJ]{BDJ}
~J.Baik, ~P.Deift and ~K.Johansson, 
\emph{On the distribution of the length of the second row 
of a Young diagram under Plancherel measure}, 
preprint, LANL E-print math.CO/9901118.

\bibitem[BR]{BR}
~J.Baik and ~E.Rains, \emph{Symmetrized increasing subsequence problems}, 
in preparation.
\bibitem[CG]{CG}
~K.Clancey and ~I.Gohberg, \emph{Factorization of Matrix Functions 
and Singular Integral Operators}, Birkh\"auser, (1981). 

\bibitem[De]{De}
~P.A.Deift, \emph{Integrable Hamiltonian systems. Dynamical systems 
and probabilistic methods in partial differential equations}, 103-138, 
in \emph{Lectures in Applied Mathematics, 31}, edited by 
~P.A.Deift,~C.D.Levermore and~C.E.Wayne, American 
Mathematical Society, Providence, RI, 1996.

\bibitem[DeZe1]{DeZe}
~J.-D.Deuschel and ~O.Zeitouni, \emph{Limiting curves for i.i.d. 
records}, Ann. Probab., \textbf{23}, 852-878, (1995).

\bibitem[DeZe2]{DeZe2}
~J.-D.Deuschel and ~O.Zeitouni, \emph{On increasing subsequences of 
i.i.d. samples}, preprint, (1997).

\bibitem[DIZ]{DIZ}
~P.A.Deift, ~A.R.Its and ~X.Zhou, \emph{Long-time Asymptotics 
for Integrable Nonlinear Wave Equations}, 
in \emph{Important Development in Soliton Theory}, 
\textbf{2nd Edition}, edited by ~A.S.Fokas and ~V.E.Zakharov, 
Springer-Verlag, to be published.

\bibitem[DKM]{DKM}
~P.A.Deift, ~T.Kriecherbauer and ~K.T-R McLaughlin, 
\emph{New Results on the Equilibrium Measure for Logarithmic potentials 
in the Presence of an External Field}, J. Approx. Theory, 
\textbf{95}, no.3, 388-475, (1998).

\bibitem[DKMVZ1]{DKMVZ}
~P.A.Deift, ~T.Kriecherbauer, ~K.T-R McLaughlin, ~S.Venakides 
and ~X.Zhou, \emph{Asymptotics for Polynomials Orthogonal with 
respect to Varying Exponential Weights}, Internat. Math. Res. 
Notices, no.~\textbf{16}, 759-782, (1997). 

\bibitem[DKMVZ2]{DKMVZ2}
~P.A.Deift, ~T.Kriecherbauer, ~K.T-R McLaughlin, ~S.Venakides
and ~X.Zhou, \emph{Strong Asymptotics for Orthogonal Polynomials with 
respect to Varying Exponential Weights via Riemann-Hilbert 
Techniques}, To appear in Comm. Pure. Appl. Math.

\bibitem[DKMVZ3]{DKMVZ3}
~P.A.Deift, ~T.Kriecherbauer, ~K.T-R McLaughlin, ~S.Venakides
and ~X.Zhou, \emph{Uniform Asymptotics for Polynomials Orthogonal with 
respect to Varying Exponential Weights and Applications to 
Universality Questions in Random Matrix Theory}, 
To appear in Comm. Pure. Appl. Math.

\bibitem[DS]{DS}
~P.Diaconis and ~M.Shahshahani, \emph{On the Eigenvalues of
Random matrices}, J. Appl. Prob. \textbf{31}, 49-61, (1994).

\bibitem[DVZ1]{DVZ2}
~P.A.Deift, ~S.Venakides and ~X.Zhou, \emph{The collisionless shock 
region for the long-time behavior of solutions of the
KdV equation}, Comm. Pure Appl. Math., \textbf{47} no. 2, 199--206, (1994).

\bibitem[DVZ2]{DVZ}
~P.A.Deift, ~S.Venakides and ~X.Zhou, \emph{New Results in Small 
Dispersion KdV by an Extension of the Steepest Descent Method for 
Riemann-Hilbert Problems}, Internat. Math. Res.
Notices, no.~\textbf{6}, 285-299, (1997).

\bibitem[DZ1]{DZ1}
~P.A.Deift and ~X.Zhou, \emph{A Steepest Descent Method for 
Oscillatory Riemman-Hilbert Problems; Asymptotics for 
the MKdV Equation}, Ann. Math., \textbf{137}, 295-368, (1993).

\bibitem[DZ2]{DZ}
~P.A.Deift and ~X.Zhou, \emph{Asymptotics for the Painlev\'e II  
Equation}, Comm. Pure. Appl. Math., \textbf{48}, 277-337, (1995).

\bibitem[ES]{ES}
~P.Erd\"os and ~G.Szekeres, \emph{A combinatorial theorem in geometry}, 
Compositio Math., \textbf{2}, 463-470, (1935).
 
\bibitem[FIK]{FIK}
~A.S.Fokas, ~A.R.Its and~V.E.Kitaev, \emph{Discrete Painlev\'e 
equations and their appearance in quantum gravity}, 
Comm. Math. Phy., \textbf{142}, 313-344, (1991).

\bibitem[FMZ]{FMZ}
~A.S.Fokas, ~U.Mugan and ~X.Zhou, \emph{On the Solvability of Painlev\'e I, III
and V}, Inverse Problems, \textbf{8}, 757-785, (1992).

\bibitem[FN]{FN}
~H.Flaschka and ~A.Newell, \emph{Monodromy and spectrum preserving 
deformations, I}, Comm. Math. Phy., \textbf{76}, no.1, 67-116, (1980).

\bibitem[FZ]{FZ}
~A.S.Fokas and ~X.Zhou, \emph{On the Solvability of Painlev\'e II and IV}, 
Comm. Math. Phy., \textbf{144}, 601-622, (1992).

\bibitem[Ge]{Ge}
~I.M.Gessel, \emph{Symmetric functions and P-recursiveness}, 
J. Combin. Theory. Ser. A, {\bf 53}, 257 - 285, (1990).

\bibitem[GK]{GK}
~I.Gohberg and ~N.Krupnik, \emph{One-Dimensional Linear Singular 
Integral Equations vol.I and II}, Operator theory, advances 
and applications ;  v. 53-54, Birkh\"auser Verlag, Basel, 1992

\bibitem[GW]{GW}
~D.J.Gross and ~E.Witten, \emph{Possible third-order phase transition 
in the large N lattice gauge theory}, Phys. Rew. D, \textbf{21}, 
446-453, (1980).

\bibitem[GWW]{GWW}
~I.M.Gessel, ~J.Weinstein and ~H.S.Wilf, \emph{Lattice walks in
$\mathbb Z^d$ and permutations with no long ascending subsequences},
Electr. J. Combin., {\bf 5(1)}, (1998).

\bibitem[Ha]{Ha}
~J.M.Hammersley, \emph{A few seedlings of research}, Proc. Sixth
Berkeley Symp. Math. Statist. and Probability, Vol. 1,
345-394, University of California Press, 1972.

\bibitem[Hi]{Hi}
~M.Hisakado, \emph{Unitary matrix models and Painlev\'e III}, 
Modern Phys. Lett. A, \textbf{11}, no.38, 3001-3010, (1996).

\bibitem[HM]{HM}
~S.P.Hastings and ~J.B.McLeod, \emph{A boundary value problem
associated with the second Painlev\'e transcendent and
the Korteweg de Vries equation}, Arch. Rational Mech. Anal.
\textbf{73}, 31-51, (1980).

\bibitem[Ka]{Ka}
~S.Kamvissis, \emph{On the Long Time Behavior of the Doubly Infinite 
Toda Lattice under Initial Data Decaying at Infinity}, 
Comm. Math. Phy., \textbf{153}, 479-519, (1993).

\bibitem[Kn]{Kn}
~D.E.Knuth, \emph{The art of computer programming},
vol. 3 : \emph{sorting and searching},
2nd ed., Addison Wesley, Reading, Mass., 1973.

\bibitem[IN]{IN}
~A.R.Its and ~V.Yu.Novokshenov, \emph{The Isomonodromic Deformation 
Method in the Theory of Painlev\'e Equations}, Lecture Notes 
in Math.\textbf{1191}, Springer-Verlag, Berlin, Heidelberg, 1986.

\bibitem[Jo1]{Johan}
~K.Johansson, \emph{The Longest Increasing Subsequence in 
a Random Permutation and a Unitary Random Matrix Model}, 
Math. Res. Lett., \textbf{5}, no.1-2, 63-82, (1998).

\bibitem[Jo2]{Jo2}
~K.Johansson, \emph{Shape fluctuations and random matrices}, 
LANL E-print math.CO/9903134.

\bibitem[Jo3]{Jo3}
~K.Johansson, \emph{Transversal fluctuations for increasing subsequences 
on the plane}, preprint, 1999.

\bibitem[JMU]{JMU}
~M.Jimbo, ~T.Miwa and ~K.Ueno, \emph{Monodromy preserving deformations 
of linear ordinary differential equations with rational coefficients, I. 
General theory and $\tau$-function}, 
Physica D, \textbf{2}, no.2, 306-352, (1981).

\bibitem[Ki]{Ki}
~J.-H. Kim, \emph{On the longest increasing subsequence of random
permutations - a concentration result}, J. Comb. Th. A, vol. \textbf{76},
148-155, (1996).

\bibitem[LS]{LS}
~B.F.Logan and ~L.A.Shepp, \emph{A variational problem
for random Young tableaux}, Advances in Math., \textbf{26},
206-222, (1977).

\bibitem[Me]{Me}
~M.L.Mehta, \emph{Random Matrices}, Second Edition, Academic Press, 
San Diago, 1991.

\bibitem[Ok]{Ok}
~A.Okounkov, \emph{Random matrices and random permutations},
preprint, 1999.

\bibitem[OPWW]{OPWW}
~A.M.Odlyzko, ~B.Poonen, ~H.Widom and ~H.S.Wilf, \emph{On the
distribution of longest increasing subsequences in random permutations},
unpublished manuscript.

\bibitem[OR]{OR}
~A.M.Odlyzko and ~E.M.Rains, \emph{On longest increasing subsequences 
in random permutations}, in preparation.

\bibitem[PS]{PS}
~V.Periwal and~D.Shevitz, \emph{Unitary-Matrix Models as Exactly Solvable 
String Theories}, Phys. Rev. Lett., {\bf 64}, 1326-1329, (1990).

\bibitem[Ra]{Rains}
~E.M.Rains, \emph{Increasing subsequences and the classical groups},
Electron. J. of Combinatorics, {\bf 5(1)}, R12, (1998).

\bibitem[Sa]{Sa}
~B.Sagan, \emph{The Symmetric Group : Representations, Combinatorial 
Algorithms, and Symmetric Functions}, Wadsworth$\&$Books/Cole, 
Pacific Grove, Calif., 1991.

\bibitem[Sc]{Sc}
~C.Schensted, \emph{Longest increasing and decreasing subsequences},
Canad. J. Math., {\bf 13}, 179 - 191, (1961).

\bibitem[Se1]{Se}
~T.Sepp\"al\"ainen, \emph{A microscopic model for the
Burgers equation and longest increasing subsequences},
Electron. J. Prob., \textbf{1}, no.5, (1996)

\bibitem[Se2]{Se2}
~T.Sepp\"al\"ainen, \emph{Large deviations for increasing sequences 
on the plane}, Probab. Theory Related Fields, \textbf{112}, 
no.2, 221-244, (1998).

\bibitem[Si]{Si}
~B.Simon, \emph{Representations of Finite and Compact Groups},
Graduate Studies in Mathematics vol. 10, American Mathematical Society,
1996.

\bibitem[ST]{ST}
~E.B.Saff and V.Totik, \emph{Logarithmic Potentials with External 
Fields}, Springer-Verlag, New York, 1997.

\bibitem[Sz1]{Sz}
~G.Szeg\"o, \emph{Orthogonal Polynomials}, American Mathematical 
Society Colloquium Publications, Vol 23, 4th Ed, New York, 1975.

\bibitem[Sz2]{Sz2}
~G.Szeg\"o, \emph{On Certain Hermitian Forms Associated with the 
Fourier Series of a Positive Function}, Comm. Seminaire Math de 
l'Univ. de Lund, tome supplementaire, dedie a Marcel Riesz, 
228-237, (1952) (or \emph{Gabor Szego : Collected Papers - Vol 3 
(1945-1972)}, 270-280, Birkh\"auser, 1982).

\bibitem[TW1]{TW}
~C.A.Tracy and ~H.Widom, \emph{Level-Spacing distributions
and the Airy kernel}, Comm. Math. Phys., \textbf{159}, 151-174, (1994).

\bibitem[TW2]{TW2}
~C.A.Tracy and ~H.Widom, \emph{Random unitary matrices, permutations 
and Painlev\'e}, preprint, LANL E-print math.CO/9811154.

\bibitem[Ul]{Ul}
~S.M.Ulam, \emph{Monte Carlo calculations in problems of 
mathematical physics}, in \emph{Modern Mathematics for the Engineers}, 
~E.F.Beckenbach, ed., McGraw-Hill, 261-281, 1961.

\bibitem[VK1]{VK}
~A.M.Vershik and ~S.V.Kerov, \emph{Asymptotics of the Plancherel
measure of the symmetric group and the limiting form of Young tables},
Soviet Math. Dokl., \textbf{18}, 527-531, (1977).

\bibitem[VK2]{VK2}
~A.M.Vershik and ~S.V.Kerov, \emph{Asymptotic behavior of the maximum 
and generic dimensions of irreducible representations of the symmetric 
group},
Functional Anal. Appl., \textbf{19}, no.1, 21-31, (1985).

\bibitem[Wi]{Wi}
~H.Widom, personal communication.

\end{thebibliography}
\end{document}